\documentclass[twoside]{article} 

\usepackage{jwn-article}
\usepackage{fullpage}
\usepackage{amsfonts}
\usepackage{turnstile}
\usepackage{mathptmx}
\usepackage{hyperref}
\usepackage{textcomp}
\usepackage{graphicx}
\usepackage{moreverb}
\usepackage{alltt}
\usepackage[all]{xy}
\usepackage{array}

\let\oldmarginpar\marginpar
\renewcommand{\marginpar}[1]{\oldmarginpar[\scriptsize\raggedleft\sf #1]{\scriptsize\raggedright\sf #1}}

\usepackage{amssymb}
\renewcommand{\ge}{\geqslant}
\renewcommand{\le}{\leqslant}

\newcommand{\cmd}[1]{\mathrm{#1}}
\newcommand{\prob}[2][]{\ensuremath{\cmd{Pr}_{#1}\left(#2\right)}}

\newcommand{\condprob}[3][]{\prob[#1]{#2\,|\,#3}}

\newcommand{\state}[1]{\ensuremath{\mathsf{#1}}}

\newcommand{\tval}{\state{T}}
\newcommand{\fval}{\state{F}}

\renewcommand{\epsilon}{\varepsilon}
\renewcommand{\phi}{\varphi}
\renewcommand{\Re}{\mathbb{R}}

\usepackage{stmaryrd}

\newcommand{\embed}[1]{\ensuremath{\left\langle#1\right\rangle}}

\newcommand{\rname}[1]{\textsc{\MakeLowercase{#1}}}

\newcommand{\bpoly}[1]{\ensuremath{\mathfrak{B}_{1}(#1)}}

\newcommand{\cspc}{\:}
\newcommand{\scond}[3][]{\ensuremath{#2 \cspc \leadsto_{#1} \cspc #3}}
\newcommand{\mcond}[3][]{\ensuremath{#2 \cspc \hookrightarrow_{#1} \cspc #3}}
\newcommand{\econd}[3][]{\ensuremath{#2 \cspc \looparrowright_{#1} \cspc #3}}
\newcommand{\fcond}[3][]{\ensuremath{#2 \cspc \vdash^*_{#1} \cspc #3}}
\newcommand{\qfcond}[4][]
           {\ensuremath{#2 \cspc : \cspc #3 \cspc \vdash^{*}_{#1} \cspc #4}}
\newcommand{\tcond}[3][]{\ensuremath{#2 \cspc \rightarrow_{#1} \cspc #3}}
\newcommand{\bcond}[3][]{\ensuremath{#2 \cspc \vdash_{#1} \cspc #3}}

\newcommand{\thd}[1]{\bfseries{#1}} 

\title{Probability Distinguishes Different Types of Conditional
  Statements}

\author{Joseph W. Norman \\ University of Michigan \\
\texttt{jwnorman@umich.edu}}
\date{September, 2014}

\begin{document}

\maketitle

\section{Introduction}

Conditional statements, including subjunctive and counterfactual
conditionals, are the source of many enduring challenges in formal
reasoning.  The language of probability can distinguish among several
different kinds of conditionals, thereby strengthening our methods of
analysis.  Here we shall use probability to define four principal
types of conditional statements: \rname{subjunctive},
\rname{material}, \rname{existential}, and \rname{feasibility}.  Each
probabilistic conditional is quantified by a fractional parameter
between zero and one that says whether it is purely affirmative,
purely negative, or intermediate in its sense.  We shall consider also
\rname{truth-functional} conditionals constructed as statements of
material implication from the propositional calculus; these constitute
a fifth principal type.  Finally there is a sixth type called
\rname{Boolean-feasibility} conditionals, which use Boole's
mathematical logic to analyze the sets of possible truth values of
formulas of the propositional calculus.  Each \rname{Truth-functional}
or \rname{Boolean-feasibility} conditional can be affirmative or
negative, with this sense indicated by a binary true/false parameter.
There are other kinds of conditional statements besides the six types
addressed here.  In particular, some conditionals ought to be regarded
as recurrence relations that generate discrete dynamical systems
\cite{norman-orthodoxy}.

Besides its principal type and the value of its (fractional or binary)
sense-parameter, each conditional statement is further characterized
by its content and by its role in analysis.  Two important aspects of
content are factuality and exception handling.  We shall consider the
factuality of a conditional relative to some proposition whose truth
value is known.  If the known proposition is included in the
conditional statement (usually as part of the antecedent clause), then
the conditional is declared `factual' relative to that proposition;
if the negation of the proposition is included then the conditional is
`antifactual' (strongly counterfactual) relative to it; and if the
proposition is omitted then the conditional is `afactual' (weakly
counterfactual) relative to it.  Like their indicative counterparts,
subjunctive conditionals may be factual or counterfactual; and
furthermore any given conditional statement may be correct or
incorrect.  These three distinct properties---mood, factuality, and
correctness---may be correlated with one another.  We shall consider
two mechanisms to address potential exceptions that may confound
conditional relationships: first, allowing the revision of old
conditional statements as new information becomes available; and
second, expressing conditionals in a cautious, defeasible manner in
the first place.  Regarding their use in analysis, there are two basic
roles for conditional statements: a conditional may be asserted as a
constraint itself, or provided as a query whose truth or falsity is to
be evaluated subject to some other set of constraints.

Recognizing the diverse types of conditional statements helps to
clarify several subtle semantic distinctions.  Using the tools of
probability and algebra, each semantically distinct conditional
statement is represented as a syntactically distinct mathematical
expression.  These expressions include symbolic probability
expressions, polynomials with real-number coefficients, sets of real
numbers defined by polynomial constraints, and systems of equations
and inequalities built from such formulas.  Various algorithmic
methods can be used to compute interesting results from conditional
statements represented in mathematical form.  These computational
methods include linear and nonlinear optimization, arithmetic with
symbolic polynomials, and manipulation of relational-database tables.
Among other benefits, these methods of analysis offer paraconsistent
procedures for logical deduction that produce such familiar results as
\emph{modus ponens}, transitivity, disjunction introduction, and
disjunctive syllogism---while avoiding any explosion of consequences
from inconsistent premises.

The proposed method of analysis is applied to several example problems
from Goodman and Adams
\cite{goodman,adams-subjunctive,adams-conditionals}.

\section{Computational Methods: 
  Probability, Optimization, and Boolean Translation}
\label{sec:methods}

To begin, let us review some computational methods that will be useful
for performing analysis.\marginpar{Provenance of methods: De~Moivre,
  Boole, Kolmogorov, Pearl, Sherali, Tuncbilek, and me}

\subsection{A Basic Parametric Probability Model}
\label{sec:model}

For the purpose of representing the various types of conditional
statements, we develop a basic probability model with two true/false
variables $A$ and $B$ symbolizing logical propositions.  We use the
real-valued parameter $x$ to specify the probability that $A$ is true;
$y$ for the conditional probability that $B$ is true given that $A$ is
true; and and $z$ for the conditional probability that $B$ is true
given that $A$ is false.  Truth and falsity are abbreviated $\tval$
and $\fval$.  In order to enforce the laws of probability, each
real-valued parameter $x$, $y$, and $z$ is constrained to lie between
zero and one.  This gives a parametric probability network that is
specified as the following graph, two input probability tables, and
some associated parameter constraints:
\begin{equation}
  \xymatrix{
    *++[o][F]{A} \ar[r] & *++[o][F]{B}
  }
  \qquad
\begin{tabular}[c]{l|l}\hline
\multicolumn{1}{l|}{$A$} & 
\multicolumn{1}{l}{$\prob[0]{{A}}$} \\ \hline\hline
$\state{T}$ & 
$x$ \\ \hline
$\state{F}$ & 
$1 - x$ \\ \hline
\end{tabular}
  \qquad
\begin{tabular}[c]{l|ll}\hline
\multicolumn{3}{l}{$\condprob[0]{{B}}{{A}}$} \\ \hline\hline
\multicolumn{1}{l|}{$A$} & 
\multicolumn{1}{l}{${B=\state{T}}$} & 
\multicolumn{1}{l}{${B=\state{F}}$} \\ \hline\hline
$\state{T}$ & 
$y$ & 
$1 - y$ \\ \hline
$\state{F}$ & 
$z$ & 
$1 - z$ \\ \hline
\end{tabular}
  \qquad
  \begin{array}{l}
    x,y,z \in \Re \\
    0 \le x \le 1 \\ 0 \le y \le 1 \\ 0 \le z \le 1
  \end{array}
  \label{eq:ab-input}
\end{equation}
There are two sets of variables in this parametric probability
network: the primary variables $A$ and $B$, which represent logical
propositions; and the parameters $x$, $y$, and $z$, which are used to
specify probabilities associated with these propositions.  For this
model, each input probability is a polynomial in the ring
$\Re[x,y,z]$.  The subscript $0$ in $\prob[0]{A}$ and
$\condprob[0]{B}{A}$ indicates that these expressions refer to
\emph{input} probabilities used to specify the probability model.  In
contrast, \emph{output} probabilities computed from these inputs are
written without a subscript.  For example the output table
$\condprob{B}{A}$ shown in Equation~\ref{eq:bga-output} was computed
from both input tables $\prob[0]{A}$ and $\condprob[0]{B}{A}$ from
Equation~\ref{eq:ab-input}.

\subsection{Symbolic Probability Inference}
\label{sec:spi}

By means of \emph{symbolic probability inference} it is possible to
compute polynomial expressions for queried probabilities from a
parametric probability network.  The essential rules of symbolic
probability inference were described well-enough several centuries ago
\cite{demoivre}.  Here some database terminology is enlisted to
describe the necessary operations (as in \cite{norman-thesis} and
\cite{norman-problogic}).  The first step is to combine the several
input probability tables into the full-joint probability distribution
over all the primary variables in the model.  For the parametric
probability network in Equation~\ref{eq:ab-input} this means
calculating the joint probabilities of both primary variables $A$ and
$B$.  The requisite operations can be viewed as a relational-database
join (designated $\Join$) after which the individual input
probabilities are aggregated by taking their algebraic products (as
though by a polynomial-aware analogue of the standard SQL aggregate
functions \texttt{SUM} and \texttt{AVERAGE}) \cite{date}:
\begin{equation}
\begin{tabular}[c]{l|l}\hline
\multicolumn{1}{l|}{$A$} & 
\multicolumn{1}{l}{$\prob[0]{{A}}$} \\ \hline\hline
$\state{T}$ & 
$x$ \\ \hline
$\state{F}$ & 
$1 - x$ \\ \hline
\end{tabular}
  \quad \Join \quad
\begin{tabular}[c]{ll|l}\hline
\multicolumn{1}{l}{$A$} & 
\multicolumn{1}{l|}{$B$} & 
\multicolumn{1}{l}{$\condprob[0]{{B}}{{A}}$} \\ \hline\hline
$\state{T}$ & 
$\state{T}$ & 
$y$ \\ \hline
$\state{T}$ & 
$\state{F}$ & 
$1 - y$ \\ \hline
$\state{F}$ & 
$\state{T}$ & 
$z$ \\ \hline
$\state{F}$ & 
$\state{F}$ & 
$1 - z$ \\ \hline
\end{tabular}
  \quad \Rightarrow \quad
\begin{tabular}[c]{r|ll|l}\hline
\multicolumn{1}{l|}{\scshape {\#}} & 
\multicolumn{1}{l}{$A$} & 
\multicolumn{1}{l|}{$B$} & 
\multicolumn{1}{l}{$\prob{{A, B}}$} \\ \hline\hline
1 & 
$\state{T}$ & 
$\state{T}$ & 
$x y$ \\ \hline
2 & 
$\state{T}$ & 
$\state{F}$ & 
$x - x y$ \\ \hline
3 & 
$\state{F}$ & 
$\state{T}$ & 
$z - x z$ \\ \hline
4 & 
$\state{F}$ & 
$\state{F}$ & 
$1 - x - z + x z$ \\ \hline
\end{tabular}
  \label{eq:ab-joint}
\end{equation}
For example the first and second elements of the full-joint
probability table $\prob{A,B}$ were computed as the products
\begin{math}
(
x
) \cdot (
y
)
\end{math}
and
\begin{math}
(
x
) \cdot (
1 - y
)
\end{math}.
The inputs in Equation~\ref{eq:ab-input} can be viewed as nothing more
or less than the specification of the full-joint probability table
$\prob{A,B}$ shown in Equation~\ref{eq:ab-joint}.  This table concerns
four elementary events defined by the various combinations of truth
and falsity of the logical propositions $A$ and $B$; it provides a
convenient factoring of the probabilities assigned to these four
elementary events.  The elements of the table $\prob{A,B}$ are subject
to the general laws of probability (each element is constrained to lie
between zero and one, and the sum of all elements is constrained to
equal one).

Each \emph{unconditioned} probability query yields the sum of selected
elements from the full-joint probability table (hence a polynomial
function of the parameters $x$, $y$, and $z$), and each
\emph{conditional} probability query yields the quotient of such sums
(hence a fractional polynomial).  For example, the unconditioned
probability that $A$ is true is computed as the sum of elements 1 and
2 of the full-joint probability table $\prob{A,B}$ from
Equation~\ref{eq:ab-joint}:
\begin{equation}
  \prob{A=\tval} \quad \Rightarrow \quad (
x y
  ) + (
x - x y
  ) \quad \Rightarrow \quad
x
  \label{eq:pr-a}
\end{equation}
The unconditioned probability that $B$ is true is the sum of elements
1 and 3 of the table $\prob{A,B}$:
\begin{equation}
  \prob{B=\tval} \quad \Rightarrow \quad (
x y
  ) + (
z - x z
  ) \quad \Rightarrow \quad
z + x y - x z
  \label{eq:pr-b}
\end{equation}
The conditional probability that $B$ is true given that $A$ is true is
computed by first recalling the definition
$\condprob{B}{A}=\prob{A,B}/\prob{A}$ from the general laws of
probability, and then evaluating the numerator and denominator
separately:
\begin{equation}
  \condprob{B=\tval}{A=\tval} \quad \Rightarrow \quad
  \frac{\prob{A=\tval,\:B=\tval}}{\prob{A=\tval}} 
  \quad \Rightarrow \quad
  \frac{
x y
  }{(
x y
    )+(
x - x y
    )
  }
  \quad \Rightarrow \quad
  \frac{
x y
  }{
x
  }
  \label{eq:bga-calc}
\end{equation}
You may be tempted to simplify this computed quotient $xy/x$ to the
elementary expression $y$, but such a premature step would discard
valuable information.  Because the constraints in
Equation~\ref{eq:ab-input} allow zero as a feasible value for $x$, the
quotient $xy/x$ could have the value $0/0$ (which is not the same
mathematical object as the expression $y$).  Avoiding premature
simplification allows a conditional probability to be recognized as
indefinite when its condition is impossible (in exactly the same sense
that the numerical quotient $0/0$ is indefinite).  For example, when
$\prob{A=\tval}=0$, we compute explicitly that
$\condprob{B=\tval}{A=\tval}$ does not have any particular real-number
value.

Here is the table of computed values for the probability query
$\condprob{B}{A}$, incorporating all four true/false combinations of
the variables $A$ and $B$:
\begin{equation}
\begin{tabular}[c]{l|ll}\hline
\multicolumn{3}{l}{$\condprob{{B}}{{A}}$} \\ \hline\hline
\multicolumn{1}{l|}{$A$} & 
\multicolumn{1}{l}{${B=\state{T}}$} & 
\multicolumn{1}{l}{${B=\state{F}}$} \\ \hline\hline
$\state{T}$ & 
$x y / x$ & 
$\left(x - x y\right) / \left(x\right)$ \\ \hline
$\state{F}$ & 
$\left(z - x z\right) / \left(1 - x\right)$ & 
$\left(1 - x - z + x z\right) / \left(1 - x\right)$ \\ \hline
\end{tabular}
\label{eq:bga-output}
\end{equation}
Note that the computed output probabilities in the table
$\condprob{B}{A}$ of Equation~\ref{eq:bga-output} are different
symbolic expressions from the input probabilities in the table
$\condprob[0]{B}{A}$ of Equation~\ref{eq:ab-input}.  The output
probabilities contain factors of $x/x$ or $(1-x)/(1-x)$, which make it
possible to catch the exception of division by zero that would be
caused by an impossible condition.

\subsection{Boolean Polynomial Translation}
\label{sec:boolean}

Boole provided a method to translate logical formulas of the
propositional calculus into polynomials with real coefficients
\cite{boole}.  Among other benefits, this formulation makes it
possible to perform logical deduction by solving equations.  Boolean
polynomial translation uses the following rules, with some notation
here modified from Boole's original text.  The rules are summarized in
Table~\ref{tbl:translation}.  Elementary truth maps to the real number
one, and elementary falsity maps to the real number zero.  Each
propositional variable $X_i$ translates into a real-valued algebraic
variable $x_i$.  Because two-valued logic mandates that each
$X_i\in\{\tval,\fval\}$, each translated $x_i$ is subject to the
constraint $x_i\in\{0,1\}$.  This constraint $x_i\in\{0,1\}$ may be
specified as the polynomial equation $x_i^2=x_i$ (which Boole called
the `fundamental law of thought' in \cite{boole}).  Propositional
functions map to arithmetical functions according to the rules shown
in Table~\ref{tbl:translation}, which were derived according to
Boole's polynomial interpolation method (which he called `function
development' in \cite{boole}).  The rules can be applied by top-down
or bottom-up parsing.  While simplifying polynomials translated from
logical formulas, any squared variable $x_i^2$ (or any higher power
$x_i^3$ and so on) can be replaced with the unadorned $x_i$ (because
$x_i^2=x_i$ for $x_i\in\{0,1\}$).

\begin{table}
  \sf
  \begin{tabular}{l|c|c} \hline
    \bfseries{Logic Function} & 
    \bfseries{Propositional Form} & \bfseries{Polynomial Form} 
    \\ \hline\hline
    Truth & $\tval$ & $1$ \\ \hline
    Falsity & $\fval$ & $0$ \\ \hline
    Atomic formula & $X_i$ & $x_i$ \\ \hline
    Negation & $\neg p$ & $1-p$ \\ \hline
    Conjunction & $p \wedge q$ & $pq$ \\ \hline
    Exclusive disjunction & $p \oplus q$ & $p+q-2pq$ \\ \hline
    Inclusive disjunction & $p \vee q$ & $p+q-pq$ \\ \hline
    Material implication & $p \rightarrow q$ & $1-p+pq$ \\ \hline
    Biconditional & $p \leftrightarrow q$ & $1-p-q+2pq$ \\ \hline
  \end{tabular}
  \caption{Rules for Boolean polynomial translation of
    propositional-calculus formulas.  Here $p$ and $q$ are polynomials
    in $\Re[x_1,\ldots,x_n]$ whose indeterminates match the
    propositional variables $X_1,\ldots,X_n$.  Because each
    $X_i\in\{\tval,\fval\}$, each $x_i\in\{0,1\}$.}
  \label{tbl:translation}
\end{table}

For example, consider the propositional-calculus formula $X \wedge (X
\rightarrow Y)$.  To prepare for translation, we declare polynomial
variables $x$ and $y$ which are subject to constraints $x\in\{0,1\}$
and $y\in\{0,1\}$.  The logical formula $X$ maps to the polynomial
$x$, and $Y$ maps to the polynomial $y$.  The material implication $x
\rightarrow y$ maps to the polynomial $1-x+xy$.  The conjunction $x
\wedge (1-x+xy)$ maps to the polynomial expression $x(1-x+xy)$, which
expands to the polynomial $x-x^2+x^2y$.  Taking advantage of the
identity $x^2=x$, this simplifies to $x-x+xy$, which then simplifies
to $xy$.  This polynomial result is a member of the ring $\Re[x,y]$.

Let us use $\bpoly{\phi}$ to denote the Boolean polynomial translation
of a propositional-calculus formula $\phi$.\footnote{The subscript 1
  indicates that elementary logical truth is mapped to the number one;
  there is an alternative translation scheme with truth mapped to the
  number zero.  Another option is to use polynomial coefficients in
  the finite field $\mathbb{F}_2=\{0,1\}$ instead of the real numbers
  $\Re$.  Numbers in $\mathbb{F}_2$ require integer arithmetic modulo
  two; there are some computational advantages to this choice.}  For
this example we write:
\begin{eqnarray}
  \bpoly{X \wedge (X \rightarrow Y)} & \Rightarrow & xy
\end{eqnarray}
understanding the real-valued variables $x$ and $y$ to be limited to
values in $\{0,1\}$.  Note that the polynomial $xy$ is also the
Boolean translation of the propositional-calculus formula $X \wedge
Y$, which has the same truth table as the formula $X \wedge (X
\rightarrow Y)$.  The polynomials generated by Boolean translation
provide a useful normal form for propositional-calculus formulas which
offers advantages over the traditional conjunctive and disjunctive
normal forms.


Boole's polynomial translation scheme offers a way to represent
unknown functions of the propositional calculus in a parametric
fashion.  For example with variables $x$ and $y$ and coefficients
$c_1, c_2, c_3, c_4$, each limited to real values in the set
$\{0,1\}$, the following polynomial expression represents an unknown
propositional-calculus function of two variables $X$ and $Y$:
\begin{equation}
  c_1 x y + c_2 x (1-y) + c_3 (1-x) y + c_4 (1-x) (1-y)
\end{equation}
To illustrate one instantiation, with $(c_1,c_2,c_3,c_4)=(1,0,1,1)$
the above expression simplifies to $1-x+xy$ indicating the statement
of material implication $X \rightarrow Y$.  Boole described this
technique in \cite{boole}.


\subsection{Embedding Formulas from the Propositional Calculus}
\label{sec:embedding}

We can extend parametric probability networks to incorporate logical
formulas from the propositional calculus, by using conditional
probability tables that mimic logical truth tables.  For example, let
us amend the probability model in Equation~\ref{eq:ab-input} to
include the formulas $A \rightarrow B$ and $A \wedge B$ from
the propositional calculus.  The truth tables of these statements of
material implication give the following conditional probability
tables:
\begin{equation}
\begin{tabular}[c]{ll|ll}\hline
\multicolumn{4}{l}{$\condprob[0]{{\embed{A \rightarrow B}}}{{A, B}}$} \\ \hline\hline
\multicolumn{1}{l}{$A$} & 
\multicolumn{1}{l|}{$B$} & 
\multicolumn{1}{l}{${\embed{A \rightarrow B}=\state{T}}$} & 
\multicolumn{1}{l}{${\embed{A \rightarrow B}=\state{F}}$} \\ \hline\hline
$\state{T}$ & 
$\state{T}$ & 
$1$ & 
$0$ \\ \hline
$\state{T}$ & 
$\state{F}$ & 
$0$ & 
$1$ \\ \hline
$\state{F}$ & 
$\state{T}$ & 
$1$ & 
$0$ \\ \hline
$\state{F}$ & 
$\state{F}$ & 
$1$ & 
$0$ \\ \hline
\end{tabular}
\qquad
\begin{tabular}[c]{ll|ll}\hline
\multicolumn{4}{l}{$\condprob[0]{{\embed{A \wedge B}}}{{A, B}}$} \\ \hline\hline
\multicolumn{1}{l}{$A$} & 
\multicolumn{1}{l|}{$B$} & 
\multicolumn{1}{l}{${\embed{A \wedge B}=\state{T}}$} & 
\multicolumn{1}{l}{${\embed{A \wedge B}=\state{F}}$} \\ \hline\hline
$\state{T}$ & 
$\state{T}$ & 
$1$ & 
$0$ \\ \hline
$\state{T}$ & 
$\state{F}$ & 
$0$ & 
$1$ \\ \hline
$\state{F}$ & 
$\state{T}$ & 
$0$ & 
$1$ \\ \hline
$\state{F}$ & 
$\state{F}$ & 
$0$ & 
$1$ \\ \hline
\end{tabular}
\end{equation}
Using these tables together with the original probability model in
Equation~\ref{eq:ab-input}, symbolic probability inference yields the
probabilities that the embedded propositional-calculus formulas are
true:
\begin{equation}
\begin{tabular}[c]{l|l}\hline
\multicolumn{1}{l|}{$\embed{A \rightarrow B}$} & 
\multicolumn{1}{l}{$\prob{{\embed{A \rightarrow B}}}$} \\ \hline\hline
$\state{T}$ & 
$1 - x + x y$ \\ \hline
$\state{F}$ & 
$x - x y$ \\ \hline
\end{tabular}
  \qquad
\begin{tabular}[c]{l|l}\hline
\multicolumn{1}{l|}{$\embed{A \wedge B}$} & 
\multicolumn{1}{l}{$\prob{{\embed{A \wedge B}}}$} \\ \hline\hline
$\state{T}$ & 
$x y$ \\ \hline
$\state{F}$ & 
$1 - x y$ \\ \hline
\end{tabular}
  \label{eq:pq-output}
\end{equation}
For convenience we can embed elementary logical truth and falsity
within a parametric probability model, using following conditional
probability tables:
\begin{equation}
\begin{tabular}[c]{l|l}\hline
\multicolumn{1}{l|}{$\embed{\tval}$} & 
\multicolumn{1}{l}{$\prob[0]{{\embed{\tval}}}$} \\ \hline\hline
$\state{T}$ & 
$1$ \\ \hline
$\state{F}$ & 
$0$ \\ \hline
\end{tabular}
\qquad
\begin{tabular}[c]{l|l}\hline
\multicolumn{1}{l|}{$\embed{\fval}$} & 
\multicolumn{1}{l}{$\prob[0]{{\embed{\fval}}}$} \\ \hline\hline
$\state{T}$ & 
$0$ \\ \hline
$\state{F}$ & 
$1$ \\ \hline
\end{tabular}
\end{equation}
These input probability tables say that truth is certainly true and
falsity is certainly false.  Note that conditioning on the truth of
embedded truth (or the falsity of embedded falsity) leaves the
polynomial formula computed for any probability expression unchanged:
this operation simply adds $1$ to the numerator and denominator of the
computed polynomial quotient.  For example, given any main variable
$A$ the conditional probability
$\condprob{A=\tval}{\embed{\tval}=\tval}$ evaluates to the same
polynomial expression as the unconditioned probability
$\prob{A=\tval}$.  Adding to consideration any variable $B$, the
conditional probability
$\condprob{A=\tval}{B=\tval,\embed{\tval}=\tval}$ yields the same
polynomial quotient as the simpler conditional probability
$\condprob{A=\tval}{B=\tval}$.

\subsection{Bounded Global Polynomial Optimization}
\label{sec:optimization}

Polynomials, such as those generated by symbolic probability inference
and those translated from the propositional calculus by Boole's
method, can be used to build optimization problems and solution sets
which will be useful for performing logical deduction.  Let us
formulate a general polynomial optimization problem with bounded
variables.  Consider a list $\mathbf{x}=(x_1,\ldots,x_n)$ of
real-valued variables with each $x_i$ bounded by finite lower and
upper limits $\alpha_i$ and $\beta_i$.  Consider also an objective
function $f$ and several constraint functions $g_1,\ldots,g_m$, each
of which is a polynomial in $\Re[\mathbf{x}]$.  These components
provide the following template for polynomial optimization problems:
\begin{equation}
  \begin{array}{rl}
    \mbox{Maximize}: & f(\mathbf{x})
    \\
    \mbox{subject to}: &
    g_1(\mathbf{x}) \ge 0, \ldots,
    g_m(\mathbf{x}) \ge 0
    \\
    \mbox{and}: &
    \alpha_1 \le x_1 \le \beta_1, \ldots,
    \alpha_n \le x_n \le \beta_n
  \end{array}
  \label{eq:pp}
\end{equation}
Certain variations on this template are allowed.  Any given
optimization problem may request either the minimum or maximum
feasible value of the objective, and each constraint may use either
equality $=$ or weak inequality $\ge$ or $\le$.  Strict inequalities
are approximated by the introduction of a small numerical constant
$\epsilon$ such as $0.001$ or $1 \times 10^{-6}$, using for example
$g(\mathbf{x}) \ge \epsilon$ to represent $g(\mathbf{x})>0$.  Any
variable $x_i$ may be restricted to take only integer values within
its range (such as $0$ and $1$ for a variable $x_i$ bounded by $0 \le
x_i \le 1$).

Optimization problems like Equation~\ref{eq:pp} can be challenging to
solve, because with nonlinear polynomials there can exist local
solutions which are not globally optimal.  The author has developed a
suitable algorithm for bounded global polynomial optimization, based
on earlier reformulation and linearization methods and dependent on a
separate mixed integer-linear programming solver \cite{norman-nlp}.
The associated software implementation was used to compute the results
presented herein.  There are some important theoretical details about
computational complexity and numerical approximation which are beyond
the scope of the present discourse.  Practically speaking, for the
small optimization problems addressed here, the author's solver takes
trivial amounts of time to compute solutions of ample numerical
precision.  Problems with inconsistent constraints are reported
directly as being infeasible.

Moving on, let us next consider a polynomial set-comprehension
expression which is related to the optimization problem in
Equation~\ref{eq:pp}.  Now, instead of seeking the minimum or maximum
feasible value of the objective $f$, we seek the set $\Phi$ of all
feasible values of the objective function subject to the constraints:
\begin{eqnarray}
  \Phi & \Leftarrow &
  \{ \: f(\mathbf{x}) \: : \:
  g(\mathbf{x}) \ge 0, \ldots g(\mathbf{x}) \ge 0;
  \alpha_1 \le x_1 \le \beta_1, \ldots,
  \alpha_n \le x_n \le \beta_n
  \: \}
  \label{eq:solset}
\end{eqnarray}
As before each constraint may use either relation $=$, $\ge$, or
$\le$; also any variable $x_i$ can be restricted to a finite set of
integer values.  By its construction, the solution set defined in
Equation~\ref{eq:solset} is a subset of the real numbers: $\Phi
\subseteq \Re$.

A solution set like Equation~\ref{eq:solset} can be described using a
pair of optimization problems like Equation~\ref{eq:pp}: one problem
to compute the minimum feasible value of the objective $f$, and
another to compute the maximum feasible value.  In the case that the
specified constraints are inconsistent, both optimization problems
will be infeasible.  Otherwise let us designate the computed minimum
and maximum feasible values of $f$ as $\alpha^*$ and $\beta^*$
respectively.  The solution set $\Phi$ must be a subset of the real
interval bounded by these computed solutions, and it must contain at
least those two points: hence $\Phi \subseteq [\alpha^*,\beta^*]$ and
$\alpha^*\in\Phi$ and $\beta^*\in\Phi$.  In the special case
$\alpha^*=\beta^*$ that the computed minimum and maximum values
coincide, then the solution set $\Phi$ is the singleton $\{\alpha^*\}$
containing that common value.

To illustrate, let us return to the basic probability model from
Equation~\ref{eq:ab-input}.  We choose as our objective $f(x,y,z)$ the
probability $\prob{B=\tval}$, which is the polynomial
\begin{math}
z + x y - x z
\end{math}
given in Equation~\ref{eq:pr-b}.  We choose as constraints the
parameter bounds given in Equation~\ref{eq:ab-input}, joined with two
additional constraints:
\begin{math}
x
  = 1
\end{math}
and
\begin{math}
x y
  =
x
\end{math}
(whose provenance will be discussed later).  Thus our goal is to find
the solution set $\Phi$ given by:
\begin{eqnarray}
  \Phi & \Leftarrow &
  \left\{ \:
z + x y - x z
  \: : \:
  x=1; \: xy=x; \:
  x,y,z \in [0,1]
  \: \right\}
  \label{eq:b-solset}
\end{eqnarray}
In order to characterize this solution set $\Phi$ we solve two
optimization problems patterned after Equation~\ref{eq:pp}:
\begin{equation}
\begin{array}{c@{\qquad}c}
\begin{array}{r@{\quad}l}
\mbox{Minimize}: & z + x y - x z \\
\mbox{subject to}:
& x = 1 \\
& x y = x \\
\mbox{and}:
& 0 \le x \le 1 \\
& 0 \le y \le 1 \\
& 0 \le z \le 1
\end{array}
&
\begin{array}{r@{\quad}l}
\mbox{Maximize}: & z + x y - x z \\
\mbox{subject to}:
& x = 1 \\
& x y = x \\
\mbox{and}:
& 0 \le x \le 1 \\
& 0 \le y \le 1 \\
& 0 \le z \le 1
\end{array}
\end{array}
\label{eq:ppb}
\end{equation}
The author's polynomial optimization solver determines that both
problems are feasible.  The computed minimum $\alpha^*$ is
\begin{math}
1.000
\end{math}
and the computed maximum $\beta^*$ is also
\begin{math}
1.000
\end{math}.
It follows from these results that the solution set $\Phi$ specified by
Equation~\ref{eq:b-solset} is a subset of the real interval
\begin{math}
[
1.000
,
1.000
]
\end{math} containing at least the point
\begin{math}
1.000
\end{math};
in other words $\Phi$ evaluates to the singleton $\{1\}$.  For the
set-comprehension expression in Equation~\ref{eq:b-solset} this result
is also evident from manual calculations.  For example, after
substituting $1$ for $x$ (according to the first constraint $x=1$) the
second constraint $xy=x$ becomes $y=1$.  Then when $1$ is substituted
for $x$ and also for $y$ the objective
\begin{math}
z + x y - x z
\end{math}
simplifies to the constant real number $1$.

\section{Types of Conditional Statements}
\label{sec:principal}

Taking advantage of the computational methods presented in
Section~\ref{sec:methods} we can develop formal definitions of several
different types of conditional statements.  These types and their
mathematical definitions are summarized in Table~\ref{tbl:principals}.
We define four principal types of conditionals in terms of constraints
on probabilities: these types are called `subjunctive', `material',
`existential', and `feasibility'.  Each of these probabilistic
conditionals includes a real-valued parameter $k$ between zero and one
that quantifies whether it is purely affirmative ($k=1$), purely
negative ($k=0$), or intermediate in its sense; the default is $k=1$
for an affirmative conditional.  Also we shall consider two types of
conditionals based on the propositional calculus.  `Truth-functional'
conditionals, the fifth principal type, use material implication.
`Boolean-feasibility' conditionals, the sixth principal type, use
Boole's algebraic method to analyze polynomial expressions translated
from propositional-calculus formulas.  Each truth-functional or
Boolean conditional includes a true/false parameter $K$, rather than a
real-valued parameter $k$, to indicate the affirmative ($K=\tval$)
versus the negative ($K=\fval$) sense; by default $K=\tval$.  The five
types of conditionals that are not called `subjunctive' are classified
as `indicative'.  Thus while there is one type of subjunctive
conditional, there are many types of indicative conditionals.

\begin{table}
  \sf
  \begin{tabular}{l|c|c|c} \hline
    \bfseries Type & 
    \bfseries Symbolic &
    \bfseries Probability or Logical Equations &
    \bfseries Affirmative ($k=1$), in Algebra
    \\ \hline\hline

    Subjunctive &
    \scond[k]{A}{B} &
    \begin{math}\begin{array}{@{}rcl@{}}
        \condprob[0]{B=\tval}{A=\tval} & = & k
    \end{array}\end{math}
    &
    \begin{math}
y
      = 1
    \end{math}
    \\ \hline

    Material &
    \mcond[k]{A}{B} &
    \begin{math}\begin{array}{@{}rcl@{}}
        \prob{A=\tval,B=\tval} & = & k \cdot \prob{A=\tval}
    \end{array}\end{math}
    &
    \begin{math}
x y
      =
x
    \end{math}
    \\ \hline

    Existential &
    \econd[k]{A}{B} &
    \begin{math}\begin{array}{@{}rcl@{}}
        \prob{A=\tval,B=\tval} & = & k \cdot \prob{A=\tval} \\
        \prob{A=\tval} & > & 0
    \end{array}\end{math}
    &
    \begin{math}
x y
      =
x
    \end{math}
    and
    \begin{math}
x
      > 0
    \end{math}
    \\ \hline

    Feasibility &
    \fcond[k]{\Gamma}{B} &
    \begin{math}\begin{array}{@{}rcl@{}}
        \{ \: \prob{B=\tval} \: : \: \Gamma \: \} & = & \{k\}
    \end{array}\end{math}
    &
    \begin{math}
      \left\{ \:
z + x y - x z
      \: : \: \Gamma, \: \Gamma_0
      \: \right\} = \{1\}
    \end{math}
    \\ \hline

    Quotient-feasibility &
    \qfcond[k]{\Gamma}{A}{B} &
    \begin{math}\begin{array}{@{}rcl@{}}
        \{ \: \condprob{B=\tval}{A=\tval} \: : \: \Gamma \: \} & = & \{k\}
    \end{array}\end{math}
    &
    \begin{math}
      \left\{ \:
x y / x
      \: : \: \Gamma, \: \Gamma_0
      \: \right\} = \{1\}
    \end{math}
    \\ \hline

    Truth-functional &
    \tcond[K]{A}{B} &
    \begin{math}\begin{array}{@{}rcl@{}}
        A \rightarrow (K \leftrightarrow B) & = & \tval
    \end{array}\end{math}
    &
    $1-a+ab = 1$
    \\ \hline

    Boolean-feasibility &
    \bcond[K]{A}{B} &
    \begin{math}\begin{array}{@{}rcl@{}}
        \{ \: B \: : \: A=\tval \: \} & = & \{K\}
    \end{array}\end{math}
    &
    \begin{math}
      \left\{ \:
      b \: : \: a=1, \: \Gamma'_0
      \: \right\} = \{1\}
    \end{math}
    \\ \hline

  \end{tabular}
  \caption{Several types of conditional statements relating true/false
    antecedent $A$ and consequent $B$, specializing the statement `If
    $A$ then $B$'.  The probabilistic conditionals are quantified by a
    fractional parameter $0 \le k \le 1$, with $k=1$ indicating
    affirmative and $k=0$ negative statements.  The symbol $\Gamma$
    denotes a set of polynomial equality and inequality constraints
    supplied by the user.  The displayed algebraic translations of
    affirmative conditionals were derived using symbolic probability
    inference from the probability model from
    Equation~\ref{eq:ab-input}, or Boolean polynomial translation from
    the propositional calculus.  The set $\Gamma_0$ consists of the
    constraints $x,y,z\in[0,1]$ required by the laws of probability
    for the model in Equation~\ref{eq:ab-input}; the set $\Gamma'_0$
    consists of the constraints $a,b\in\{0,1\}$ required by two-valued
    logic.}
  \label{tbl:principals}
\end{table}

In the following sections, the six principal types of conditionals are
defined mathematically.  The four probabilistic conditionals use
probability expressions as intermediate representations, and the two
others use the propositional calculus as an intermediate form.
Ultimately each type of conditional statement is translated into a
system of equations and inequalities involving polynomials with real
coefficients.  These algebraic formulas can be used for analysis.  Due
to their different algebraic definitions, the various types of
conditionals exhibit distinctive patterns of behavior (in terms of
what they say about conditional probabilities, how affirmative and
negative conditionals interact with one another, what information each
conditional statement provides about its antecedent, and what happens
with false antecedents).  These properties are discussed in
Section~\ref{sec:properties}.

Note that conditional statements need not have any temporal or causal
significance.  It is not necessary that the consequent event should
come after the antecedent event in time, nor is it necessary that the
antecedent should be a cause of the consequent.

\subsection{Subjunctive Conditionals}
\label{sec:subjunctive}

The \rname{subjunctive} conditional $\scond[k]{A}{B}$ with antecedent
term $A$, consequent term $B$, and fractional parameter $k$ is defined
by the equation that the input probability that $B$ is true given that
$A$ is true must equal $k$:
\begin{eqnarray}
  \condprob[0]{B=\tval}{A=\tval} & = & k
  \label{eq:subjunctive}
\end{eqnarray}
The \rname{subjunctive} conditional defined by
Equation~\ref{eq:subjunctive} says `If $A$ then $B$' in the following
more specific sense: `If there were any $A$, then $k$ would also be
$B$'.  With $k=1$ this becomes the affirmative statement
$\scond{A}{B}$ that if there were any $A$, then all would also be $B$.
With $k=0$ this becomes the negative statement $\scond[0]{A}{B}$ that
if there were any $A$, then none would also be $B$.  Using the basic
probability model from Equation~\ref{eq:ab-input}, the definition in
Equation~\ref{eq:subjunctive} becomes the following polynomial
equation:
\begin{eqnarray}
y
  & = & k
\end{eqnarray}
Hence
\begin{math}
y
  = 1
\end{math}
for the affirmative conditional and
\begin{math}
y
  = 0
\end{math}
for its negative counterpart.

Even when there is no input table $\condprob[0]{B}{A}$ in a particular
probability model relating the antecedent $A$ and consequent $B$, it
may still be possible to express a subjunctive conditional using that
model.  The essential requirement is that the polynomial quotient
computed as the value of the conditional probability
$\condprob{B=\tval}{A=\tval}$ must contain its denominator as a factor
in its numerator.  In the probability network graph the antecedent
must be an ancestor of the consequent.  The appropriate polynomial
expression to constrain is then constructed by eliminating this common
term from the numerator and denominator (without assuming it is
nonzero).  That is, the polynomial quotient computed for
$\condprob{B=\tval}{A=\tval}$ must have the form $pq/q$ for some
polynomials $p$ and $q$, where $q$ is also the computed value of
$\prob{A=\tval}$.  The subjunctive conditional is formed by
constraining the factor $p$ alone to equal the sense-parameter $k$.

Multiple antecedent terms can be accommodated in a \rname{subjunctive}
conditional by including additional variables in the condition part of
the input probability from Equation~\ref{eq:subjunctive}.  Thus the
\rname{subjunctive} conditional $\scond[k]{\{\:A_1,\ldots,A_m\:\}}{B}$ is
defined by the equation:
\begin{eqnarray}
  \condprob[0]{B=\tval}{A_1=\tval,\ldots,A_m=\tval} & = & k
\end{eqnarray}

For all types of conditionals it is acceptable to invert the truth or
falsity of any antecedent terms $A$ or $A_j$ and/or of the consequent
term $B$.  Thus for example to express the \rname{subjunctive}
interpretations of the conditionals `If $A$ then not-$B$' and `If
$A_1$ and not-$A_2$ then $B$' we would use:
\begin{eqnarray}
  \scond{A}{\overline{B}} & \equiv & 
  \condprob[0]{B=\fval}{A=\tval} = 1 \\
  \scond{\{\:A_1,\:\overline{A_2}\:\}}{B} & \equiv &
  \condprob[0]{B=\tval}{A_1=\tval,A_2=\fval} = 1
\end{eqnarray}
Note that there are two ways to oppose the affirmative conditional
$\scond{A}{B}$: first by inverting the truth of the consequent term
$B$, as in $\scond{A}{\overline{B}}$; second by using the fractional
parameter value $k=0$ instead of $k=1$, as in $\scond[0]{A}{B}$.  Both
approaches lead to equivalent polynomial constraints.  It is a
different thing to negate an entire conditional: the negation
$\neg(\scond[k]{A}{B})$ means the constraint
$\condprob[0]{B=\tval}{A=\tval} \neq k$.

\subsection{Material Conditionals}
\label{sec:material}

The \rname{material} conditional $\mcond[k]{A}{B}$ with antecedent
term $A$, consequent term $B$, and fractional parameter $k$ is defined
by the equation that the probability that $A$ and $B$ are both true
must equal $k$ times the probability that $A$ is true:
\begin{eqnarray}
  \prob{A=\tval,B=\tval} & = & k \cdot \prob{A=\tval}
  \label{eq:material}
\end{eqnarray}
The \rname{material} conditional defined by Equation~\ref{eq:material}
says `If $A$ then $B$' in the following more specific sense: `Either
there are no $A$, or $k$ of the $A$ are also $B$'.  With $k=1$ this is
the affirmative statement $\mcond{A}{B}$ that either all $A$ are also
$B$, or there are no $A$.  With $k=0$ this is the negative statement
$\mcond[0]{A}{B}$ that either no $A$ are also $B$, or there are no
$A$.  The \rname{material} conditional with an intermediate value of
its parameter $k$ asserts that either there are no $A$, or the
specified fraction $k$ of $A$ are also $B$.  Using the basic
probability model from Equation~\ref{eq:ab-input} the definition in
Equation~\ref{eq:material} becomes the polynomial equation:
\begin{eqnarray}
x y
 & = & k
x
\end{eqnarray}
Hence
\begin{math}
x y
  =
x
\end{math}
for the affirmative \rname{material} conditional and
\begin{math}
x y
  = 0
\end{math}
for its negative counterpart.  

Multiple antecedent terms can be accommodated in a \rname{material}
conditional by including additional variables in the probability
expressions from Equation~\ref{eq:material}.  Thus the
\rname{material} conditional $\mcond[k]{\{\:A_1,\ldots,A_m\:\}}{B}$ is
defined by the equation:
\begin{eqnarray}
  \prob{A_1=\tval,\ldots,A_m=\tval,B=\tval} & = & 
  k \cdot \prob{A_1=\tval,\ldots,A_m=\tval}
\end{eqnarray}

\subsection{Existential Conditionals}
\label{sec:existential}

The \rname{existential} conditional $\econd[k]{A}{B}$ adds to its
\rname{material} counterpart the requirement that the probability that
the antecedent $A$ is true cannot be zero.  Therefore the
\rname{existential} conditional with antecedent term $A$, consequent
term $B$, and fractional parameter $k$ is defined by the following
equation and inequality:
\begin{equation}
  \begin{array}{rcl}
    \prob{A=\tval,B=\tval} & = & k \cdot \prob{A=\tval} \\
    \prob{A=\tval} & > & 0
  \end{array}
  \label{eq:existential}
\end{equation}
The \rname{existential} conditional defined by
Equation~\ref{eq:existential} says `If $A$ then $B$' in the following
more specific sense: `There are some $A$, of which $k$ are also $B$'.
With $k=1$ this is the affirmative statement $\econd{A}{B}$ that all
$A$ are $B$, and there are some $A$.  With $k=0$ this is the negative
statement $\econd[0]{A}{B}$ that no $A$ are $B$, and there are some
$A$.  Intermediate values of $k$ assert that $B$ is true for the given
proportion of cases in which $A$ is also true.  Using the basic
probability model from Equation~\ref{eq:ab-input} the definition in
Equation~\ref{eq:existential} becomes the polynomial equation and
inequality:
\begin{equation}
  \begin{array}{rcl}
x y
    & = & k
x
    \\
x
    & > & 0
  \end{array}
\end{equation}

Equation~\ref{eq:existential} can be modified to accommodate multiple
antecedent terms in an \rname{existential} conditional.  The
\rname{existential} conditional $\econd[k]{\{\:A_1,\ldots,A_n\:\}}{B}$
is defined by:
\begin{equation}
  \begin{array}{rcl}
    \prob{A_1=\tval,\ldots,A_m=\tval,B=\tval} & = & 
    k \cdot \prob{A_1=\tval,\ldots,A_m=\tval} \\
    \prob{A_1=\tval,\ldots,A_m=\tval} & > & 0
  \end{array}
\end{equation}

\subsection{Feasibility Conditionals}
\label{sec:feasibility}

The \rname{feasibility} conditional $\fcond[k]{\Gamma}{B}$ with
antecedent $\Gamma$, consequent $B$, and parameter $k$ is defined by
the equation that $k$ is the only feasible value of the probability
that $B$ is true, subject to the polynomial equality and inequality
constraints in the set $\Gamma$:
\begin{eqnarray}
  \fcond[k]{\Gamma}{B} & \equiv &
  \left\{ \: \prob{B=\tval} \: : \: \Gamma \: \right\}
  = \{k\}
  \label{eq:feasibility}
\end{eqnarray}
The \rname{feasibility} conditional defined by
Equation~\ref{eq:feasibility} says `If $\Gamma$ then $B$' in the
following more specific sense: `Subject to the constraints in
$\Gamma$, the only feasible value of the probability that $B$ is true
is $k$'.  The symbolic probability inference methods from
Section~\ref{sec:spi} and the polynomial optimization methods from
Section~\ref{sec:optimization} can be used to evaluate the solution
set on the left-hand side of Equation~\ref{eq:feasibility}, and
thereby determine whether or not the equation is satisfied.  The
general constraints $\Gamma_0$ required by the laws of probability
must be included during analysis (restricting each probability to lie
between zero and one, and constraining appropriate sums of probability
to equal one).

Using the probability model from Equation~\ref{eq:ab-input} and
symbolic probability inference, the definition in
Equation~\ref{eq:feasibility} for the \rname{feasibility} conditional
$\fcond[k]{\Gamma}{B}$ becomes the following equation about sets of
real numbers:
\begin{eqnarray}
  \left\{ \:
z + x y - x z
  \: : \: \Gamma; \: 
  x,y,z \in [0,1]
  \: \right\}
  & = & \{k\}
\end{eqnarray}
Here the general constraints
$\Gamma_0=\{x\in[0,1],y\in[0,1],z\in[0,1]\}$ required by the laws of
probability for the model in Equation~\ref{eq:ab-input} have been
included.  

The antecedent set $\Gamma$ may include arbitrary polynomial equality
and inequality constraints, perhaps derived from other probability
expressions and probabilistic conditionals.  For example we might use
the constraint-set
\begin{math}
  \{ \: \prob{A=\tval}=1, \: \condprob[0]{B=\tval}{A=\tval}=1 \: \}
\end{math}
to provide the premises that $A$ is certainly true and that the
affirmative \rname{subjunctive} conditional $\scond{A}{B}$ from
Equation~\ref{eq:subjunctive} is true.  After symbolic probability
inference this example set $\Gamma$ consists of the polynomial
constraints
\begin{math}
  \{ \:
x
  = 1, \:
y
  = 1
  \: \}
\end{math}.

Note that unlike the parts of the other types of conditionals, the
antecedent $\Gamma$ and consequent $B$ of a \rname{feasibility}
conditional have different data types: $\Gamma$ is a set of polynomial
constraints using the real-valued parameters of the probability model
in use, whereas $B$ is a true/false variable that is one of the
primary variables in the probability model.  However we can use
$\fcond[k]{A}{B}$ to abbreviate the \rname{feasibility} conditional
$\fcond[k]{\{\prob{A=\tval}=1\}}{B}$ whose antecedent set $\Gamma$
consists of the solitary constraint that $A$ must certainly be true.

\subsubsection{Quotient-Feasibility Conditionals}
\label{sec:quotient-feasibility}

As a variation on the theme we can formulate a
\rname{quotient-feasibility} conditional $\qfcond[k]{\Gamma}{A}{B}$
which says that, subject to the constraints in the set $\Gamma$, the
only feasible value of the computed conditional probability that $B$
is true given that $A$ is true is $k$.  Thus the definition:
\begin{eqnarray}
  \qfcond[k]{\Gamma}{A}{B} & \equiv &
  \left\{ \: \condprob{B=\tval}{A=\tval} \: : \: \Gamma \: \right\}
  = \{k\}
  \label{eq:quotient-feasibility}
\end{eqnarray}
It is necessary to use \emph{fractional} polynomial optimization to
evaluate directly the solution set in this definition, since its
objective is a quotient of polynomials.  A suitable optimization
method is discussed in Appendix~\ref{sec:fractional}.  Notably by this
fractional polynomial optimization method, it is considered infeasible
for the denominator of the objective function to equal zero.
Therefore if the denominator is otherwise constrained to equal zero,
the entire problem is considered infeasible by the solver.  Otherwise,
the computed solution set describes the possible values of the
objective function when the denominator is not zero and the given
constraints $\Gamma$ are satisfied.  This arrangement means that the
conditional $\qfcond[k]{\Gamma}{A}{B}$ could be satisfied without
excluding zero as a possible value for the denominator
$\prob{A=\tval}$ subject to the constraints $\Gamma$; it is just that
other values besides zero would need to be possible too.

Alternatively, the \rname{quotient-feasibility} conditional can be
evaluated using two different solution sets.  In order for
the right-hand side of Equation~\ref{eq:quotient-feasibility} to be
satisfied, all three of the following must be true:
\begin{equation}
  \begin{array}{rcl}
    \left\{ \: \prob{A=\tval} \: : \: \Gamma \: \right\} & \neq & \emptyset
    \\
    \left\{ \: \prob{A=\tval} \: : \: \Gamma \: \right\} & \neq & \{0\}
    \\
    \left\{ \: \prob{A=\tval,B=\tval} - k \cdot \prob{A=\tval} \: : \: 
    \Gamma \: \right\} & = & \{0\}
  \end{array}
  \label{eq:quotient-feasibility-reformulated}
\end{equation}
In other words: the constraints $\Gamma$ must be feasible; it must be
feasible that the denominator probability $\prob{A=\tval}$ is not
zero, subject to $\Gamma$; and it must be true that the ratio
$\prob{A=\tval,B=\tval}/\prob{A=\tval}$ equals $k$ whenever it is
defined (that is, when $\prob{A=\tval}$ is not zero), also subject to
$\Gamma$.  This reformulation follows from the laws of probability.
Using the probability model in Equation~\ref{eq:ab-input}, the
\rname{quotient-feasibility} conditional
\begin{math}
  \qfcond[k]{\Gamma}{A}{B}
\end{math}
is given by the following three equations about sets of real numbers:
\begin{equation}
  \{\:
x
  \: : \: \Gamma; \: x,y,z\in[0,1] \: \} \; \neq \; \emptyset
  , \qquad
  \{\:
x
  \: : \: \Gamma; \: x,y,z\in[0,1] \: \} \; \neq \; \{0\}
  , \qquad
  \{\:
x y
  - k
x
  \: : \: \Gamma; \: x,y,z\in[0,1] \: \} \; = \; \{0\}
\end{equation}
These equations are similar to the \rname{existential} conditional
$\econd[k]{A}{B}$ from Equation~\ref{eq:existential}.

The \rname{quotient-feasibility} conditional
$\qfcond[k]{\Gamma}{\embed{\tval}}{B}$ whose antecedent is the
embedded propositional-calculus formula for elementary truth is
equivalent to the simple \rname{feasibility} conditional
$\fcond[k]{\Gamma}{B}$ defined above.

\subsection{Truth-Functional Conditionals}
\label{sec:truth-functional}

The affirmative \rname{truth-functional} conditional $\tcond{A}{B}$
with antecedent $A$ and consequent $B$ is defined as the logical
equation that the corresponding statement of material implication from
the propositional calculus is true:\footnote{The
  \rname{truth-functional} conditional is a special case of the
  probabilistic \rname{material} conditional.  Note $\prob{\embed{A
      \rightarrow B}=\tval}=1$ gives same constraint as
  $\prob{A=\tval,B=\tval}=\prob{A=\tval}$, and
  $\prob{\embed{A\rightarrow\neg B}=\tval}=1$ same as
  $\prob{A=\tval,B=\tval}=0$.  Convenient to calculate with
  $\mathbb{F}_2$ after Boolean translation of \rname{truth-functional}
  version.}
\begin{eqnarray}
  \tcond{A}{B} & \equiv &
  A \rightarrow B = \tval
  \label{eq:truth-functional}
\end{eqnarray}
Here $A$ and $B$ denote formulas from the propositional calculus;
these may be atomic formulas (propositional variables) or compound
formulas.  The affirmative \rname{truth-functional} conditional has
the same meaning as the affirmative ($k=1$) probabilistic
\rname{material} conditional: `Either $A$ is false or $B$ is true (or
both)'.

When $A$ and $B$ are atomic formulas, Boole's polynomial translation
method from Section~\ref{sec:boolean} yields the following polynomial
version of Equation~\ref{eq:truth-functional} representing the
conditional $\tcond{A}{B}$:
\begin{eqnarray}
  1-a+ab & = & 1
\end{eqnarray}
The real-valued variables $a$ and $b$ are subject to the constraints
$a\in\{0,1\}$ and $b\in\{0,1\}$.

A \rname{truth-functional} conditional is negated by inverting the
sense of its consequent term $B$.  Hence the negative
\rname{truth-functional} conditional $\tcond{A}{\overline{B}}$ is
defined by the equation that the propositional-calculus statement that
$A$ materially implies the negation of $B$ is true:
\begin{eqnarray}
  \tcond{A}{\overline{B}} & \equiv &
  A \rightarrow \neg B = \tval
  \label{eq:truth-functional-neg}
\end{eqnarray}
When $A$ and $B$ are atomic formulas, Boolean translation provides the
following polynomial equation for the negative
\rname{truth-functional} conditional $\tcond{A}{\overline{B}}$:
\begin{eqnarray}
  1-ab & = & 1
\end{eqnarray}
We can express the affirmative and negative \rname{truth-functional}
conditionals in an integrated way with the introduction of a logical
parameter $K$ which can be either true or false, using $K=\tval$ to
indicate an affirmative conditional and $K=\fval$ to indicate a
negative conditional:
\begin{eqnarray}
  \tcond[K]{A}{B} & \equiv &
  A \rightarrow (K \leftrightarrow B) = \tval
  \label{eq:truth-functional-k}
\end{eqnarray}
Note that the biconditional $\tval \leftrightarrow B$ has the same
truth table as $B$ itself, and the biconditional $\fval
\leftrightarrow B$ has the same truth table as the negation $\neg B$.
Boolean translation maps the true/false value $K$ to a real number $k$
which is either $0$ or $1$.  When $A$ and $B$ are atomic formulas, the
Boolean polynomial translation of the definiens in
Equation~\ref{eq:truth-functional-k} is as follows:
\begin{eqnarray}
  1-ak-ab+2kab & = & 1
\end{eqnarray}
Now $k=1$ designates an affirmative conditional and $k=0$ designates a
negative conditional.  Note that whereas the parameter $k$ for the
probabilistic conditionals above takes continuous values in the real
interval $[0,1]$, the parameter $k$ for \rname{truth-functional}
conditionals is limited to the integers $\{0,1\}$.

When the antecedent $A$ and consequent $B$ of a
\rname{truth-functional} conditional are not atomic formulas but
instead general propositional-calculus functions of some other
propositional variables $X_1,\ldots,X_n$, it is necessary to compute
the appropriate Boolean polynomial translation using the rules in
Table~\ref{tbl:translation}.  In this case the polynomial
translation of the \rname{truth-functional} conditional defined by
Equation~\ref{eq:truth-functional-k} is given by:
\begin{eqnarray}
  \tcond[K]{A}{B} & \equiv &
  \bpoly{A \rightarrow (K \leftrightarrow B)} = 1
\end{eqnarray}
using the value of $K$ supplied by the user (with $K=\tval$ for an
affirmative conditional and $K=\fval$ for a negative one).  The
translated propositional-calculus formula on the left-hand side of the
above equation will be a polynomial function of the real-valued
variables $x_1,\ldots,x_n$ corresponding to the propositional
variables used by the formulas $A$ and $B$.

Additional antecedent terms can be accommodated through logical
conjunction.  Considering antecedent formulas $A_1,\ldots,A_m$, the
\rname{truth-functional} conditional with consequent $B$ and
sense-parameter $K$ is defined by:
\begin{eqnarray}
  \tcond[K]{\{A_1,\ldots,A_m\}}{B} & \equiv &
  \bpoly{(A_1 \wedge \cdots \wedge A_m) \rightarrow (K \leftrightarrow B)} = 1
\end{eqnarray}

\subsection{Boolean-Feasibility Conditionals}
\label{sec:boolean-feasibility}

The \rname{Boolean-feasibility} conditional $\bcond[K]{A}{B}$ with
antecedent $A$, consequent $B$, and parameter $K$ is defined by the
constraint that $K$ is the only feasible truth value for $B$, subject
to the constraint that $A$ is true:
\begin{eqnarray}
  \bcond[K]{A}{B} & \equiv &
  \left\{ \:
  B \: : \: A=\tval
  \: \right\}
   = \{K\}
  \label{eq:boolean-feasibility}
\end{eqnarray}
Here $A$ and $B$ are formulas of the propositional calculus, and $K$
is either elementary truth value $\tval$ or $\fval$.  The
\rname{Boolean-feasibility} conditional defined by
Equation~\ref{eq:boolean-feasibility} says `If $A$ then $B$' in the
following more specific sense: `Assuming that $A$ is true, then $B$
must be $K$'.  The affirmative \rname{Boolean-feasibility} conditional
is indicated by $K=\tval$ and its negative counterpart by $K=\fval$.
The definiens of Equation~\ref{eq:boolean-feasibility} relates two
sets of elementary truth values, each of which is a (non-strict)
subset of $\{\tval,\fval\}$.  Note that
Equation~\ref{eq:boolean-feasibility} imposes both `positive' and
`negative' criteria on the possible truth values of the consequent $B$
subject to the constraint that the antecedent $A$ is true.  Positively
speaking, there must be at least one feasible case in which $A=\tval$
and $B=K$.  Negatively speaking, there must be no feasible case in
which $A=\tval$ and $B=\neg K$.  If there are no feasible cases in
which $A=\tval$ at all then the set $\{B:A=\tval\}$ is empty and the
conditional $\bcond[K]{A}{B}$ fails.  Here a `case' means a valuation
of the propositional variables involved.  We can have a qualitative
system in which we mark some valuations as possible and others as
impossible; or we can assign quantitative probabilities to the
possible valuations.  Implicit in this framing is the idea that there
may be a separate declaration of which valuations of the propositional
variables are feasible and which are not.  There are several possible
ways to supply such prior information, for example by adding
additional terms to the antecedent part of the conditional or by
introducing a probability distribution over the set of valuations of
the propositional variables.

Boolean polynomial translation of the definiens of
Equation~\ref{eq:boolean-feasibility} gives the following equation:
\begin{eqnarray}
  \left\{ \:
  \bpoly{B} \: : \: \bpoly{A}=1
  ; \: \Gamma_0
  \: \right\}
  & = & \{k\}
  \label{eq:boolean-feasibility-p}
\end{eqnarray}
where $k=1$ indicates an affirmative conditional and $k=1$ indicates a
negative conditional.  Reflecting two-valued logic, the set $\Gamma_0$
includes a constraint that limits the value of each variable to either
zero or one:
\begin{eqnarray}
  \Gamma_0 & \Leftarrow & \{x_1\in\{0,1\},\ldots,x_n\in\{0,1\}\}
\end{eqnarray}
Following Boole each constraint $x_i\in\{0,1\}$ could be specified as
$x_i^2=x_i$.  Equation~\ref{eq:boolean-feasibility-p} relates two sets
of real numbers, each a (non-strict) subset of the set $\{0,1\}$.
Note that there are four possible values for the solution set included
in Equation~\ref{eq:boolean-feasibility-p}: the set $\{0\}$; the set
$\{1\}$; the set $\{0,1\}$; and the empty set $\emptyset$.  This set
of solution sets can be exploited to develop a system of modal logic
\cite{norman-orthodoxy}.

Multiple antecedent terms can be accommodated in
\rname{Boolean-feasibility} conditionals through conjunction or
equivalently through additional constraints.  Thus we define:
\begin{eqnarray}
  \bcond[K]{\{A_1,\ldots,A_m\}}{B} & \equiv &
  \{ \: B \: : \: (A_1 \wedge \cdots \wedge A_m)=\tval \:\} = \{K\} \\
  & \equiv &
  \{ \: B \: : \: A_1=\tval,\ldots,\: A_m=\tval \:\} = \{K\}
\end{eqnarray}
The Boolean polynomial translation method from
Section~\ref{sec:boolean} could be applied to the conjunction $(A_1
\wedge \cdots \wedge A_m)$ following the first version, or applied to
the separate antecedent formulas $A_1$ through $A_m$ following the
second version.

\subsection{Relationships Among the Various Types of Conditionals}

Let us consider two sorts of relationships among the several types of
conditionals just defined.  First, the conditionals based on the
propositional calculus turn out to be closely related to probabilistic
conditionals: \rname{truth-functional} conditionals are special cases
of \rname{material} conditionals, and \rname{Boolean-feasibility}
conditionals are special cases of \rname{quotient-feasibility}
conditionals.  Second, there is a hierarchy among the simple
probabilistic conditionals: \rname{existential} conditionals subsume
\rname{subjunctive} ones, which in turn subsume \rname{material} ones.
\rname{Feasibility} conditionals stand outside this hierarchy; they
can be used as a metalevel statements to reason about conditionals of
other types.

\subsubsection{Truth-Functional and Material Conditionals}

The \rname{truth-functional} conditional $\tcond[K]{A}{B}$ defined by
Equation~\ref{eq:truth-functional} is a special case of the
probabilistic \rname{material} conditional $\mcond[k]{A}{B}$ from
Equation~\ref{eq:material} in the following sense.  Recall that the
affirmative ($K=\tval$) \rname{truth-functional} conditional is
defined as the equation that the propositional-calculus statement of
material implication $A \rightarrow B$ is true; similarly the negative
($K=\fval$) \rname{truth-functional} conditional is defined as the
equation that $A \rightarrow \neg B$ is true.  These
propositional-calculus formulas can be embedded within a probability
network by the technique from Section~\ref{sec:embedding}.  It happens
that the equation stating that the embedded formula $A \rightarrow B$
is certainly true is the same as the equation that specifies the
affirmative ($k=1$) \rname{material} conditional (according to
Definition~\ref{eq:material}):
\begin{eqnarray}
  \prob{\embed{A \rightarrow B}=\tval}=1 & \equiv &
  \prob{A=\tval,B=\tval}=\prob{A=\tval}
  \label{eq:material-equiv-affirmative}
\end{eqnarray}
Likewise the equation stating that the embedded formula $A \rightarrow
\neg B$ is true is the same as the equation specifying the negative
($k=0$) \rname{material} conditional:
\begin{eqnarray}
  \prob{\embed{A \rightarrow \neg B}=\tval}=1 & \equiv &
  \prob{A=\tval,B=\tval}=0
  \label{eq:material-equiv-negative}
\end{eqnarray}
In other words, the affirmative \rname{truth-functional} conditional
$\tcond{A}{B}$ specifies the same equation as the affirmative
\rname{material} conditional $\mcond{A}{B}$; and the negative
\rname{truth-functional} conditional $\tcond[\fval]{A}{B}$
(equivalently $\tcond{A}{\overline{B}}$) specifies the same equation
as the negative \rname{material} conditional $\mcond[0]{A}{B}$
(equivalently $\mcond{A}{\overline{B}}$).

The general relationships from Equations
\ref{eq:material-equiv-affirmative} an
\ref{eq:material-equiv-negative} hold using any parametric probability
model that includes the main variables $A$ and $B$ and additional main
variables for the embedded propositional-calculus formulas $A
\rightarrow B$ and $A \rightarrow \neg B$.  For example using the
probability model in Equation~\ref{eq:ab-input}, both affirmative
statements give the constraint
\begin{math}
x y = x
\end{math}
(symbolic probability inference yields
\begin{math}
  \prob{\embed{A \rightarrow B}=\tval} \Rightarrow
1 - x + x y
\end{math}).
Likewise using this same probability model both negative statements
give the constraints
\begin{math}
x y = 0
\end{math}
(symbolic probability inference yields
\begin{math}
  \prob{\embed{A \rightarrow \neg B}=\tval} \Rightarrow
1 - x y
\end{math}).  

Note that intermediate values of the sense-parameter $k$ allow
probabilistic \rname{material} conditionals to express constraints
that have no direct propositional-calculus counterparts.  For example
by Equation~\ref{eq:material} the \rname{material} conditional
$\mcond[0.5]{A}{B}$ specifies the equation
\begin{math}
  \prob{A=\tval,B=\tval} = 0.5 \cdot \prob{A=\tval}
\end{math},
 which does not have a simple relationship to the equation
\begin{math}
  \prob{\embed{A \rightarrow B}=\tval}=0.5
\end{math},
nor to the equation
\begin{math}
  \prob{\embed{A \rightarrow \neg B}=\tval}=0.5
\end{math}.  Using the probability model in Equation~\ref{eq:ab-input}
the respective polynomial equations are:
\begin{equation}
x y = 0.5 x
  , \qquad
0.5 + x y = x
  , \qquad
0.5 = x y
\end{equation}
These three equations express the respective ideas: first, `Either $A$
is false or half the time $B$ is true' (also, `Either $A$ is false or
half the time $B$ is false'); second, `Half the time, it is the case
that either $A$ is false or $B$ is true'; and third, `Half the time,
it is the case that either $A$ is false or $B$ is false'.  These are
some halfway points in various truth-functional interpretations of the
conditionals `If $A$ then $B$' and `If $A$ then not-$B$'.

\subsubsection{Boolean-Feasibility and Feasibility Conditionals}

The \rname{Boolean-feasibility} conditional $\bcond[K]{A}{B}$ defined
by Equation~\ref{eq:boolean-feasibility} is a special case of the
probabilistic \rname{quotient-feasibility} conditional from
Equation~\ref{eq:quotient-feasibility}, using an empty set $\Gamma$ of
additional constraints.  In this case the
affirmative \rname{quotient-feasibility} conditional
\begin{math}
  \qfcond{\emptyset}{A}{B}
\end{math}
and its affirmative \rname{Boolean-feasibility} counterpart
$\bcond{A}{B}$ both make the similar demand that, assuming that $A$ is
certainly true, truth must be the only feasible value for $B$:
\begin{eqnarray}
  \{\: \condprob{B=\tval}{A=\tval} \: \} & = & 1 \\
  \{\: B \: : \: A=\tval \: \} & = & \{\tval\}
\end{eqnarray}
Likewise the negative \rname{quotient-feasibility} conditional
\begin{math}
  \qfcond[0]{\emptyset}{A}{B}
\end{math}
(equivalently
\begin{math}
  \qfcond{\emptyset}{A}{\overline{B}}
\end{math})
and its \rname{Boolean-feasibility} counterpart $\bcond[\fval]{A}{B}$
(equivalently $\bcond{A}{\overline{B}}$) both make the similar demand
that, assuming that $A$ is certainly true, falsity must be the only
feasible value for $B$:
\begin{eqnarray}
  \{\: \condprob{B=\tval}{A=\tval} \: \} & = & 0 \\
  \{\: B \: : \: A=\tval \: \} & = & \{\fval\}
\end{eqnarray}
Recall from Section~\ref{sec:quotient-feasibility} and
Appendix~\ref{sec:fractional} that an optimization problem with the
fractional polynomial objective $\condprob{B=\tval}{A=\tval}$, meaning
the quotient $\prob{A=\tval,B=\tval}/\prob{A=\tval}$, is considered
infeasible if zero is the only feasible value of the objective's
denominator $\prob{A=\tval}$.  Hence both types of conditionals
$\bcond[K]{A}{B}$ and $\qfcond[k]{\emptyset}{A}{B}$ are false if it is
certain \emph{a priori} that the antecedent $A$ cannot be true,
regardless of the value of the respective parameter $K$ or $k$.
However if it is feasible but not necessary that $\prob{A=\tval}=0$
(meaning $0 \in \{\: \prob{A=\tval} \:\}$ and $\{\: \prob{A=\tval} \:\}
\neq \{0\}$), then the indefinite value $0/0$ is a feasible value of
the computed conditional probability $\condprob{B=\tval}{A=\tval}$
even if a conditional of either type is satisfied.  If desired, an
explicit constraint $\prob{A=\tval}>0$ could be added in order to
forbid this possibility of an indefinite quotient.

Propositional-calculus conditionals ($\rname{truth-functional}$ and
$\rname{Boolean-feasibility}$) allow some results to be computed by
techniques which are different from the algebraic and numerical
calculations required by their probabilistic counterparts
($\rname{material}$ and $\rname{Boolean-feasibility}$).  For small
problems the propositional-calculus formulations support direct and
exhaustive enumeration of various sets of formulas and valuations.
For larger problems, Boolean translation of propositional-calculus
formulas (now into polynomials with coefficients in the finite field
$\mathbb{F}_2$ instead of the real numbers $\Re$) facilitates
automated search through finite sets of equivalence classes of logical
formulas.

\subsubsection{Hierarchy Among Probabilistic Conditionals}

Next let us discuss the subsumption relationships among conditionals
of three of the probabilistic types (\rname{subjunctive},
\rname{material}, and \rname{existential}), using their affirmative
variants.  Consider two binary distinctions: first $\prob{A=\tval}=0$
versus $\prob{A=\tval}>0$; and second
$\condprob[0]{B=\tval}{A=\tval}<1$ versus
$\condprob[0]{B=\tval}{A=\tval}=1$.  Using the probability model in
Equation~\ref{eq:ab-input}, let us denote the first distinction as
\begin{math}
x
  = 0
\end{math}
versus
\begin{math}
x
  > 0
\end{math},
and the second distinction as
\begin{math}
y
  < 1
\end{math}
versus
\begin{math}
y
  = 1
\end{math}.
Integrating these two binary distinctions gives four possible
configurations:
\begin{equation}
  \begin{array}{|c|c|} \hline
x
    = 0, \:
y
    = 1
    &
x
    > 0, \:
y
    = 1
    \\ \hline
x
    = 0, \:
y
    < 1
    &
x
    > 0, \:
y
    < 1
    \\ \hline
  \end{array}
\end{equation}
The \rname{material} conditional defined by Equation~\ref{eq:material}
is satisfied in three of these four configurations (either
\begin{math}
x
  = 0
\end{math},
or
\begin{math}
y
  = 1
\end{math});
it is the least restrictive statement.  The \rname{subjunctive}
conditional defined by Equation~\ref{eq:subjunctive} is satisfied in
two of these four configurations
(\begin{math}
y
  = 1
\end{math});
it is intermediate.  The \rname{existential} conditional defined
by Equation~\ref{eq:existential} is the most restrictive statement,
being satisfied by only one of the four configurations
(\begin{math}
x
  > 0
\end{math}
and
\begin{math}
y
  = 1
\end{math}).

Therefore we can identify the affirmative forms of the three basic
probabilistic conditionals by the configurations of the values
\begin{math}
x
\end{math}
(representing $\prob{A=\tval}$) and
\begin{math}
y
\end{math}
(representing $\condprob[0]{B=\tval}{A=\tval}$) consistent with each:
\begin{equation}
  \begin{array}{|c|c|} \hline
    \multicolumn{2}{|c|}{\rname{material} \quad \mcond{A}{B}} 
    \\ \hline\hline
x
    = 0, \:
y
    = 1
    &
x
    > 0, \:
y
    = 1
    \\ \hline
x
    = 0, \:
y
    < 1
    &
    \\ \hline
  \end{array}
  \qquad
  \begin{array}{|c|c|} \hline
    \multicolumn{2}{|c|}{\rname{subjunctive} \quad \scond{A}{B}} 
    \\ \hline\hline
x
    = 0, \:
y
    = 1
    &
x
    > 0, \:
y
    = 1
    \\ \hline
    &
    \\ \hline
  \end{array}
  \qquad
  \begin{array}{|c|c|} \hline
    \multicolumn{2}{|c|}{\rname{existential} \quad \econd{A}{B}} 
    \\ \hline\hline
    \phantom{x=0,\:y=1}
    &
x
    > 0, \:
y
    = 1
    \\ \hline
    &
    \\ \hline
  \end{array}
\end{equation}
You can see regarding these affirmative conditionals that if the
\rname{existential} conditional holds then the \rname{subjunctive}
conditional must also, and in turn if the \rname{subjunctive}
conditional holds then the \rname{material} conditional must also.

\section{Distinctive Features of Conditionals}
\label{sec:properties}

Let us compare a few features of the various types of conditionals
just defined.  We shall examine four features:
\begin{itemize}
\item \emph{Conditional probability}: what a probabilistic conditional
  `If $A$ then $B$' says about the conditional probability that $B$ is
  true given that $A$ is true.
\item \emph{Consistency of opposites}: what happens when opposing
  conditionals `If $A$ then $B$' and `If $A$ then not-$B$' are
  asserted together.
\item \emph{Existential import}: what information the conditionals `If
  $A$ then $B$' and `If $A$ then not-$B$' (considered separately or
  together) convey about the truth of their antecedent $A$.
\item \emph{False antecedents}: what a false antecedent $A$ means for
  the conditionals `If $A$ then $B$' and `If $A$ then not-$B$'.
\end{itemize}
Different types of conditionals exhibit different patterns of
behavior.  No particular set of features is intrinsically right or
wrong.  Instead it behooves the user to choose the type of conditional
statement that provides the desired features for any given instance of
analysis.  Table~\ref{tbl:properties} summarizes the features which
are discussed.

\begin{table}
  \sf
  \begin{tabular}{l|llll} \hline
    \thd{Type} & 
    \begin{tabular}{@{}l@{}}
      \thd{Cond.\ Prob.} \\ $\condprob{B=\tval}{A=\tval}$
    \end{tabular}
    &
    \begin{tabular}{@{}l@{}}
      \thd{Opposites} \\ $k=0$, $k=1$
    \end{tabular}
    &
    \begin{tabular}{@{}l@{}}
      \thd{Existential} \\ \thd{Import}
    \end{tabular}
    &
    \begin{tabular}{@{}l@{}}
      \thd{False Antecedent} \\ $\prob{A=\tval}=0$
    \end{tabular}
    \\ \hline\hline
    Subjunctive & $k$ or $0/0$ & inconsistent & none & 
    $\scond[k]{A}{B}$ unaffected
    \\ \hline
    Material & $k$ or $0/0$ & consistent & indirect & 
    $\mcond[k]{A}{B}$ true
    \\ \hline
    Existential & $k$ & inconsistent & direct & 
    $\econd[k]{A}{B}$ false
    \\ \hline
    Feasibility & $k$ or $0/0$ & inconsistent & direct &
    $\fcond[k]{A}{B}$ false
    \\ \hline
  \end{tabular}
  \caption{Features of various types of probabilistic conditionals
    with antecedent $A$ and consequent $B$: the possible values of the
    computed conditional probability shown; the result of asserting
    opposite (affirmative and negative) conditionals simultaneously;
    the type of existential import; and the result of an antecedent
    that is certainly false.}
  \label{tbl:properties}
\end{table}

\subsection{Implications for Conditional Probability}

Each probabilistic conditional with antecedent $A$, consequent $B$,
and parameter $k$ says in a slightly different way that the
conditional probability that $B$ is true given that $A$ is true equals
$k$.  Recall from Equation~\ref{eq:bga-calc} in Section~\ref{sec:spi}
that the conditional probability $\condprob{B=\tval}{A=\tval}$ is
computed from the probability model in Equation~\ref{eq:ab-input} as
the following quotient:
\begin{eqnarray}
  \condprob{B=\tval}{A=\tval}
  & \Rightarrow &
  \frac{
x y
  }{
x
  }
  \label{eq:condprob}
\end{eqnarray}
Here
\begin{math}
x
\end{math}
is the input $\prob[0]{A=\tval}$ and
\begin{math}
y
\end{math}
is the input $\condprob[0]{B=\tval}{A=\tval}$, both from
Equation~\ref{eq:ab-input}.  There are three different ways to fix the
value of this quotient
\begin{math}
x y / x
\end{math}
to equal some specified constant value $k$: the \rname{subjunctive}
conditional $\scond[k]{A}{B}$ says
\begin{math}
y
  = k
\end{math}; 
the \rname{material} conditional $\mcond[k]{A}{B}$ says
\begin{math}
x y
  = k
x
\end{math}; 
and the \rname{existential} conditional $\econd[k]{A}{B}$ says both
\begin{math}
x y
  = k
x
\end{math}
and
\begin{math}
x
  > 0
\end{math}.
The \rname{subjunctive} and \rname{material} constraints allow two
possible results for the computed value of
$\condprob{B=\tval}{A=\tval}$: either this value is $k$, or it is
undefined (because of division by zero):
\begin{eqnarray}
  \condprob{B=\tval}{A=\tval} & \Rightarrow &
  \left\{
  \begin{array}{rl}
    k, & \mbox{if $\prob{A=\tval}>0$} \\
    0/0, & \mbox{if $\prob{A=\tval}=0$}
  \end{array}
  \right.
  \label{eq:condprob-2}
\end{eqnarray}
On the other hand the \rname{existential} constraint requires that the
computed conditional probability must have the definite value $k$
(since having the denominator equal zero is specifically forbidden):
\begin{eqnarray}
  \condprob{B=\tval}{A=\tval} & \Rightarrow & k
\end{eqnarray}

In contrast to their \rname{material} and \rname{existential}
counterparts, probabilistic \rname{feasibility} conditionals express a
different idea of `conditional' probability---now using the algebraic
operation of constraint instead of the arithmetical operation of
division as the mechanism to implement conditioning.  Thus instead of
considering the possible values of the conditional-probability
expression $\condprob{B=\tval}{A=\tval}$, for the \rname{feasibility}
conditional we consider the possible values of the
unconditioned-probability expression $\prob{B=\tval}$ subject to the
constraint that the unconditioned-probability expression
$\prob{A=\tval}$ must equal $1$.  Following
Equation~\ref{eq:feasibility}, the \rname{feasibility} version of the
conditional `If $A$ then $B$' is defined using the template 
$\fcond[k]{\Gamma}{B}$ with the antecedent constraint-set
$\Gamma\Leftarrow\{\:\prob{A=\tval}=1\:\}$ containing the assertion
that $A$ must certainly be true:
\begin{eqnarray}
  \fcond[k]{\{\:\prob{A=\tval}=1\:\}}{B} & \equiv &
  \{\: \prob{B=\tval} \: : \: \prob{A=\tval}=1 \:\} = \{k\}
\end{eqnarray}
This particular \rname{feasibility} conditional can be written in
abbreviated form as $\fcond[k]{A}{B}$.  There are two possible
computed values for the solution set that is constrained by
$\fcond[k]{A}{B}$.  If it is feasible that $\prob{A=\tval}=1$
(satisfying also the general constraints of probability for the model
in use), then the solution set will have the computed value $\{k\}$;
otherwise the computed solution set will be empty.  Thus in parallel
with Equation~\ref{eq:condprob-2} above we have two possible
consequences of the asserted \rname{feasibility} conditional
$\fcond[k]{A}{B}$:
\begin{eqnarray}
  \{\: \prob{B=\tval} \: : \: \prob{A=\tval}=1 \:\} & \Rightarrow &
  \left\{
  \begin{array}{rl}
    \{k\}, & \mbox{if $\prob{A=\tval}=1$ is feasible} \\
    \emptyset, & 
    \begin{array}[t]{@{}l@{}}
      \mbox{if $\prob{A=\tval}=1$ is infeasible} \\
      \mbox{(but $\prob{A=\tval}>0$ is feasible)}
    \end{array}
  \end{array}
  \right.
  \label{eq:solset-2}
\end{eqnarray}
After symbolic probability inference using the probability model in
Equation~\ref{eq:ab-input}, the \rname{feasibility} conditional
$\fcond[k]{A}{B}$ says:
\begin{eqnarray}
  \{\:
z + x y - x z
  \: : \:
x
  = 1; \:
  x,y,z\in[0,1]
  \:\} & = & \{k\}
\end{eqnarray}
Note that after substituting $1$ for
\begin{math}
x
\end{math}
the objective
\begin{math}
z + x y - x z
\end{math}
simplifies to the polynomial
\begin{math}
y
\end{math}.
Thus using the probability model in Equation~\ref{eq:ab-input}, the
conditional $\fcond[k]{A}{B}$ requires that
\begin{math}
y
  = k
\end{math}
and that it is feasible that
\begin{math}
x
  = 1
\end{math}.

\subsection{Consistency of Opposites}

Whether it is consistent or inconsistent to assert opposing
conditionals like `If $A$ then $B$' and `If $A$ then not-$B$' depends
upon the type of the conditional statements.\footnote{Note that with
  every type of conditional, it is different to negate the entire
  conditional statement versus negating just its consequent.  That is,
  the statement `Not (If $A$ then $B$)' is mathematically distinct
  from the statement `If $A$ then not-$B$' in every interpretation.}
For four of the six principal types of conditionals, it is
inconsistent in a simple algebraic sense to assert opposing
conditionals; for the other two types it is consistent.  For
\rname{subjunctive} conditionals as defined by
Equation~\ref{eq:subjunctive} the affirmative statement $\scond{A}{B}$
(meaning $\scond[1]{A}{B}$) and the negative statement
$\scond[0]{A}{B}$ together specify the following pair of equations:
\begin{eqnarray}
  \condprob[0]{B=\tval}{A=\tval} & = & 1 \\
  \condprob[0]{B=\tval}{A=\tval} & = & 0
\end{eqnarray}
meaning that the same input value must simultaneously equal the real
number $1$ and the real number $0$; this is obviously impossible.  For
example using the probability model from Equation~\ref{eq:ab-input}
the relevant equations are
\begin{math}
y
  = 1
\end{math}
and
\begin{math}
y
  = 0
\end{math},
which cannot be satisfied simultaneously by any real value of $y$.

For \rname{existential} conditionals as defined by
Equation~\ref{eq:existential} the affirmative statement $\econd{A}{B}$
and the negative statement $\econd[0]{A}{B}$ together specify the
following system of equations and inequalities:
\begin{eqnarray}
  \prob{A=\tval,B=\tval} & = & \prob{A=\tval} \\
  \prob{A=\tval,B=\tval} & = & 0 \\
  \prob{A=\tval} & > & 0
\end{eqnarray}
This system is algebraically inconsistent.  The first two equations
require $\prob{A=\tval}$ to equal $0$ yet the final inequality
requires $\prob{A=\tval}$ to be strictly greater than $0$.  For
example using the probability model from Equation~\ref{eq:ab-input}
the relevant constraints are
\begin{math}
x y
  = 
x
\end{math},
\begin{math}
x y
  =  0
\end{math}, and
\begin{math}
x
  > 0
\end{math},
which cannot be satisfied by any real values of $x$ and $y$.

For \rname{feasibility} conditionals the affirmative statement
$\fcond{A}{B}$ and the negative statement $\fcond[0]{A}{B}$ as defined
by Equation~\ref{eq:feasibility} together
specify the following pair of equations:
\begin{eqnarray}
  \{\: \prob{B=\tval} \: : \: \prob{A=\tval}=1 \:\} & = & \{1\} \\
  \{\: \prob{B=\tval} \: : \: \prob{A=\tval}=1 \:\} & = & \{0\}
\end{eqnarray}
meaning that the same set of real numbers must simultaneously equal
$\{1\}$ and $\{0\}$, which is algebraically impossible.  For example
using the probability model from Equation~\ref{eq:ab-input} the
relevant equations are:
\begin{eqnarray}
  \{\:
z + x y - x z
  \: : \:
x
  = 1; \:
  x,y,z\in[0,1]
  \:\} & = & \{1\}
  \\
  \{\:
z + x y - x z
  \: : \:
x
  = 1; \:
  x,y,z\in[0,1]
  \:\} & = & \{0\}
\end{eqnarray}
There are no real values $(x,y,z)$ that satisfy both equations; that
would require
\begin{math}
y
  = 1
\end{math}
and
\begin{math}
y
  = 0
\end{math}
simultaneously.

Likewise for \rname{Boolean-feasibility} conditionals the affirmative
statement $\bcond{A}{B}$ and the negative statement
$\bcond[\fval]{A}{B}$ (meaning $\bcond{A}{\neg B}$) together specify
the following pair of equations:
\begin{eqnarray}
  \{\: B \: : \: A=\tval \:\} & = & \{\tval\} \\
  \{\: B \: : \: A=\tval \:\} & = & \{\fval\}
\end{eqnarray}
This system is inconsistent because the same set of truth values
cannot simultaneously equal $\{\tval\}$ and $\{\fval\}$.  Using the
Boolean polynomial translation method from Section~\ref{sec:boolean},
these equations about logical formulas become equations about
polynomial formulas:
\begin{eqnarray}
  \{\: b \: : \: a=1; \: a^2=a,\: b^2=b \:\} & = & \{1\} \\
  \{\: b \: : \: a=1; \: a^2=a,\: b^2=b \:\} & = & \{0\}
\end{eqnarray}
There are no real values $(a,b)$ that satisfy both equations; that
would require $b=1$ and $b=0$ simultaneously.

In contrast, opposing conditionals of the remaining two types are not
inconsistent with one another.  For \rname{material} conditionals as
defined by Equation~\ref{eq:material} the affirmative statement
$\mcond{A}{B}$ and the negative statement $\mcond[0]{A}{B}$ together
specify the following pair of equations using probability expressions
(which map to polynomials using the symbolic probability inference
method of Section~\ref{sec:spi}):
\begin{eqnarray}
  \prob{A=\tval,B=\tval} & = & \prob{A=\tval} \\
  \prob{A=\tval,B=\tval} & = & 0
\end{eqnarray}
These equations are satisfied simultaneously when $\prob{A=\tval}=0$.
For example using the probability model from
Equation~\ref{eq:ab-input} the relevant equations are
\begin{math}
x y
  = 
x
\end{math}
and
\begin{math}
x y
  =  0
\end{math},
which are both satisfied when
\begin{math}
x
  = 0
\end{math}.
For \rname{truth-functional} conditionals as defined by
Equation~\ref{eq:truth-functional} the affirmative statement
$\tcond{A}{B}$ and the negative statement $\tcond[\fval]{A}{B}$
together specify the following pair of equations using
propositional-calculus formulas (which map to polynomials using the
Boolean translation method of Section~\ref{sec:boolean}):
\begin{eqnarray}
  (A \rightarrow B) & = & \tval \\
  (A \rightarrow \neg B) & = & \tval
\end{eqnarray}
Both formulas $A \rightarrow B$ and $A \rightarrow \neg B$ are true
when $A$ is false.  For example after Boolean translation these two
equations about logical formulas become a system of equations about
polynomials with real coefficients:
\begin{eqnarray}
  1-a+ab & = & 1 \\
  1-a+a(1-b) & = & 1 \\
  a^2 & = & a \\
  b^2 & = & b
\end{eqnarray}
You can see by inspection that this system is satisfied for $a=0$ and
any value of $b$.  The number $0$ of course represents logical falsity
in Boole's translation scheme.

\subsection{Existential Import}

Something interesting happens when opposing conditionals of the
\rname{material} or \rname{truth-functional} types are combined: they
tell us that the antecedent must be false.  In the \rname{material}
case the opposing conditionals $\mcond{A}{B}$ (meaning
$\mcond[1]{A}{B}$) and $\mcond[0]{A}{B}$ are satisfied simultaneously
if and only if $\prob{A=\tval}=0$ (meaning that the antecedent $A$ is
certainly false).  In fact this conclusion would follow from any
material conditionals $\mcond[k_1]{A}{B}$ and $\mcond[k_2]{A}{B}$ with
different sense-parameters $k_1 \neq k_2$.  In the
\rname{truth-functional} case the opposing conditionals $\tcond{A}{B}$
and $\tcond[\fval]{A}{B}$ (meaning $\tcond{A}{\neg{B}}$) are satisfied
simultaneously if and only if $A=\fval$ (meaning that the antecedent
$A$ is identically false).  Let us call this phenomenon `indirect
existential import': whereas no individual conditional statement
limits the possible values of the antecedent, pairs of opposite
conditionals indeed do so.

In contrast, \rname{existential}, \rname{feasibility}, and
\rname{Boolean-feasibility} conditionals possess `direct existential
import'.  If an individual conditional of one of these types is true,
then it must be at least \emph{possible} for the antecedent to be true
(though it need not be \emph{certain} that the antecedent is true).
An \rname{existential} conditional $\econd[k]{A}{B}$ requires directly
$\prob{A=\tval}>0$, stating that the probability that the antecedent
$A$ is true must be strictly greater than zero.  A \rname{feasibility}
conditional $\fcond[k]{\Gamma}{B}$ requires it to be possible for the
antecedent constraints $\Gamma$ to be satisfied.  Using
$\{\prob{A=\tval}=1\}$ as the antecedent constraint-set $\Gamma$, the
resulting \rname{feasibility} conditional $\fcond[k]{A}{B}$ requires
that it is feasible for the the antecedent $A$ to be certainly true
(which is different from the requirement $\prob{A=\tval}>0$).
Likewise the \rname{Boolean-feasibility} conditional $\bcond{A}{B}$
requires that it is possible for its antecedent $A$ to be true.

Finally, \rname{subjunctive} conditionals have no intrinsic
existential import of either the direct or the indirect kind.  In
order to express a \rname{subjunctive} conditional the probability
model must have the antecedent $A$ as an ancestor of the consequent
$B$ in the probability network graph, as in the model in
Equation~\ref{eq:ab-input}.  With this setup the input
$\condprob[0]{B=\tval}{A=\tval}$ does not concern the probability
$\prob{A=\tval}$ at all; when $\prob{A=\tval}$ is calculated by
symbolic probability inference (as described in Section~\ref{sec:spi})
the input $\condprob[0]{B=\tval}{A=\tval}$ cancels out algebraically.
A particular probability model could introduce direct existential
import to a \rname{subjunctive} conditional by using the same
parameters to specify the input $\condprob[0]{B=\tval}{A=\tval}$ and
for example the input $\prob[0]{A=\tval}$, or by including constraints
relating these input values.  In the absence of such features,
\rname{subjunctive} conditionals remain free of existential import.

\subsection{Consequences of False Antecedents}

Conditionals of the various types behave differently when their
antecedents are false.  Usually, the truth or falsity of its
antecedent $A$ has no effect on the truth or falsity of a
\rname{subjunctive} conditional $\scond{A}{B}$.  As stated in the last
section the input $\condprob[0]{B=\tval}{A=\tval}$ constrained by
$\scond{A}{B}$ is unrelated to the computed probability
$\prob{A=\tval}$, unless the probability model includes specific
constraints (or common parameters) that link these two values.

With $\prob{A=\tval}=0$, a \rname{material} conditional
$\mcond[k]{A}{B}$ using any value of $k$ must be true.  In this case
Equation~\ref{eq:material} defining the conditional $\mcond[k]{A}{B}$
simplifies to the trivial requirement $0=0$.  Likewise with $A=\fval$
both \rname{truth-functional} conditionals $\tcond{A}{B}$ and
$\tcond{A}{\neg B}$ must be true.  Using the propositional calculus,
both equations $(\fval \rightarrow B) = \tval$ and $(\fval \rightarrow
\neg B) = \tval$ described by Equation~\ref{eq:truth-functional} are
correct, because each of the included formulas of material implication
simplifies to the elementary truth-value $\tval$.

Conditionals of the remaining types are always false when their
antecedents are false.  With $\prob{A=\tval}=0$ the constraint
$\prob{A=\tval}>0$ required by Equation~\ref{eq:existential} for the
\rname{existential} conditional $\econd[k]{A}{B}$ is violated.  Also
with $\prob{A=\tval}=0$ the constraint $\prob{A=\tval}=1$ included in
the solution set from Equation~\ref{eq:feasibility} for the
\rname{feasibility} conditional $\fcond[k]{A}{B}$ (meaning
$\fcond[k]{\{\prob{A=\tval}=1\}}{B}$) is violated, leading to an empty
solution set which therefore cannot equal any set $\{k\}$.  Similarly
with $A=\fval$ the constraint $A=\tval$ in the solution set from
Equation~\ref{eq:boolean-feasibility} for the
\rname{Boolean-feasibility} conditional $\bcond[K]{A}{B}$ is violated,
leading to an empty solution set that cannot equal $\{\tval\}$ nor
$\{\fval\}$.




\section{Details: Factuality, Basic Measures, and Revision}

Let us now address a few details of the mathematical framework for
analyzing conditional statements which was just introduced: factuality
and counterfactuality; correctness and exception handling; and diverse
basic measures underlying probabilities.  There are several related
properties to appreciate about conditional statements.  We have
already considered \emph{mood}, using the distinction between
subjunctive and several kinds of indicative conditionals.  These other
properties are conceptually separate.  For example, a conditional
statement can be factual or counterfactual, independently of whether
it is subjunctive or indicative.  Moreover, a conditional could be
correct or incorrect, whatever combination of mood and factuality it
has.  What is believed to be correct or incorrect may change over time
as the knowledge of the analyst changes.  Finally, different
conditional statements, even those describing the same set of events
or propositions, may in fact concern measures of different basic
properties: observed frequency, subjective belief, causal propensity,
theoretical symmetry, and so on.  All of these features can be mixed
and matched in order to express precisely the semantic message that
the analyst wishes to express.

\subsection{Factuality and Counterfactuality}

Let us say that the `factuality' of a conditional statement pertains
to specific identified facts, and describes whether each fact or its
negation is included explicitly in the conditional statement.  We
recognize three states of factuality for a conditional statement
relative to a given fact: \emph{factuality} if the fact appears in the
conditional; \emph{antifactuality} if the negation of the given fact
appears in the conditional; and \emph{afactuality} if neither the fact
nor its negation appears in the conditional.  Afactuality and
antifactuality are our two kinds of counterfactuality; the former is
weaker and the latter is stronger.  A single conditional could be both
factual and antifactual with respect to some particular fact, if both
the fact and its negation were to appear in that conditional (though
legal, that would be an odd construction to make).  Although
factuality generally concerns antecedent part of the conditional in
question, we shall also include the consequent part in our
deliberations.

For example if the truth of some logical proposition $C$ were accepted
as a fact, then all of these conditionals would be factual with
respect to this fact $C=\tval$:
\begin{equation}
  \mcond{\{\overline{A},C\}}{B} \qquad
  \scond{C}{B} \qquad
  \tcond{C}{B} \qquad
  \bcond{\{A_1,\overline{A_2},C\}}{B} \qquad
  \fcond{\{A,C\}}{\overline{B}} \qquad
  \tcond{A}{C}
\end{equation}
However if the known fact were instead $A=\tval$, then relative to
this fact: the first conditional would be antifactual; the second,
third, and fourth conditionals would be afactual; and the fifth and
sixth conditionals would be factual.

As mentioned above, \emph{factuality} and \emph{mood} are two
different properties.  In principle any combination of \{factual,
afactual, antifactual\} and \{subjunctive, indicative\} can occur
within a single conditional statement.  However certain other features
of conditionals may correlate with these two properties.  For example,
an antifactual indicative conditional (whose antecedent contains the
negation of a known fact) would have to be correct if its type were
\rname{material} or \rname{truth-functional}, but incorrect if it had
any other type.  Anyway, it is important not to confuse
subjunctiveness with counterfactuality; they are two different
properties.

\subsection{Diverse Basic Measures}

Probabilities can be understood as ratios of basic measures (in the
measure-theoretic sense).  Basic measures may refer to many different
properties of objects or events.  Therefore, probabilistic
conditionals may refer to different properties of the propositions
that they concern.  Even \rname{truth-functional} and
\rname{Boolean-feasibility} conditionals can be understood to invoke
probabilities and hence diverse basic measures.  It can be very
important to understand which basic measure is described by a
particular conditional statement.

\begin{table}
  \sf
  \begin{tabular}{l|ccrrrcr} \hline
    \thd{Name} & 
    \thd{Cu-plating?} &
    \thd{Edge} &
    \thd{Value} &
    \thd{Mass} &
    \thd{2013 Production} &
    \thd{Lincolns} &
    \thd{Fraction Cu}
    \\ \hline\hline
    Penny & yes & smooth & 1~\textcent & 2.500~g &
    7070.00 M coins & 1~image & 2.50~\%
    \\ \hline
    Nickel & no & smooth & 5~\textcent & 5.000~g &
    1223.04 M coins & 0~images & 75.00~\%
    \\ \hline
    Dime & no & reeded & 10~\textcent & 2.268~g &
    2112.00 M coins & 0~images & 91.77~\%
    \\ \hline
    Quarter & no & reeded & 25~\textcent & 5.670~g &
    1455.20 M coins & 0~images & 91.77~\%
    \\ \hline
  \end{tabular}
  \caption{Selected properties of common United States coins, using
    data from \url{http://www.usmint.gov}.  Here `Cu' abbreviates the
    element copper and `M' stands for million.  Each penny has one
    image of President Lincoln on it; coins of the other denominations
    have none.  Reeds are the small ridges along the edges of some
    coins.}
  \label{tbl:coins}
\end{table}

For a simple illustration, consider certain properties of several
U.S. coins which are listed in Table~\ref{tbl:coins}.  Let us first
ask, what is the probability of a dime?  Using face monetary value as
the basic measure, we compute the following measures for the relevant
events:
\begin{eqnarray}
  \mu(\{\state{dime}\}) & \Rightarrow & 10\;\mbox{\textcent} \\
  \mu(\{\state{penny},\state{nickel},\state{dime},\state{quarter}\})
  & \Rightarrow & 41\;\mbox{\textcent}
\end{eqnarray}
The probability of the event 
\begin{math}
  \{\state{dime}\}
\end{math}
relative to the universal event
\begin{math}
  \{\state{penny},\state{nickel},\state{dime},\state{quarter}\}
\end{math}
(generally denoted $\Omega$) is given by the ratio of these measures:
\begin{math}
  10 \mbox{\textcent} / 41 \textrm{\textcent}
\end{math},
which is approximately $0.244$.  Different basic measures give
different numerical probabilities.  For example using mass as the
basic measure, $\prob{\{\state{dime}\}}$ evaluates to
$2.268\,\mathrm{g}/15.438\,\mathrm{g}$ which is approximately
$0.147$.  Or, using the cardinality of the set that constitutes each
event as its measure, the probability $\prob{\{\state{dime}\}}$, which
in fact designates the conditional probability
\begin{math}
  \condprob{\{\state{dime}\}}
           {\{\state{penny},\state{nickel},\state{dime},\state{quarter}\}}
\end{math},
evaluates to $1/4$.

Note a few important features of probabilities: probabilities are
inherently conditional (the universal set $\Omega$ containing all
elementary events under consideration gives the default denominator);
basic measures have units of measure (cents, grams, and so on), which
cancel out during division; measures are assigned to \emph{sets} of
elementary events; measures are by definition additive with the union
of disjoint sets; the empty set must have measure zero; other sets of
elementary events may also have measure zero.  Note also that not
every numerical property behaves like a measure.  For example the
fraction of copper is not additive with set union: although a dime
contains $91.77\%$ copper and a quarter also contains $91.77\%$
copper, it is not the case that the set containing both a dime and a
quarter itself contains $183.54\%$ copper.  Moving on, in the special
case that $\mu(\Omega)=0$ all probabilities for the system in question
have the indefinite value $0/0$.  Uncountable or countably infinite
sets of elementary events require special attention as well.

Next let us ask, what is the conditional probability of copper plating
(denoted $B$), given smooth edge (denoted $A$)?  Using monetary value
as the basic measure gives $\condprob{B=\tval}{A=\tval} \Rightarrow
1/6$.  Using Lincoln images as the basic measure gives
$\condprob{B=\tval}{A=\tval} \Rightarrow 1$.  Returning our attention
to conditionals, let us observe that the affirmative \rname{material}
conditional $\mcond{A}{B}$ (defined by Equation~\ref{eq:material} with
$k=1$) is correct using Lincoln images as the basic measure, but
incorrect using monetary value as the basic measure.  The choice of
basic measure can affect the truth or falsity of a conditional
statement.

Note that for some very important basic measures we do not know quite
what the appropriate units of measure are.  For example, what is the
unit of subjective belief?  Or of causal propensity?  One potential
basic measure for historical events is an indicator-type measure of
what actually happened: we might assign one (anonymous) unit of
measure to the event that actually happened, and zero units to every
event that did not happen.  Here we see one benefit of the
construction of probability: it is possible to reason about ratios of
measures whose units we do not know, since those units cancel out
during the division used to compute the ratios.  We must take care to
recognize the potential for division by zero, and to handle that
exception appropriately.

We come to an important point.  Different types of conditionals tend
to use different basic measures.  \rname{Subjunctive} conditionals
tend to use basic measures such as subjective belief, symmetry
(symmetrical cases, as in probability $1/2$ for a coin to land on
heads because it has two symmetrical sides), or causal propensity.  In
contrast, the various indicative conditionals tend to use basic
measures such as occurrence (e.g.\ the indicator just discussed),
observed cases (whose ratio gives frequency), or physical properties
such as mass.

\subsection{Correctness and Revision for Exceptions}

Any given conditional statement may be correct or incorrect (true or
false).  In some cases that correctness or incorrectness derives
purely from the syntax of the statement.  For example, the
\rname{Boolean-feasibility} conditional $\bcond{A}{\overline{A}}$,
which says `Assuming that $A$ is true then $A$ must be false'
according to Equation~\ref{eq:boolean-feasibility}, is intrinsically
incorrect without regard to which logical proposition or real-world
event the variable $A$ represents.\footnote{You may be tempted to
  insert a self-reference here.  That would be a fine way to formulate
  a recurrence relation that defines a discrete dynamical system,
  amenable to algebraic analysis as discussed in
  \cite{norman-orthodoxy}.}  In other cases the correctness or
incorrectness of a conditional statement depends upon the referents of
its variables.  For example, the \rname{Boolean-feasibility}
conditional $\bcond{A}{B}$ might or might not be true, depending on
which propositions or events the variables $A$ and $B$ represent.

Critically, an analyst's impression about the correctness of any given
conditional statement may change over time as that analyst's knowledge
changes.  It is appropriate that our models should change as our
knowledge does.  It is important to allow some mechanisms for revision
or exception handling.  Thus we might supplement the statistician
Box's observation that ``All models are wrong, but some are useful''
\cite{box} with a paraphrase of Shakespeare: Some models are born
wrong, and some have wrongness thrust upon them.

Considering the diverse basic measures above, it is natural that some
of these measures should change over time.  The physical universe
changes.  And the subjective beliefs of an analyst might change too.
Here it is helpful to think in terms of absolute units of belief
rather than ratios of them.  A common motif in problems is learning
about an exception to a rule.  This may seem tricky in terms of
probabilities, but it is more straightforward in terms of basic
measures.  

Also note that it may not be possible or practical to verify whether
any given conditional statement is true or false, because the
underlying basic measure cannot practically be measured.  We have as
yet no instruments to test subjective belief.  It is a matter of
scientific experimentation and statistical analysis to establish
measures of causal propensity.  Even symmetry may take some nontrivial
investigation to understand.

\subsubsection{Revision of Old Conditionals}

Our models may change as our knowledge does.  Consider two options
which may be useful.  First, have the ability to retract an old
conditional statement upon learning about an exception (and then
assert a new one).  This requires nonmonotonicity.  Second, formulate
(some or all) conditionals as defeasible statements in the first place
(using nullifier terms in the antecedents), so that any statement can
be rendered ineffectual without being retracted (`defeated without
being deleted').  Now we need only the ability to add fixed values for
propositional variables at run time (and to answer queries conditioned
on the values of the nullifier terms, so that they can be ignored when
it is desired to consider the default case).

\subsubsection{Defeasible Conditionals}

One option is to include anonymous `nullifier' terms to represent
otherwise unexpected exceptions.  For example instead of
$\bcond{A}{B}$ we could use $\bcond{(A \wedge \neg Z)}{B}$ to say `If
$A$ then $B$, unless $Z$'.  Another option is to use probabilities
other than zero and one, for example always $0.001$ and $0.999$
instead (or a fancier solution such as $\delta$ instead of $0$, and
$1-\delta$ instead of $1$, with $0 \le \delta \le 0.001$).

\section{Algebraic Deduction}

A conditional of any principal type can be interpreted as a statement
of logical deduction.  Among the various types, \rname{Feasibility}
and \rname{Boolean-feasibility} conditionals are especially
well-suited for this purpose: they reproduce intuitive patterns of
inference such as \emph{modus ponens}; they avoid explosions of
consequences from inconsistent premises; and they provide frameworks
within which conditionals of other types may be included.  In this
section we shall investigate how certain methods of algebra can be
used to compute logical deductions involving probabilistic,
truth-functional, and Boolean conditionals.  As we shall see, the
solutions to various systems of polynomial equalities and inequalities
provide interesting results.

\subsection{Probabilistic Deduction}

Probabilistic \rname{feasibility} conditionals can be interpreted as
statements of logical deduction.  Using this interpretation we can
read the following affirmative \rname{feasibility} conditional from
Equation~\ref{eq:feasibility} (with $k=1$):
\begin{eqnarray}
  \fcond{\Gamma}{B} & \equiv &
  \left\{ \: \prob{B=\tval} \: : \: \Gamma \: \right\}
  = \{1\}
  \label{eq:prob-conseq}
\end{eqnarray}
as the deductive statement that `$B$ is a consequence of $\Gamma$'.
In this interpretation we consider the members of the antecedent set
$\Gamma$ to be premises, and the consequent $B$ to be the conclusion.
Recalling the earlier description this deductive statement also says
that `Subject to the constraints in $\Gamma$, the probability that $B$
is true must be exactly $1$'.  As with all conditionals, any
particular deductive statement $\fcond{\Gamma}{B}$ itself may be true
or false.  The methods of Section \ref{sec:spi} and
\ref{sec:optimization} allow us to compute which is the case.  Let us
say that this formulation provides one of many possible notions of
consequentiality involving logical propositions.

For example, to demonstrate \emph{modus ponens} using the probability
model from Equation~\ref{eq:ab-input}, let us specify as the
antecedent set $\Gamma$ the constraint $\prob{A=\tval}$ that $A$ is
certainly true, and some version of the affirmative conditional `If
$A$ then $B$'.  Using the \rname{material} version of this affirmative
conditional gives the constraint
$\prob{A=\tval,B=\tval}=\prob{A=\tval}$ from
Equation~\ref{eq:material}.  We identify $\prob{B=\tval}$ as the
objective function.  The resulting solution set is given by:
\begin{eqnarray}
  \Phi & \Leftarrow &
  \left\{ \:
  \prob{B=\tval} \: : \:
  \prob{A=\tval}=1, \:
  \prob{A=\tval,B=\tval}=\prob{A=\tval}
  \: \right\}
\end{eqnarray}
After symbolic probability inference, this set-comprehension
expression becomes precisely one presented in
Equation~\ref{eq:b-solset}.  As demonstrated in
Section~\ref{sec:optimization} the computed minimum solution
\begin{math}
  \alpha^* =
1.000
\end{math}
and the maximum solution
\begin{math}
  \beta^* =
1.000
\end{math}.
By these calculations the solution set $\Phi \Rightarrow \{1\}$.  In
other words, we have computed that the following equation about sets
of real numbers is satisfied:
\begin{eqnarray}
  \left\{ \:
  \prob{B=\tval} \: : \:
  \prob{A=\tval}=1, \:
  \prob{A=\tval,B=\tval}=\prob{A=\tval}
  \: \right\}
  & = & \{1\}
\end{eqnarray}
It follows that the equivalent \rname{feasibility} conditional from
Equation~\ref{eq:prob-conseq} is true:
\begin{equation}
  \fcond{\{ \: A, \: \mcond{A}{B} \: \}}{B}
\end{equation}
Therefore it is correct to say that, assuming the premises that $A$ is
certainly true and that the affirmative \rname{material} conditional
`If $A$ then $B$' holds, it is a consequence that $B$ must certainly
be true.  For this example, it happens that substituting either the
\rname{subjunctive} or the \rname{existential} interpretation of the
inner conditional `If $A$ then $B$' would also give a true outer
\rname{feasibility} conditional; both of the following statements are
correct:
\begin{equation}
  \begin{array}{r}
    \fcond{\{ \: A, \: \scond{A}{B} \: \}}{B} \\
    \fcond{\{ \: A, \: \econd{A}{B} \: \}}{B}
  \end{array}
\end{equation}

\subsubsection{Probabilistic Conditionals as Consequences}
\label{sec:conditionals-consequences}

We can use probability equations similar to the \rname{feasibility}
conditional in Equation~\ref{eq:feasibility} in order to determine
whether or not certain proposed \rname{subjunctive}, \rname{material},
or \rname{existential} conditionals are correct subject to given sets
of constraints.  Let us focus on affirmative conditionals (with
fractional parameter $k=1$).  Consider the following solution sets
which make use of a given set $\Gamma$ of input constraints, with
certain probability expressions involving the primary variables $A$
and $B$ as objective functions and additional constraints:
\begin{eqnarray}
  \Phi & \Leftarrow &
  \{ \: \condprob[0]{B=\tval}{A=\tval} \: : \: \Gamma \: \}
  \label{eq:c-phi}
  \\
  \Psi & \Leftarrow &
  \{ \: \prob{A=\tval}-\prob{A=\tval,B=\tval} \: : \: \Gamma \: \}
  \label{eq:c-psi}
  \\
  \Upsilon & \Leftarrow &
  \{ \: \prob{A=\tval} \: : \: 
  \prob{A=\tval}=\prob{A=\tval,B=\tval}, \:  
  \Gamma \: \}
  \label{eq:c-upsilon}
\end{eqnarray}
If (and only if) the input constraints $\Gamma$ are inconsistent (when
considered together with the default probability constraints
$\Gamma_0$), then all three solution sets $\Phi$, $\Psi$, and
$\Upsilon$ will evaluate to the empty set.  Otherwise the values of
these solution sets reveal which probabilistic conditional statements
are necessarily correct and which are possibly correct.

The affirmative \rname{subjunctive} conditional $\scond{A}{B}$ defined
by Equation~\ref{eq:subjunctive} requires
$\condprob[0]{B=\tval}{A=\tval}=1$, hence $\Phi=\{1\}$ in terms of
Equation~\ref{eq:c-phi}.  Thus it is \emph{necessary} that
$\scond{A}{B}$ is correct if the computed value of $\Phi$ is precisely
the set $\{1\}$.  If $\Phi$ merely contains $1$ then it is
\emph{possible} that $\scond{A}{B}$ is correct.  If $\Phi$ does not
contain $1$ then it is \emph{impossible} that $\scond{A}{B}$ is
correct.  Likewise, the affirmative \rname{material} conditional
$\mcond{A}{B}$ defined by Equation~\ref{eq:material} requires
$\prob{A=\tval}=\prob{A=\tval,B=\tval}$, hence $\Psi=\{0\}$ in terms
of Equation~\ref{eq:c-psi}.  It follows that $\mcond{A}{B}$ is
necessarily correct if $\Psi=\{0\}$; possibly correct if $0 \in \Psi$;
and impossibly correct if $0 \notin \Psi$.  The \rname{existential}
conditional $\econd{A}{B}$ from Equation~\ref{eq:existential} adds to
its \rname{material} counterpart the requirement $\prob{A=\tval}>0$.
In terms of Equation~\ref{eq:c-upsilon} this means the solution set
$\Upsilon$ must contain only nonzero values: hence $\Upsilon \neq
\emptyset$ and $0 \notin \Upsilon$.  It follows that $\econd{A}{B}$ is
necessarily correct if $\Upsilon \neq \emptyset$ and $0 \notin
\Upsilon$; possibly correct if $\Upsilon \neq \emptyset$ and $\Upsilon
\neq \{0\}$; and impossibly correct if $\Upsilon=\emptyset$ or
$\Upsilon=\{0\}$.


\subsection{Boolean Deduction}
\label{sec:boolean-deduction}

\rname{Boolean feasibility} conditionals can be interpreted as
statements of logical deduction.  This is in effect what Boole did
this throughout his \emph{Laws of Thought} by drawing logical
conclusions from the solutions to polynomial equations \cite{boole}.
Considering an antecedent $A$ and a consequent $B$ which are formulas
of the propositional calculus (using some propositional variables
$X_1,\ldots,X_n$), let us read the affirmative conditional
$\bcond{A}{B}$ from Equation~\ref{eq:boolean-feasibility} as the
deductive statement `$B$ is a consequence of $A$':
\begin{eqnarray}
  \bcond{A}{B} & \equiv &
  \left\{ \: B \: : \: A=\tval \: \right\}
  = \{\tval\}
  \label{eq:boole-conseq}
\end{eqnarray}
Recalling the earlier description, this statement $\bcond{A}{B}$ also
says that `Assuming that $A$ is true, then $B$ must be true'.  After
Boolean polynomial translation, this definition becomes:
\begin{eqnarray}
  \bcond{A}{B} & \equiv &
  \left\{ \: \bpoly{B} \: : \: \bpoly{A}=1; \: 
  x_1\in\{0,1\},\ldots,x_n\in\{0,1\}
  \: \right\}
  = \{1\}
\end{eqnarray}
The set $\Gamma_0$ of constraints expanded into the above definition
includes a constraint $x_i\in\{0,1\}$ for each real-valued variable
$x_i$ translated from a propositional variable $X_i$.  Note that the
affirmative \rname{Boolean-feasibility} conditional $\bcond{A}{B}$ in
Equation~\ref{eq:boole-conseq} is not the same as the affirmative
\rname{truth-functional} conditional from
Equation~\ref{eq:truth-functional}.  In particular, when the
antecedent formula $A$ is unsatisfiable (identically false) then the
relation $\tcond{A}{B}$ is true but the relation $\bcond{A}{B}$ is
not.

We could also allow modal statements of consequentiality, such as `$B$
is a potential consequence of $A$' meaning
$\tval\in\{\:B\::\:A=\tval\:\}$, thus satisfied by either solution set
$\{\tval\}$ or $\{\tval,\fval\}$ but not by $\{\fval\}$ or $\emptyset$.

Boolean deduction as defined by Equation~\ref{eq:boole-conseq}
reproduces the familiar pattern of \emph{modus ponens}.  Considering
propositional variables $X$ and $Y$, let us take as the antecedent $A$
the assertion that the formulas $X$ and $X \rightarrow Y$ are both
true: hence the conjunction $X \wedge (X \rightarrow Y)$.  Does it
follow that $Y$ must also be true?  Using the formula $Y$ as the
consequent $B$, Equation~\ref{eq:boole-conseq} says that the
Boolean-deductive statement $\bcond{X \wedge (X \rightarrow Y)}{Y}$
signifies the following equation about sets of truth values:
\begin{eqnarray}
  \left\{ \: X \: : \: X \wedge (X \rightarrow Y) = \tval; \: 
  \: \right\}
  & = & \{\tval\}
\end{eqnarray}
After Boolean translation according to Section~\ref{sec:boolean}, this
becomes an equation about sets of real numbers:
\begin{eqnarray}
  \left\{ \: x \: : \: xy = 1; \: 
  x \in \{0,1\}, \: y \in \{0,1\}
  \: \right\}
  & = & \{1\}
\end{eqnarray}
Section~\ref{sec:boolean} showed the Boolean polynomial translation
\begin{math}
  \bpoly{X \wedge (X \rightarrow Y)} \Rightarrow xy
\end{math}.
Algebraic analysis confirms that the value of the set-comprehension
expression above is indeed the set $\{1\}$.  The corresponding
optimization problems:
\begin{equation}
\begin{array}{r@{\quad}l}
\mbox{Minimize}: & y \\
\mbox{subject to}:
& x y = 1 \\
\mbox{and}:
& x \in \{{0, 1}\} \\
& y \in \{{0, 1}\} \\

\end{array}
\qquad
\begin{array}{r@{\quad}l}
\mbox{Maximize}: & y \\
\mbox{subject to}:
& x y = 1 \\
\mbox{and}:
& x \in \{{0, 1}\} \\
& y \in \{{0, 1}\} \\

\end{array}
\end{equation}
have computed solutions
\begin{math}
  \alpha^* =
1.000
\end{math}
and
\begin{math}
  \beta^* =
1.000
\end{math}.
Thus it is demonstrated that the statement of Boolean
deduction
\begin{math}
  \bcond{X \wedge (X \rightarrow Y)}{Y}
\end{math},
interpreted according to Equation~\ref{eq:boole-conseq}, is correct.
This statement says that, assuming the premise that the formula $X
\wedge (X \rightarrow Y)$ is true, then it is a consequence that the
formula $Y$ must be true also.

Logical deduction based on \rname{Boolean-feasibility} conditionals
exhibits several desirable properties.  This methodology provides a
paraconsistent system of logic (for formulas of the propositional
calculus) that avoids explosion from inconsistent premises, yet
retains the intuitive properties of disjunction introduction,
disjunctive syllogism, and transitivity.

\subsubsection{Boolean Consequences and Sets of Valuations}

There is an important difference between the \rname{truth-functional}
conditionals from Section~\ref{sec:truth-functional} and the
\rname{Boolean-feasibility} conditionals from
Section~\ref{sec:boolean-feasibility}.  Whereas each
\rname{truth-functional} conditional operates on individual valuations
of the propositional variables in use, each
\rname{Boolean-feasibility} conditional operates on set of such
valuations.  Understanding this feature we can develop algebraic
problems to determine when proposed \rname{Boolean-feasibility}
conditionals are true subject to given constraints.

Considering a system with some list $(X_1,\ldots,X_n)$ of
propositional variables, let us say that a \emph{valuation} or
\emph{valuation vector} consists of an assignment of a truth value
$K_i\in\{\tval,\fval\}$ to each propositional variable $X_i$; thus a
vector $(K_1,\ldots,K_n)$ of truth values.  With $n$ propositional
variables there are $2^n$ possible valuation vectors.  For example
with $2$ propositional variables $(A,B)$ there are $4$ possible
valuations $(\tval,\tval)$, $(\tval,\fval)$, $(\fval,\tval)$, and
$(\fval,\fval)$.  The affirmative \rname{truth-functional} conditional
$\tcond{A}{B}$ from Equation~\ref{eq:truth-functional} operates on one
valuation at a time.  For example with $(A,B)=(\tval,\tval)$ the
conditional $\tcond{A}{B}$ is true and with $(A,B)=(\tval,\fval)$ the
conditional $\tcond{A}{B}$ is false.  A logical truth table summarizes
such results.  In this case:
\begin{equation}
  \begin{array}{cc|c} \hline
    A & B & \tcond{A}{B} \\ \hline\hline
    \tval & \tval & \tval \\ \hline
    \tval & \fval & \fval \\ \hline
    \fval & \tval & \tval \\ \hline
    \fval & \fval & \tval \\ \hline
  \end{array}
\end{equation}

In contrast, the affirmative \rname{Boolean-feasibility} conditional
from Equations \ref{eq:boolean-feasibility} and \ref{eq:boole-conseq}
properly operates on one \emph{set} of valuations at a time.  For
example, given the valuation $(A,B)=(\tval,\tval)$, is the affirmative
\rname{Boolean-feasibility} conditional $\bcond{A}{B}$ true or false?
Well, that depends on which other valuations are possible.  If the
valuation $(A,B)=(\tval,\fval)$ were also possible then the solution
set $\{B:A=\tval\}$ from Equation~\ref{eq:boolean-feasibility} would
contain $\fval$ as well as $\tval$, hence the solution set
$\{\tval,\fval\}$.  This solution set $\{\tval,\fval\}$ does not equal
the set $\{\tval\}$ required by the affirmative conditional; thus
$\bcond{A}{B}$ would be false.  On the other hand, if the valuation
$(A,B)=(\tval,\fval)$ were impossible, then the solution set
$\{B:A=\tval\}$ from Equation~\ref{eq:boolean-feasibility} would
contain only the value $\tval$, hence the solution set $\{\tval\}$
which would make the affirmative conditional $\bcond{A}{B}$ true.

We can use a modified truth table to describe the sets of variable
valuations that make a given \rname{Boolean-feasibility} conditional
true.  To this end let us evaluate the set comprehension expression
for the conditional $\bcond{A}{B}$ from
Equation~\ref{eq:boolean-feasibility} four separate times, using as
input the singleton set $\{(K_1,K_2)\}$ tabulated by each row of the
table:
\begin{equation}
  \begin{array}{cc|c} \hline
    A & B & \{B:A=\tval\} \\ \hline\hline
    \tval & \tval & \{\tval\} \\ \hline
    \tval & \fval & \{\fval\} \\ \hline
    \fval & \tval & \emptyset \\ \hline
    \fval & \fval & \emptyset \\ \hline
  \end{array}
\end{equation}
The last two sets are empty because the constraint $A=\tval$ is
violated.  These values give instructions for constructing the sets of
valuations of $(A,B)$ that satisfy the conditional $\bcond{A}{B}$:
each set must include the point $(A,B)=(\tval,\tval)$; it must not
include the point $(\tval,\fval)$; and it may or may not include the
points $(\fval,\tval)$ and $(\fval,\fval)$.  These instructions are
like the four-part set comprehension expressions that Boole used in
\cite{boole}.  However in this context, in general the mandatory part
of the set of valuations is specified as one or more members of an
identified subset.  Anyway for this example, $4$ of the $16$ possible
sets of valuations match the instructions:
\begin{equation}
  \{ \: (\tval,\tval) \: \}, \qquad
  \{ \: (\tval,\tval), \: (\fval,\tval) \: \}, \qquad
  \{ \: (\tval,\tval), \: (\fval,\fval) \: \}, \qquad
  \{ \: (\tval,\tval), \: (\fval,\tval), \: (\fval,\fval) \: \}
\end{equation}
For each of these four sets of valuation vectors $(A,B)$ it is true
that the set of values of $B$ subject to the constraint $A=\tval$
equals the set $\{\tval\}$.  Thus each set of valuations satisfies the
relation $\{B:A=\tval\}=\{1\}$ that defines the affirmative
\rname{Boolean-feasibility} conditional $\bcond{A}{B}$ according to
Equation~\ref{eq:boolean-feasibility}.  Note that only one of these
four sets of valuations consists of the valuations that satisfy the
material-implication formula $A \rightarrow B$.

\subsubsection{Probabilistic Expression of Boolean Feasibility}
\label{sec:prob-boole}

There is a different way to specify the set of sets of valuations that
satisfy a given \rname{Boolean-feasibility} conditional.  Let us adopt
a new probability model in which the joint probability table
$\prob[0]{A,B}$ is input directly:
\begin{equation}
\begin{tabular}[c]{ll|l}\hline
\multicolumn{1}{l}{$A$} & 
\multicolumn{1}{l|}{$B$} & 
\multicolumn{1}{l}{$\prob[0]{{A, B}}$} \\ \hline\hline
$\state{T}$ & 
$\state{T}$ & 
$x_{1}$ \\ \hline
$\state{T}$ & 
$\state{F}$ & 
$x_{2}$ \\ \hline
$\state{F}$ & 
$\state{T}$ & 
$x_{3}$ \\ \hline
$\state{F}$ & 
$\state{F}$ & 
$x_{4}$ \\ \hline
\end{tabular}
  \label{eq:ab-linear}
\end{equation}
The constraints $0 \le x_i \le 1$ and
\begin{math}
x_{1}
  +
x_{2}
  +
x_{3}
  +
x_{4}
  = 1
\end{math}
are added to enforce the laws of probability.  We embed the
propositional-calculus formulas $A \rightarrow B$ and $A \wedge B$ by
the method of Section~\ref{sec:embedding}:
\begin{equation}
\begin{tabular}[c]{ll|ll}\hline
\multicolumn{4}{l}{$\condprob[0]{{\embed{A \wedge \neg B}}}{{A, B}}$} \\ \hline\hline
\multicolumn{1}{l}{$A$} & 
\multicolumn{1}{l|}{$B$} & 
\multicolumn{1}{l}{${\embed{A \wedge \neg B}=\state{T}}$} & 
\multicolumn{1}{l}{${\embed{A \wedge \neg B}=\state{F}}$} \\ \hline\hline
$\state{T}$ & 
$\state{T}$ & 
$0$ & 
$1$ \\ \hline
$\state{T}$ & 
$\state{F}$ & 
$1$ & 
$0$ \\ \hline
$\state{F}$ & 
$\state{T}$ & 
$0$ & 
$1$ \\ \hline
$\state{F}$ & 
$\state{F}$ & 
$0$ & 
$1$ \\ \hline
\end{tabular}
  \qquad
\begin{tabular}[c]{ll|ll}\hline
\multicolumn{4}{l}{$\condprob[0]{{\embed{A \wedge B}}}{{A, B}}$} \\ \hline\hline
\multicolumn{1}{l}{$A$} & 
\multicolumn{1}{l|}{$B$} & 
\multicolumn{1}{l}{${\embed{A \wedge B}=\state{T}}$} & 
\multicolumn{1}{l}{${\embed{A \wedge B}=\state{F}}$} \\ \hline\hline
$\state{T}$ & 
$\state{T}$ & 
$1$ & 
$0$ \\ \hline
$\state{T}$ & 
$\state{F}$ & 
$0$ & 
$1$ \\ \hline
$\state{F}$ & 
$\state{T}$ & 
$0$ & 
$1$ \\ \hline
$\state{F}$ & 
$\state{F}$ & 
$0$ & 
$1$ \\ \hline
\end{tabular}
\end{equation}
Parametric probability analysis gives the following results for
$\prob{\embed{A \rightarrow B}}$ and $\prob{\embed{A \wedge B}}$:
\begin{equation}
\begin{tabular}[c]{l|l}\hline
\multicolumn{1}{l|}{$\embed{A \wedge \neg B}$} & 
\multicolumn{1}{l}{$\prob{{\embed{A \wedge \neg B}}}$} \\ \hline\hline
$\state{T}$ & 
$x_{2}$ \\ \hline
$\state{F}$ & 
$x_{1} + x_{3} + x_{4}$ \\ \hline
\end{tabular}
  \qquad
\begin{tabular}[c]{l|l}\hline
\multicolumn{1}{l|}{$\embed{A \wedge B}$} & 
\multicolumn{1}{l}{$\prob{{\embed{A \wedge B}}}$} \\ \hline\hline
$\state{T}$ & 
$x_{1}$ \\ \hline
$\state{F}$ & 
$x_{2} + x_{3} + x_{4}$ \\ \hline
\end{tabular}
\end{equation}
The output table $\prob{A,B}$ is identical to the input table
$\prob[0]{A,B}$.
It happens in this model that all polynomials are linear functions of
the parameters
\begin{math}
x_{1}
  , \ldots,
x_{4}
\end{math};
hence conditional probability queries yield fractional linear
functions.

For the affirmative \rname{Boolean-feasibility} conditional
$\bcond{A}{B}$ we can define the set of satisfactory sets of
valuations of $(A,B)$ by two probability constraints involving the
embedded logical formulas:
\begin{eqnarray}
  \prob{\embed{A \wedge \neg B}=\tval} & = & 0
  \label{eq:animplb}
  \\
  \prob{\embed{A \wedge B}=\tval} & > & 0 
  \label{eq:aandb}
\end{eqnarray}
Let us call these constraints the `negative' and `positive' criteria,
respectively.  Regarding the probability model defined by the input
table $\prob[0]{A,B}$ above, you can see that these two equations are
necessary and sufficient conditions for the probabilistic
\rname{feasibility} conditional $\fcond{A}{B}$ to hold.  Note that the
embedded event $\embed{A \wedge \neg B}=\tval$ is equivalent to the
embedded events $\embed{A \rightarrow B}=\fval$ and $\embed{A
  \nrightarrow B}=\tval$.  The material nonimplication relation $A
\nrightarrow B$ means the negation $\neg(A \rightarrow B)$.

We consider a set of valuations to be satisfactory if each of its
members (corresponding to an element of the table $\prob{A,B}$) can be
assigned a nonzero probability such that these constraints are
satisfied.  After symbolic probability inference these become the
algebraic constraints
\begin{math}
x_{2}
  = 0
\end{math}
and
\begin{math}
x_{1}
  > 0
\end{math}.
These relations are useful because they can be used as constraints or
as objectives in optimization problems, in order to assert or test
whether the conditional $\bcond{A}{B}$ holds.  For example to query
whether the conditional $\bcond{A}{B}$ is true given some set $\Gamma$
of logical constraints, we can evaluate the solution sets:
\begin{eqnarray}
  \Phi & \Leftarrow &
  \{ \: \prob{\embed{A \wedge \neg B}=\tval} \: : \: \Gamma \: \}
  \\
  \Psi & \Leftarrow &
  \{ \: \prob{\embed{A \wedge B}=\tval} \: : \:
  \prob{\embed{A \wedge \neg B}=\tval}=0, \: \Gamma \: \}
\end{eqnarray}
If $\Phi=\{0\}$ and $\Psi \neq \emptyset$ and $0 \notin \Psi$ then
given $\Gamma$ it is necessary that $\bcond{A}{B}$ is true; if $0 \in
\Phi$ and $\Psi \neq \emptyset$ and $\Psi \neq \{0\}$ then given
$\Gamma$ it is possible that $\bcond{A}{B}$ is true; and if $0 \notin
\Phi$ or $\Psi=\{0\}$ or $\Psi=\emptyset$ then given $\Gamma$ it is
impossible that $\bcond{A}{B}$ is true.  When the full-joint
probability distribution concerning the propositional variables has
been specified using a single input table as $\prob[0]{A,B}$ above,
the desired solutions can be computed by linear optimization (there
are no nonlinear polynomials involved, as is the general case for
parametric probability models).

Note that there are several ways for a set of constraints $\Gamma$ to
fail to confirm a conditional $\bcond{A}{B}$.  It could be that the
constraints $\Gamma$ are inconsistent, in which case we do not
consider any conditional statement to hold.  Next there is the `sin of
commission' that it is possible (probability greater than zero) for
the consequent $B$ to be false while the antecedent $A$ is true.
Finally there is the `sin of omission' that it is impossible
(probability zero) for the consequent $B$ to be true while the
antecedent $A$ is true.  We can consider these mechanisms of failure
in a modal way by the above analysis (thus computing whether it is
necessary, possible, or impossible for $\bcond{A}{B}$ to hold subject
to the constraints $\Gamma$).

\subsubsection{Disjunction Introduction and Disjunctive Syllogism}

Logical deduction based on \rname{Boolean-feasibility} conditionals
provides disjunction introduction and disjunctive syllogism.  To check
the former we evaluate the conditional $\bcond{A}{A \vee B}$; for the
latter we evaluate $\bcond{\{\:A \vee B,\:\neg A\:\}}{B}$.  We clarify
that restrictions on the joint prior probabilities of logical
propositions may invalidate these or other deductions.  For example if
$A=B$ then disjunctive syllogism fails.  To be explicit about prior
probabilities we could use probabilistic \rname{feasibility}
conditionals instead of \rname{Boolean-feasibility} conditionals, with
the aid of embedded propositional-calculus functions to express the
desired premises.

\subsubsection{Transitivity of Boolean Deduction}
\label{sec:transitivity}

To illustrate transitivity, let us investigate the relationship between the
\rname{Boolean-feasibility} conditionals $\bcond{A}{B}$ and
$\bcond{B}{C}$ asserted as premises, and the conditional
$\bcond{A}{C}$ queried as an objective.  We desire to evaluate the
solution set:
\begin{eqnarray}
  \Phi_0 & \Leftarrow &
  \{ \: (\bcond{A}{C}) \: : \: 
  (\bcond{A}{B})=\tval, \: (\bcond{B}{C})=\tval \: \}
\end{eqnarray}
If the computed value of $\Phi_0$ is $\{\tval\}$ then we will declare
$\bcond{A}{C}$ a necessary consequence of the premises $\bcond{A}{B}$
and $\bcond{B}{C}$, according to Equation~\ref{eq:boole-conseq}.  It
turns out that we can formulate a series of linear optimization
problems whose solutions will reveal the computed value of the
truth-value set $\Phi_0$.

Consider a new probability model with the following input table
$\prob[0]{A,B,C}$, including parameter constraints
\begin{math}
  \Gamma_0 \Leftarrow 
  \left\{ \: 0 \le y_i \le 1, \; \sum_i y_i = 1 \: \right\}
\end{math}
to enforce the laws of probability:
\begin{equation}
\begin{tabular}[c]{lll|l}\hline
\multicolumn{1}{l}{$A$} & 
\multicolumn{1}{l}{$B$} & 
\multicolumn{1}{l|}{$C$} & 
\multicolumn{1}{l}{$\prob[0]{{A, B, C}}$} \\ \hline\hline
$\state{T}$ & 
$\state{T}$ & 
$\state{T}$ & 
$y_{1}$ \\ \hline
$\state{T}$ & 
$\state{T}$ & 
$\state{F}$ & 
$y_{2}$ \\ \hline
$\state{T}$ & 
$\state{F}$ & 
$\state{T}$ & 
$y_{3}$ \\ \hline
$\state{T}$ & 
$\state{F}$ & 
$\state{F}$ & 
$y_{4}$ \\ \hline
$\state{F}$ & 
$\state{T}$ & 
$\state{T}$ & 
$y_{5}$ \\ \hline
$\state{F}$ & 
$\state{T}$ & 
$\state{F}$ & 
$y_{6}$ \\ \hline
$\state{F}$ & 
$\state{F}$ & 
$\state{T}$ & 
$y_{7}$ \\ \hline
$\state{F}$ & 
$\state{F}$ & 
$\state{F}$ & 
$y_{8}$ \\ \hline
\end{tabular}
\end{equation}
Symbolic probability inference with embedded propositional-calculus
formulas yields:
\begin{eqnarray}
\begin{tabular}[c]{l|l}\hline
\multicolumn{1}{l|}{$\embed{A \wedge B}$} & 
\multicolumn{1}{l}{$\prob{{\embed{A \wedge B}}}$} \\ \hline\hline
$\state{T}$ & 
$y_{1} + y_{2}$ \\ \hline
$\state{F}$ & 
$y_{3} + y_{4} + y_{5} + y_{6} + y_{7} + y_{8}$ \\ \hline
\end{tabular}
  & &
\begin{tabular}[c]{l|l}\hline
\multicolumn{1}{l|}{$\embed{A \rightarrow B}$} & 
\multicolumn{1}{l}{$\prob{{\embed{A \rightarrow B}}}$} \\ \hline\hline
$\state{T}$ & 
$y_{1} + y_{2} + y_{5} + y_{6} + y_{7} + y_{8}$ \\ \hline
$\state{F}$ & 
$y_{3} + y_{4}$ \\ \hline
\end{tabular}
  \\
\begin{tabular}[c]{l|l}\hline
\multicolumn{1}{l|}{$\embed{B \wedge C}$} & 
\multicolumn{1}{l}{$\prob{{\embed{B \wedge C}}}$} \\ \hline\hline
$\state{T}$ & 
$y_{1} + y_{5}$ \\ \hline
$\state{F}$ & 
$y_{2} + y_{3} + y_{4} + y_{6} + y_{7} + y_{8}$ \\ \hline
\end{tabular}
  & &
\begin{tabular}[c]{l|l}\hline
\multicolumn{1}{l|}{$\embed{B \rightarrow C}$} & 
\multicolumn{1}{l}{$\prob{{\embed{B \rightarrow C}}}$} \\ \hline\hline
$\state{T}$ & 
$y_{1} + y_{3} + y_{4} + y_{5} + y_{7} + y_{8}$ \\ \hline
$\state{F}$ & 
$y_{2} + y_{6}$ \\ \hline
\end{tabular}
  \\
\begin{tabular}[c]{l|l}\hline
\multicolumn{1}{l|}{$\embed{A \wedge C}$} & 
\multicolumn{1}{l}{$\prob{{\embed{A \wedge C}}}$} \\ \hline\hline
$\state{T}$ & 
$y_{1} + y_{3}$ \\ \hline
$\state{F}$ & 
$y_{2} + y_{4} + y_{5} + y_{6} + y_{7} + y_{8}$ \\ \hline
\end{tabular}
  & &
\begin{tabular}[c]{l|l}\hline
\multicolumn{1}{l|}{$\embed{A \rightarrow C}$} & 
\multicolumn{1}{l}{$\prob{{\embed{A \rightarrow C}}}$} \\ \hline\hline
$\state{T}$ & 
$y_{1} + y_{3} + y_{5} + y_{6} + y_{7} + y_{8}$ \\ \hline
$\state{F}$ & 
$y_{2} + y_{4}$ \\ \hline
\end{tabular}
\end{eqnarray}

Following Equations \ref{eq:aandb} and \ref{eq:animplb} Let us use the
following set $\Gamma$ of probability constraints to express the
premises that $\bcond{A}{B}$ and $\bcond{B}{C}$ are both true:
\begin{eqnarray}
  \Gamma & \Leftarrow &
  \left\{ \:
  \prob{\embed{A \wedge B}=\tval} > 0, \:
  \prob{\embed{A \rightarrow B}=\fval} = 0, \:
  \prob{\embed{B \wedge C}=\tval} > 0, \:
  \prob{\embed{B \rightarrow C}=\fval} = 0
  \: \right\}
\end{eqnarray}
Parametric probability analysis yields the corresponding algebraic
expressions:
\begin{eqnarray}
  \Gamma & \Rightarrow &
  \left\{ \:
y_{1} + y_{2}
  > 0, \:
y_{3} + y_{4}
  = 0, \:
y_{1} + y_{5}
  > 0, \:
y_{2} + y_{6}
  = 0
  \: \right\}
\end{eqnarray}
First we query the feasible values of 
\begin{math}
  \prob{\embed{A \rightarrow C}=\fval}
\end{math}
subject to the constraints $\Gamma$ and $\Gamma_0$.  The paired linear
optimization problems:
\begin{equation}
\begin{array}{r@{\quad}l}
\mbox{Minimize}: & y_{2} + y_{4} \\
\mbox{subject to}:
& y_{1} + y_{2} + y_{3} + y_{4} + y_{5} + y_{6} + y_{7} + y_{8} = 1 \\
& y_{1} + y_{2} \ge \epsilon \\
& y_{3} + y_{4} = 0 \\
& y_{1} + y_{5} \ge \epsilon \\
& y_{2} + y_{6} = 0 \\
\mbox{and}:

& 0 \le y_{1} \le 1 \\
& 0 \le y_{2} \le 1 \\
& 0 \le y_{3} \le 1 \\
& 0 \le y_{4} \le 1 \\
& 0 \le y_{5} \le 1 \\
& 0 \le y_{6} \le 1 \\
& 0 \le y_{7} \le 1 \\
& 0 \le y_{8} \le 1 \\
& \epsilon = 0.001
\end{array}
\qquad
\begin{array}{r@{\quad}l}
\mbox{Maximize}: & y_{2} + y_{4} \\
\mbox{subject to}:
& y_{1} + y_{2} + y_{3} + y_{4} + y_{5} + y_{6} + y_{7} + y_{8} = 1 \\
& y_{1} + y_{2} \ge \epsilon \\
& y_{3} + y_{4} = 0 \\
& y_{1} + y_{5} \ge \epsilon \\
& y_{2} + y_{6} = 0 \\
\mbox{and}:

& 0 \le y_{1} \le 1 \\
& 0 \le y_{2} \le 1 \\
& 0 \le y_{3} \le 1 \\
& 0 \le y_{4} \le 1 \\
& 0 \le y_{5} \le 1 \\
& 0 \le y_{6} \le 1 \\
& 0 \le y_{7} \le 1 \\
& 0 \le y_{8} \le 1 \\
& \epsilon = 0.001
\end{array}
\end{equation}
have computed solutions
\begin{math}
  \alpha^* =
0.000
\end{math}
and
\begin{math}
  \beta^* =
0.000
\end{math}.
Thus it was computed that the solution set:
\begin{eqnarray}
  \{ \:
  \prob{\embed{A \rightarrow C}=\fval}
  \: : \: \Gamma
  \: \}
  & \Rightarrow & \{0\}
\end{eqnarray}
Using similar technique, linear optimization reveals that the solution set:
\begin{eqnarray}
  \{ \:
  \prob{\embed{A \wedge C}=\tval}
  \: : \: 
  \prob{\embed{A \rightarrow C}=\fval}=0, \:
  \Gamma
  \: \}
  & \subseteq &
  [
0.001
    ,
1.000
  ]
\end{eqnarray}
(with both end points
\begin{math}
0.001
\end{math}
and
\begin{math}
1.000
\end{math}
feasible), indicating that the queried objective must be strictly
greater than zero given the provided constraints (the precise
numerical results reflect our choice of the small constant
$\epsilon$).  In other words, these optimization results reveal that,
subject to the constraints that both conditionals $\bcond{A}{B}$ and
$\bcond{B}{C}$ are true, it is necessary that both of the following
constraints are satisfied:
\begin{eqnarray}
  \prob{\embed{A \rightarrow C}=\fval} & = & 0 \\
  \prob{\embed{A \wedge C}=\tval} & > & 0
\end{eqnarray}
These two constraints are precisely the criteria under which the
\rname{Boolean-feasibility} conditional $\bcond{A}{C}$ is true.  Hence
truth is the only feasible value for the conditional $\bcond{A}{C}$,
subject to the premises $(\bcond{A}{B})=\tval$ and
$(\bcond{B}{C})=\tval$; in other words $\Phi_0 \Rightarrow \{\tval\}$.
Thus we have computed using parametric probability and linear
optimization that the following metalevel \rname{Boolean-feasibility}
conditional is correct:
\begin{equation}
  \bcond{\{\:(\bcond{A}{B}),\:(\bcond{B}{C})\:\}}{(\bcond{A}{C})}
\end{equation}
This confirms the transitivity of \rname{Boolean-feasibility}
conditionals interpreted as statements of logical deduction.

Marginal probabilities:
\begin{equation}
\begin{tabular}[c]{l|l}\hline
\multicolumn{1}{l|}{$A$} & 
\multicolumn{1}{l}{$\prob{{A}}$} \\ \hline\hline
$\state{T}$ & 
$y_{1} + y_{2} + y_{3} + y_{4}$ \\ \hline
$\state{F}$ & 
$y_{5} + y_{6} + y_{7} + y_{8}$ \\ \hline
\end{tabular}
  \qquad
\begin{tabular}[c]{l|l}\hline
\multicolumn{1}{l|}{$B$} & 
\multicolumn{1}{l}{$\prob{{B}}$} \\ \hline\hline
$\state{T}$ & 
$y_{1} + y_{2} + y_{5} + y_{6}$ \\ \hline
$\state{F}$ & 
$y_{3} + y_{4} + y_{7} + y_{8}$ \\ \hline
\end{tabular}
  \qquad
\begin{tabular}[c]{l|l}\hline
\multicolumn{1}{l|}{$C$} & 
\multicolumn{1}{l}{$\prob{{C}}$} \\ \hline\hline
$\state{T}$ & 
$y_{1} + y_{3} + y_{5} + y_{7}$ \\ \hline
$\state{F}$ & 
$y_{2} + y_{4} + y_{6} + y_{8}$ \\ \hline
\end{tabular}
\end{equation}
Computed joint probabilities:
\begin{equation}
\begin{tabular}[c]{ll|l}\hline
\multicolumn{1}{l}{$A$} & 
\multicolumn{1}{l|}{$B$} & 
\multicolumn{1}{l}{$\prob{{A, B}}$} \\ \hline\hline
$\state{T}$ & 
$\state{T}$ & 
$y_{1} + y_{2}$ \\ \hline
$\state{T}$ & 
$\state{F}$ & 
$y_{3} + y_{4}$ \\ \hline
$\state{F}$ & 
$\state{T}$ & 
$y_{5} + y_{6}$ \\ \hline
$\state{F}$ & 
$\state{F}$ & 
$y_{7} + y_{8}$ \\ \hline
\end{tabular}
  \qquad
\begin{tabular}[c]{ll|l}\hline
\multicolumn{1}{l}{$B$} & 
\multicolumn{1}{l|}{$C$} & 
\multicolumn{1}{l}{$\prob{{B, C}}$} \\ \hline\hline
$\state{T}$ & 
$\state{T}$ & 
$y_{1} + y_{5}$ \\ \hline
$\state{T}$ & 
$\state{F}$ & 
$y_{2} + y_{6}$ \\ \hline
$\state{F}$ & 
$\state{T}$ & 
$y_{3} + y_{7}$ \\ \hline
$\state{F}$ & 
$\state{F}$ & 
$y_{4} + y_{8}$ \\ \hline
\end{tabular}
  \qquad
\begin{tabular}[c]{ll|l}\hline
\multicolumn{1}{l}{$A$} & 
\multicolumn{1}{l|}{$C$} & 
\multicolumn{1}{l}{$\prob{{A, C}}$} \\ \hline\hline
$\state{T}$ & 
$\state{T}$ & 
$y_{1} + y_{3}$ \\ \hline
$\state{T}$ & 
$\state{F}$ & 
$y_{2} + y_{4}$ \\ \hline
$\state{F}$ & 
$\state{T}$ & 
$y_{5} + y_{7}$ \\ \hline
$\state{F}$ & 
$\state{F}$ & 
$y_{6} + y_{8}$ \\ \hline
\end{tabular}
\end{equation}

\clearpage
\section{Examples}

Here are a few problems from the literature, analyzed by the proposed
methodology \cite{goodman,adams-conditionals,adams-subjunctive}.
These clever and well-crafted problems raise myriad issues in modeling
and analysis, all of which we can address using the methods introduced
above.


\subsection{Hot Buttered Conditionals}

We begin with Goodman's first two problems from \emph{Fact,
  Fiction, and Forecast} \cite{goodman}.  From his page~4:
\begin{quote}
  What, then, is the \emph{problem} about counterfactual conditionals?
  Let us confine ourselves to those in which antecedent and consequent
  are inalterably false---as, for example, when I say of a piece of
  butter that was eaten yesterday, and that had never been heated,
  \begin{quote}
    If that piece of butter had been heated to 150$^\circ$~F., it would
    have melted.
  \end{quote}
  Considered as truth-functional compounds, all counterfactuals are of
  course true, since their antecedents are false.  Hence
  \begin{quote}
    If that piece of butter had been heated to 150$^\circ$~F., it would
    not have melted
  \end{quote}
  would also hold.  Obviously something different is intended, and the
  problem is to define the circumstances under which a given
  counterfactual holds while the opposing counterfactual with the
  contradictory consequent fails to hold.
\end{quote}
Probability makes it easy to express `something different'.  Let us
use the probability model from Equation~\ref{eq:ab-input}, with $A$
denoting the event that the considered piece of butter was heated, and
$B$ denoting the event that the piece of butter melted.  We can
distill four statements from Goodman's problem description:
\begin{equation}
  \begin{array}{l|c|c|l} \hline
    \textsc{id} & \textsc{probability} & \textsc{polynomial} &
    \textsc{description} \\ \hline\hline
    \mathcal{S}_1 & \prob{A=\tval}=0 &
x
    = 0 &
    \mbox{It was not heated}
    \\ \hline
    \mathcal{S}_2 & \prob{B=\tval}=0 &
z + x y - x z
    = 0 &
    \mbox{It did not melt}
    \\ \hline
    \mathcal{S}_3 & \condprob[0]{B=\tval}{A=\tval}=1 &
y
    = 1 &
    \mbox{It would melt if heated}
    \\ \hline
    \mathcal{S}_4 & \condprob[0]{B=\tval}{A=\tval}=0 &
y
    = 0 &
    \mbox{It would not melt if heated}
    \\ \hline
  \end{array}
  \label{eq:butter}
\end{equation}
Here the premises $\mathcal{S}_1$ and $\mathcal{S}_2$ assert the facts
that the butter was not heated and that it did not melt.  The
conditionals $\mathcal{S}_3$ and $\mathcal{S}_4$ are counterfactual
(more specifically, antifactual) relative to the fact $\mathcal{S}_1$
that $A$ is known to be false.  More importantly, $\mathcal{S}_3$ and
$\mathcal{S}_4$ are also subjunctive: they concern the input
probability $\condprob[0]{B}{A}$ which does not involve the
probability $\prob{A}$.  The constraints in Equation~\ref{eq:butter}
behave exactly as Goodman intended.  It is consistent to assert either
subjunctive/counterfactual statement $\mathcal{S}_3$ or
$\mathcal{S}_4$ along with the facts $\mathcal{S}_1$ and
$\mathcal{S}_2$.  However it is inconsistent to assert both of the
opposite subjunctives together.

We can also consider Goodman's statement:
\begin{quote}
  Since that butter did not melt, it wasn't heated to 150$^\circ$~F
\end{quote}
Let us assume Goodman intends that the affirmative subjunctive stating
that the butter would have melted had it been heated is in effect
here.  In order to confirm the quoted statement of deduction, we can
query the feasible values of $\prob{A=\fval}$ subject to the
constraints $\prob{B=\tval}=0$ and $\condprob[0]{B=\tval}{A=\tval}=1$.
Using symbolic probability inference to find polynomial formulas
for these probability expressions reveals the solution set:
\begin{eqnarray}
  \left\{ \:
1 - x
  \: : \:
z + x y - x z
  = 0; \:
y
  = 1; \:
  x,y,z \in [0,1]
  \: \right\}
  & \Rightarrow & \{1\}
\end{eqnarray}
Here are the optimization problems generated to evaluate this solution
set:
\begin{equation}
\begin{array}{r@{\quad}l}
\mbox{Minimize}: & 1 - x \\
\mbox{subject to}:
& z + x y = x z \\
& y = 1 \\
\mbox{and}:
& 0 \le x \le 1 \\
& 0 \le y \le 1 \\
& 0 \le z \le 1
\end{array}
  \qquad
\begin{array}{r@{\quad}l}
\mbox{Maximize}: & 1 - x \\
\mbox{subject to}:
& z + x y = x z \\
& y = 1 \\
\mbox{and}:
& 0 \le x \le 1 \\
& 0 \le y \le 1 \\
& 0 \le z \le 1
\end{array}
\end{equation}
Solving the optimization problems gives computed minimum
\begin{math}
  a^* =
1.000
\end{math}
and maximum
\begin{math}
  b^* =
1.000
\end{math},
hence the solution set $\{1\}$ containing just the real number one.
In other words, it is a probabilistic consequence of the premises that
the butter did not melt ($\prob{B=\tval}=0$) and that the butter would
have melted had it been heated ($\condprob[0]{B=\tval}{A=\tval}=1$)
that the butter certainly was not heated ($\prob{A=\fval} \in \{1\}$).
It might be clearer to rephrase Goodman's statement of this deductive
result as:
\begin{quote}
  Since that butter did not melt, and it would have melted had it been
  heated---therefore it must not have been heated to 150$^\circ$~F.
\end{quote}
Separately, the facts that $A$ and $B$ are both false allow
us to deduce the subjunctive conditional
$\condprob[0]{B=\tval}{A=\fval}=0$ (that is,
\begin{math}
z
= 0
\end{math}).
In other words, the stated facts require that the piece of butter
certainly would not have melted, had it not been heated.  


\subsection{Conspiracy Theories}

We continue with Adams's problems from \cite{adams-subjunctive}
involving indicative and subjunctive conditionals.  For these
examples, probabilistic \rname{subjunctive} conditionals provide the
desired meanings.  There are several alternative ways to render the
indicative statements: as \rname{material}, \rname{existential}, or
\rname{feasibility} conditionals.  This author's intuition is that the
\rname{feasibility} interpretation is best; from this perspective the
indicative conditionals are regarded as statements of deduction (which
happen to be incorrect, for the examples that Adams provided).

\subsubsection{Murder Most Subjunctive}

For this example we are asked to consider a murder victim V and two
suspects A and B, with the evidence strongly favoring A as the
culprit.  We are given the following subjunctive conditional Su1 and
its indicative counterpart I1:
\begin{enumerate}
  \item[(Su1)] If A hadn't murdered V, B would not have either.
  \item[(I1)] If A didn't murder V, then B didn't either.
\end{enumerate}
Adams suggests that Su1 is `justified' whereas I1 is `unjustified'.
Let us calculate.  We identify the following propositions as the main
variables in our probability model:
\begin{quote}
  $A:$ A murdered V, \qquad
  $B:$ B murdered V, \qquad
  $V:$ V was murdered
\end{quote}
Adams's description suggests that $V$ requires either $A$ or $B$ (or
both)---a murder needs a murderer.  Hence we set the input
$\prob[0]{V|A,B}$ as though for the embedded propositional-calculus
formula $\embed{A \vee B}$.  The probability network is as follows:
\begin{equation}
  \vcenter{\xymatrix{
    *++[o][F]{A} \ar[r] \ar[dr] & *++[o][F]{B} \ar[d] \\
    & *++[o][F]{V}
  }}
  \qquad
\begin{tabular}[c]{l|l}\hline
\multicolumn{1}{l|}{$A$} & 
\multicolumn{1}{l}{$\prob[0]{{A}}$} \\ \hline\hline
$\state{T}$ & 
$x$ \\ \hline
$\state{F}$ & 
$1 - x$ \\ \hline
\end{tabular}
  \qquad
\begin{tabular}[c]{l|ll}\hline
\multicolumn{3}{l}{$\condprob[0]{{B}}{{A}}$} \\ \hline\hline
\multicolumn{1}{l|}{$A$} & 
\multicolumn{1}{l}{${B=\state{T}}$} & 
\multicolumn{1}{l}{${B=\state{F}}$} \\ \hline\hline
$\state{T}$ & 
$y$ & 
$1 - y$ \\ \hline
$\state{F}$ & 
$z$ & 
$1 - z$ \\ \hline
\end{tabular}
  \qquad
\begin{tabular}[c]{ll|ll}\hline
\multicolumn{4}{l}{$\condprob[0]{{V}}{{A, B}}$} \\ \hline\hline
\multicolumn{1}{l}{$A$} & 
\multicolumn{1}{l|}{$B$} & 
\multicolumn{1}{l}{${V=\state{T}}$} & 
\multicolumn{1}{l}{${V=\state{F}}$} \\ \hline\hline
$\state{T}$ & 
$\state{T}$ & 
$1$ & 
$0$ \\ \hline
$\state{T}$ & 
$\state{F}$ & 
$1$ & 
$0$ \\ \hline
$\state{F}$ & 
$\state{T}$ & 
$1$ & 
$0$ \\ \hline
$\state{F}$ & 
$\state{F}$ & 
$0$ & 
$1$ \\ \hline
\end{tabular}
  \label{eq:abv-input}
\end{equation}
with each parameter $x$, $y$, and $z$ constrained to the real interval
$[0,1]$.

Consider the following, which includes two ways to render the
indicative I1 (\rname{material} or \rname{existential}):
\begin{equation}
  \begin{array}{l|c|c|l} \hline
    \textsc{id} & \textsc{probability} & \textsc{polynomial} &
    \textsc{description} \\ \hline\hline
    \textrm{Su1} & \condprob[0]{B=\tval}{A=\fval}=0 &
z
    = 0 &
    \mbox{If A hadn't, B wouldn't have}
    \\ \hline
    \textrm{I1-M} & 
    \prob{V=\tval,A=\fval,B=\tval}=0 &
z - x z
    = 0 &
    \mbox{Either A did or B didn't}
    \\ \hline
    \textrm{I1-E} & 
    \begin{array}{@{}r@{}}
      \prob{V=\tval,A=\fval,B=\tval}=0 \\ 
      \prob{V=\tval,A=\fval}>0
    \end{array}
    &
    \begin{array}{@{}r@{}}
z - x z
    = 0 \\
z - x z
    > 0
    \end{array}
    &
    \begin{array}{@{}l@{}}
      \mbox{If A didn't then B didn't} \\
      \mbox{(nonzero chance that A didn't)}
    \end{array}
    \\ \hline
  \end{array}
  \label{eq:abv}
\end{equation}
Both indicative constraints assert
$\condprob{B=\tval}{V=\tval,A=\fval}=0$, or equivalently
$\condprob{B=\fval}{V=\tval,A=\fval}=1$ (it is the same to say
`certainly-not true' as `certainly false'); the \rname{material}
interpretation I1-M also allows the indefinite value $0/0$ for these
conditional probabilities.  Anyway under the \rname{material}
probabilistic interpretation I1-M, we would accept proposition I1 as a
true statement and conclude that A must have murdered V, if in fact V
was murdered.  Under the \rname{existential} probabilistic
interpretation I1-E, we would reject proposition I1 as inconsistent: it
asserts that there is zero probability that B murdered V and A didn't,
and simultaneously that this same probability is strictly greater than
zero.  Note what the indicative I1 says when it is interpreted as a
\rname{material} conditional:
\begin{enumerate}
  \item[(I1)] Either A murdered V, or B didn't murder V
\end{enumerate}

We can also interpret I1 as a statement of deduction, let us say I1-F,
claiming that `Assuming that V was murdered and A didn't do
it---therefore B didn't either.'  This statement of deduction asserts
that the following solution set is $\{1\}$ (indicating that $B$ is
certainly false, given the premise):
\begin{eqnarray}
  \left\{ \:
  \prob{B=\fval} \: : \:
  \prob{V=\tval,A=\fval} = 1
  \: \right\}
  & = & \{1\}
  \label{eq:i1-feasibility}
\end{eqnarray}
We shall see presently that this assertion is unsatisfiable.  After
substituting the results of symbolic probability inference the
requested solution set is given by:
\begin{eqnarray}
  \left\{ \:
1 - z - x y + x z
  \: : \:
z - x z
  = 1 ; \:
  x,y,z \in [0,1]
  \: \right\}
  & \Rightarrow & \{0\}
\end{eqnarray}
In this case the only feasible value of the objective is $0$.  It
follows that Equation~\ref{eq:i1-feasibility} could be satisfied only
if $0=1$.  Therefore the correct deductive statement is the opposite
of I1-F, namely `Assuming V was murdered and A didn't do it---therefore
B must have' or more succinctly `If A didn't murder V, then B did'.
We would get the same result having started from a slightly different
set comprehension expression that splits the constraints on $A$ and
$V$ into two equations:
\begin{equation}
  \left\{ \:
  \prob{B=\fval} \: : \:
  \prob{V=\tval}=1; \:
  \prob{A=\fval}=1
  \: \right\}
\end{equation}

All of the indicative interpretations I1-M, I1-E, and I1-F are factual
with respect to the fact that V was murdered (`factual' meaning that
the conditions explicitly include the known event $V=\tval$).  On the
other hand the subjunctive interpretation of Su1 is afactual (that is,
weakly counterfactual) about the fact that V was murdered (`afactual'
meaning that the known event $V=\tval$ does not appear in the
condition part of the input table $\condprob[0]{B}{A}$).

The following output probability tables were used to generate the
results just presented:
\begin{equation}
\begin{tabular}[c]{r|ll|l}\hline
\multicolumn{1}{l|}{\scshape {\#}} & 
\multicolumn{1}{l}{$V$} & 
\multicolumn{1}{l|}{$A$} & 
\multicolumn{1}{l}{$\prob{{V, A}}$} \\ \hline\hline
1 & 
$\state{T}$ & 
$\state{T}$ & 
$x$ \\ \hline
2 & 
$\state{T}$ & 
$\state{F}$ & 
$z - x z$ \\ \hline
3 & 
$\state{F}$ & 
$\state{T}$ & 
$0$ \\ \hline
4 & 
$\state{F}$ & 
$\state{F}$ & 
$1 - x - z + x z$ \\ \hline
\end{tabular}
  \qquad
\begin{tabular}[c]{r|lll|l}\hline
\multicolumn{1}{l|}{\scshape {\#}} & 
\multicolumn{1}{l}{$V$} & 
\multicolumn{1}{l}{$A$} & 
\multicolumn{1}{l|}{$B$} & 
\multicolumn{1}{l}{$\prob{{V, A, B}}$} \\ \hline\hline
1 & 
$\state{T}$ & 
$\state{T}$ & 
$\state{T}$ & 
$x y$ \\ \hline
2 & 
$\state{T}$ & 
$\state{T}$ & 
$\state{F}$ & 
$x - x y$ \\ \hline
3 & 
$\state{T}$ & 
$\state{F}$ & 
$\state{T}$ & 
$z - x z$ \\ \hline
4 & 
$\state{T}$ & 
$\state{F}$ & 
$\state{F}$ & 
$0$ \\ \hline
5 & 
$\state{F}$ & 
$\state{T}$ & 
$\state{T}$ & 
$0$ \\ \hline
6 & 
$\state{F}$ & 
$\state{T}$ & 
$\state{F}$ & 
$0$ \\ \hline
7 & 
$\state{F}$ & 
$\state{F}$ & 
$\state{T}$ & 
$0$ \\ \hline
8 & 
$\state{F}$ & 
$\state{F}$ & 
$\state{F}$ & 
$1 - x - z + x z$ \\ \hline
\end{tabular}
\end{equation}
\begin{equation}
\begin{tabular}[c]{l|l}\hline
\multicolumn{1}{l|}{$V$} & 
\multicolumn{1}{l}{$\prob{{V}}$} \\ \hline\hline
$\state{T}$ & 
$x + z - x z$ \\ \hline
$\state{F}$ & 
$1 - x - z + x z$ \\ \hline
\end{tabular}
  \qquad
\begin{tabular}[c]{l|l}\hline
\multicolumn{1}{l|}{$A$} & 
\multicolumn{1}{l}{$\prob{{A}}$} \\ \hline\hline
$\state{T}$ & 
$x$ \\ \hline
$\state{F}$ & 
$1 - x$ \\ \hline
\end{tabular}
  \qquad
\begin{tabular}[c]{l|l}\hline
\multicolumn{1}{l|}{$B$} & 
\multicolumn{1}{l}{$\prob{{B}}$} \\ \hline\hline
$\state{T}$ & 
$z + x y - x z$ \\ \hline
$\state{F}$ & 
$1 - z - x y + x z$ \\ \hline
\end{tabular}
\end{equation}
Also note:
\begin{equation}
\begin{tabular}[c]{r|lll|l}\hline
\multicolumn{1}{l|}{\scshape {\#}} & 
\multicolumn{1}{l}{$A$} & 
\multicolumn{1}{l}{$V$} & 
\multicolumn{1}{l|}{$B$} & 
\multicolumn{1}{l}{$\condprob{{V, B}}{{A}}$} \\ \hline\hline
1 & 
$\state{T}$ & 
$\state{T}$ & 
$\state{T}$ & 
$x y / x$ \\ \hline
2 & 
$\state{T}$ & 
$\state{T}$ & 
$\state{F}$ & 
$\left(x - x y\right) / \left(x\right)$ \\ \hline
3 & 
$\state{T}$ & 
$\state{F}$ & 
$\state{T}$ & 
$\left(0\right) / \left(x\right)$ \\ \hline
4 & 
$\state{T}$ & 
$\state{F}$ & 
$\state{F}$ & 
$\left(0\right) / \left(x\right)$ \\ \hline
5 & 
$\state{F}$ & 
$\state{T}$ & 
$\state{T}$ & 
$\left(z - x z\right) / \left(1 - x\right)$ \\ \hline
6 & 
$\state{F}$ & 
$\state{T}$ & 
$\state{F}$ & 
$\left(0\right) / \left(1 - x\right)$ \\ \hline
7 & 
$\state{F}$ & 
$\state{F}$ & 
$\state{T}$ & 
$\left(0\right) / \left(1 - x\right)$ \\ \hline
8 & 
$\state{F}$ & 
$\state{F}$ & 
$\state{F}$ & 
$\left(1 - x - z + x z\right) / \left(1 - x\right)$ \\ \hline
\end{tabular}
\end{equation}

Consider also the matter of explaining inconsistency and revising the
probability model to address it.  For example, if it is asserted that
Su1 holds but that V certainly has been murdered and that A certainly
is innocent, what can we do to avoid the implied contradiction?  One
course might be to adjust probabilities, for example to interpret Su1
as the constraint
\begin{math}
z
  \approx 0
\end{math} instead of as a strict equality
(allowing for example
\begin{math}
z
  = \delta_1
\end{math}
for some small constant $\delta_1$).  Another revision might be to change
$\condprob[0]{V=\tval}{A=\fval,B=\fval}$ from $0$ to some small
constant $\delta_2$ to represent the idea that someone other than A or
B might have murdered V.

\subsubsection{Kennedy and Oswald}

We next consider the problems in \cite{adams-subjunctive} about the
late President John F.\ Kennedy and his assassin Lee Harvey Oswald.
It is generally accepted as historical fact that Kennedy was shot and
killed by Oswald in Dallas in 1963.  Adams provided two pairs of
conditionals.  The first pair:
\begin{enumerate}
  \item[(Su2)] If Oswald hadn't shot Kennedy in Dallas, then no one
    else would have.
  \item[(I2)] If Oswald didn't shoot Kennedy in Dallas, then no one
    else did.
\end{enumerate}
And the second pair:
\begin{enumerate}
  \item[(Su3)] If Oswald hadn't shot Kennedy in Dallas,
    Kennedy would be alive today.
  \item[(I3)] If Oswald didn't shoot Kennedy in Dallas,
    then Kennedy is alive today.
\end{enumerate}
In order to analyze these conditionals, we revise the probability
model from Equation~\ref{eq:abv-input} to change the referents of the
variables $A$, $B$, and $V$:
\begin{quote}
  $A:$ Oswald shot Kennedy, \qquad
  $B:$ Someone else shot Kennedy, \qquad
  $V:$ Kennedy was shot
\end{quote}
We add a new variable $L$:
\begin{quote}
  $L:$ Kennedy is alive today
\end{quote}
The updated probability-network graph and the additional input table
for $L$ are as follows:
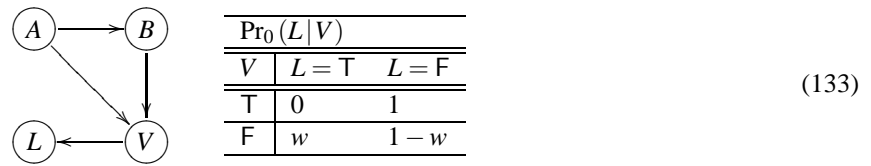
\begin{equation}
  \vcenter{\xymatrix{
    *++[o][F]{A} \ar[r] \ar[dr] & *++[o][F]{B} \ar[d] \\
    *++[o][F]{L }& *++[o][F]{V} \ar[l]
  }}
  \qquad
\begin{tabular}[c]{l|ll}\hline
\multicolumn{3}{l}{$\condprob[0]{{L}}{{V}}$} \\ \hline\hline
\multicolumn{1}{l|}{$V$} & 
\multicolumn{1}{l}{${L=\state{T}}$} & 
\multicolumn{1}{l}{${L=\state{F}}$} \\ \hline\hline
$\state{T}$ & 
$0$ & 
$1$ \\ \hline
$\state{F}$ & 
$w$ & 
$1 - w$ \\ \hline
\end{tabular}
  \label{eq:abvl-input}
\end{equation}
The input table $\condprob[0]{L}{V}$ says that Kennedy certainly would
not be alive today had he been shot; however had he not been shot,
there is some probability $w$ that he would now be alive.  As with the
other parameters $x$, $y$, and $z$, we have the constraint $0 \le w
\le 1$.

In terms of probability, Su2 and I2 are just like Su1 and I1 above.
The subjunctive conditional Su2 is the probability constraint
$\condprob[0]{B=\tval}{A=\fval}=0$, equivalently
$\condprob[0]{B=\fval}{A=\fval}=1$, which yields the polynomial
constraint
\begin{math}
z
  = 0
\end{math}.  As with I1, the indicative statement I2 may be
interpreted in a variety of ways: as a \rname{material},
\rname{existential}, or \rname{feasibility} conditional.  The
\rname{existential} interpretation I2-E offers the infeasible
constraints $\condprob{B=\tval}{V=\tval,A=\fval}=0$ 
and $\prob{V=\tval,A=\fval}>0$, which together assert the incompatible
premises that there is nonzero probability that Kennedy was shot by
someone other than Oswald, yet that in that case no one else could
have done it.  The \rname{material} interpretation I2-M also forbids
the event that someone else shot Kennedy if Kennedy was not shot by
Oswald, but allows that there could be zero probability that Kennedy
was shot by someone other than Oswald; in this zero-probability case
the conditional probabilities $\condprob{B=\tval}{V=\tval,A=\fval}$
and $\condprob{B=\fval}{V=\tval,A=\fval}$ have the indeterminate value
$0/0$.  The \rname{feasibility} interpretation I2-F provides the
incorrect deductive statement that it is a probabilistic consequence
of the premise that Kennedy was shot, and not by Oswald, that someone
else must not have shot Kennedy either.  The correct deduction is that
indeed someone else must have shot Kennedy in that case.

Moving along we now consider statements Su3 and I3.  Following
Table~\ref{tbl:principals} the subjunctive conditional Su3 should be
rendered as the constraint $\condprob[0]{L=\tval}{A=\fval}=1$ which
says that the input probability that Kennedy is alive given that
Oswald had not shot him must be one.  However in our probability model
there is no such input table $\condprob[0]{L}{A}$.  We can find an
alternative means to express this particular subjunctive conditional
Su3 using members of the input tables $\condprob[0]{L}{V}$ and
$\condprob[0]{B}{A}$ from Equations \ref{eq:abv-input} and
\ref{eq:abvl-input}.  Following the logic of the example, Su3 requires
both $\condprob[0]{B=\fval}{A=\fval}=1$ and
$\condprob[0]{L=\tval}{V=\fval}=1$: the first saying that no one else
would have shot Kennedy had Oswald not done it, and the second saying
that Kennedy would still be alive had he not been shot.  These
probability constraints yield the polynomial constraints:
\begin{math}
z
  = 0
\end{math}
and
\begin{math}
w
  = 1
\end{math}.
Because the variables $w$ and $z$ are constrained to lie between zero
and one, it is equivalent to provide the following single constraint
to represent the subjunctive conditional Su3:
\begin{eqnarray}
w
(
1 - z
) & = & 1
\end{eqnarray}
Note the polynomial quotient computed by symbolic probability
inference for the conditional probability that Kennedy
would be alive today, given that Oswald had not shot him:
\begin{eqnarray}
  \condprob{L=\tval}{A=\fval} & \Rightarrow &
\frac{w - x w - z w + x z w}{1 - x}
\end{eqnarray}
This quotient factors into the following expression, whose numerator
and denominator are products of input probabilities:
\begin{equation}
\frac{
(
1 - x
)
(
1 - z
)
(
1
)
(
w
)
}{
1 - x
}
\end{equation}
Here the denominator is $\prob[0]{A=\fval}$ and the numerator is the
product:
\begin{equation}
  \prob[0]{A=\fval} \cdot
  \condprob[0]{B=\fval}{A=\fval} \cdot
  \condprob[0]{V=\fval}{A=\fval,B=\fval} \cdot
  \condprob[0]{L=\tval}{V=\fval}
\end{equation}
Eliminating the input probability $\prob[0]{A=\fval}$ from the
numerator and denominator (without assuming it must be nonzero!)
yields the desired input-probability expression to constrain for the
subjunctive conditional Su3.

There are several ways to interpret the indicative conditional I3.
\begin{equation}
  \begin{array}{l|c|c|l} \hline
    \textsc{id} & \textsc{probability} & \textsc{polynomial} &
    \textsc{description} \\ \hline\hline
    \textrm{Su3} & 
    \condprob[0]{L=\tval}{V=\fval} \cdot
    \condprob[0]{B=\fval}{A=\fval} = 1 &
w
    (
1 - z
    ) = 1 &
    \mbox{If O. hadn't shot, K. would be alive}
    \\ \hline
    \textrm{I3-M} & 
    \prob{V=\tval,A=\fval,L=\tval}=\prob{V=\tval,A=\fval} &
0
    =
z - x z
    &
    \mbox{Either O. shot or K. is alive}
    \\ \hline
    \textrm{I3-E} & 
    \begin{array}{@{}c@{}}
      \prob{V=\tval,A=\fval,L=\tval}=\prob{V=\tval,A=\fval} \\
      \prob{V=\tval,A=\fval}>0
    \end{array}
    &
    \begin{array}{@{}c@{}}
0
    =
z - x z
    \\
z - x z
    > 0
    \end{array}
    &
    \begin{array}{@{}l@{}}
      \mbox{If O. didn't shoot then K. is alive} \\
      \mbox{(nonzero chance that O. didn't)}
    \end{array}
    \\ \hline
  \end{array}
\end{equation}
Consider the joint probability:
\begin{equation}
\begin{tabular}[c]{r|lll|l}\hline
\multicolumn{1}{l|}{\scshape {\#}} & 
\multicolumn{1}{l}{$V$} & 
\multicolumn{1}{l}{$A$} & 
\multicolumn{1}{l|}{$L$} & 
\multicolumn{1}{l}{$\prob{{V, A, L}}$} \\ \hline\hline
1 & 
$\state{T}$ & 
$\state{T}$ & 
$\state{T}$ & 
$0$ \\ \hline
2 & 
$\state{T}$ & 
$\state{T}$ & 
$\state{F}$ & 
$x$ \\ \hline
3 & 
$\state{T}$ & 
$\state{F}$ & 
$\state{T}$ & 
$0$ \\ \hline
4 & 
$\state{T}$ & 
$\state{F}$ & 
$\state{F}$ & 
$z - x z$ \\ \hline
5 & 
$\state{F}$ & 
$\state{T}$ & 
$\state{T}$ & 
$0$ \\ \hline
6 & 
$\state{F}$ & 
$\state{T}$ & 
$\state{F}$ & 
$0$ \\ \hline
7 & 
$\state{F}$ & 
$\state{F}$ & 
$\state{T}$ & 
$w - x w - z w + x z w$ \\ \hline
8 & 
$\state{F}$ & 
$\state{F}$ & 
$\state{F}$ & 
$1 - x - z - w + x z + x w + z w - x z w$ \\ \hline
\end{tabular}
\end{equation}

Interpreted as a statement of deduction, the indicative I3-F asserts
that one is the only feasible value of the probability that Kennedy is
still alive, given that he was shot and not by Oswald:
\begin{eqnarray}
  \left\{ \:
  \prob{L=\tval} \: : \:
  \prob{V=\tval,A=\fval} = 1
  \: \right\}
  & = & \{1\}
\end{eqnarray}
But this equation is incorrect; the actual solution set is $\{0\}$.
Here is the set-comprehension expression after symbolic probability
inference:
\begin{eqnarray}
  \left\{ \:
w - x w - z w + x z w
  \: : \:
z - x z
  = 1; \:
  w,x,y,z \in [0,1]
  \: \right\}
  & \Rightarrow & \{0\}
\end{eqnarray}
You can see by inspection that the constraints require $x=0$ and
$z=1$; with these substitutions the objective simplifies to the
constant $0$.  The same result would come from specifying the premise
as two probability constraints instead of one:
\begin{eqnarray}
  \left\{ \:
  \prob{L=\tval} \: : \:
  \prob{V=\tval} = 1, \:
  \prob{A=\fval} = 1
  \: \right\}
  & \Rightarrow & \{0\}
\end{eqnarray}

Note that omitting the fact $V=\tval$ that Kennedy was shot yields a
different result for the indicative conditional I3-F (which thereby
becomes afactual---weakly counterfactual---with respect to Kennedy's
shooting).  In this case we evaluate the solution set:
\begin{equation}
  \left\{ \:
  \prob{L=\tval} \: : \:
  \prob{A=\fval} = 1
  \: \right\}
\end{equation}
The algebraic expression is the following:
\begin{eqnarray}
  \left\{ \:
w - x w - z w + x z w
  \: : \:
1 - x
  = 1; \:
  w,x,y,z \in [0,1]
  \: \right\}
  & \Rightarrow & [0,1]
\end{eqnarray}
In this case the constraints require $x=0$, with which substitution
the objective $\prob{L=\tval}$ simplifies to the product $w(1-z)$.
Subject to the constraint that Oswald did not shoot him, but allowing
that he may not have been shot at all, the probability that Kennedy is
alive today depends on the probability $z$ that someone else would
have shot him (had Oswald not done so) and the probability $w$ that he
would still be alive had he not been shot.  Absent other constraints
on the values of these parameters, the value of the objective $w(1-z)$
could have any value between zero and one.  In this case I3 still
would not be a correct statement of deduction, because the computed
solution set $[0,1]$ is not exactly $\{1\}$.  However it would be a
correct statement of probabilistic deduction that:
\begin{quote}
  Ignoring the fact that Kennedy was shot, if Oswald didn't shoot
  Kennedy in Dallas, then Kennedy \emph{may or may not be} alive today.
\end{quote}
Here `may or may not be' reflects the idea that it is feasible for the
queried probability to be anywhere between zero and one (including
these end points), subject to the given constraints.

For clarity we might adopt the convention to preface factual
conditionals with a phrase like, `Accepting the fact that \ldots' when
stating them in natural language.  Thus we might clarify which
indicative I3 is intended.  For example, there is the factual
\rname{material} conditional interpretation of I3:
\begin{quote}
  Accepting the fact that Kennedy was shot, either Oswald shot Kennedy
  in Dallas or Kennedy is alive today (or both).
\end{quote}
And there is a corresponding afactual \rname{material} conditional:
\begin{quote}
  Ignoring the fact that Kennedy was shot, either Oswald shot Kennedy
  in Dallas or Kennedy is alive today (or both).
\end{quote}
Both of these happen to be correct given the information provided.

\subsubsection{Soft, What Conditional Breaks?}

We now consider the role of observed evidence.  We are given the
following conditionals:
\begin{enumerate}
  \item[(Su4)] $X$ is soft at time $t$ $=_{\mathrm{df}}$ if $X$ should
    be (were, had been, depending on the relation of $t$ to the
    present) subject to moderate deforming pressure at time $t$, then
    it would be (would have been) significantly deformed.
  \item[(I4)] If $X$ is (was) subject to moderate deforming pressure
    at time $t$, then it will be (is, was) significantly deformed.
\end{enumerate}
We are asked to ``\ldots suppose that we have observed that at time
$t$, $X$ was not deformed, but that we don't know whether it was
subject to deforming stress at that time.''  Given this observation,
we are asked to evaluate the truth of the conditionals Su4 and I4.  In
order to analyze this problem let us return to the basic probability
model from Equation~\ref{eq:ab-input}.  Let us assume that we are
concerned with just one moment in time $t$.  We update the referents
of the main variables $A$ and $B$ to the following propositions:
\begin{quote}
  $A:$ $X$ was subject to moderate deforming pressure, \qquad
  $B:$ $X$ was significantly deformed
\end{quote}
The observation that $X$ was not deformed constitutes the evidence
that $B$ is certainly false, in other words $\prob{B=\tval}=0$ or
equivalently $\prob{B=\fval}=1$.  By symbolic probability inference
this evidence becomes the constraint
\begin{math}
z + x y = x z
\end{math}.

What can we say about the conditionals Su4 and I4 given this evidence?
Following Equation~\ref{eq:subjunctive} the affirmative
\rname{subjunctive} conditional Su4 corresponds to the equation
$\condprob[0]{B=\tval}{A=\tval}=1$.  Therefore, in order to evaluate
the conditional Su4, let us analyze the set $\Phi_1$ of feasible
values for the input probability $\condprob[0]{B=\tval}{A=\tval}$
subject to the given evidence $\prob{B=\tval}=0$ along with the
general constraints $\Gamma_0$ of the probability model:
\begin{eqnarray}
  \Phi_1 & \Leftarrow &
  \left\{ \:
  \condprob[0]{B=\tval}{A=\tval} \: : \:
  \prob{B=\tval}=0, \: \Gamma_0
  \: \right\}
\end{eqnarray}
Using symbolic probability inference with the probability model from
Equation~\ref{eq:ab-input}, we generate the following pair of
polynomial optimization problems to characterize this solution set
$\Phi_1$:
\begin{equation}
\begin{array}{r@{\quad}l}
\mbox{Minimize}: & y \\
\mbox{subject to}:
& z + x y = x z \\
\mbox{and}:
& 0 \le x \le 1 \\
& 0 \le y \le 1 \\
& 0 \le z \le 1
\end{array}
\qquad
\begin{array}{r@{\quad}l}
\mbox{Maximize}: & y \\
\mbox{subject to}:
& z + x y = x z \\
\mbox{and}:
& 0 \le x \le 1 \\
& 0 \le y \le 1 \\
& 0 \le z \le 1
\end{array}
\end{equation}
The computed minimum and maximum solutions are
\begin{math}
  \alpha^* =
0.000
\end{math}
and
\begin{math}
  \beta^* =
1.000
\end{math},
from which it follows that the solution set $\Phi_1 \subseteq [0,1]$
with $\Phi_1$ including at least the points $0$ and $1$.  In other
words, the evidence that $B$ is certainly false does not constrain the
value of the objective $\condprob[0]{B=\tval}{A=\tval}$ to any
particular value.  Therefore the provided evidence does not tell us
anything about whether the subjunctive Su4 is true or false, leaving
it possible that Su4 ``might well be justified'' by the evidence (to
borrow Adams's phrasing).

In order to check whether the affirmative \rname{material} or
conditional interpretation of I4 might hold given the provided
evidence $\prob{B=\tval}=0$, we can evaluate the following solution
set $\Phi_2$:
\begin{eqnarray}
  \Phi_2 & \Leftarrow &
  \{ \:
  \prob{A=\tval}-\prob{A=\tval,B=\tval} \: : \:
  \prob{B=\tval}=0, \:
  \Gamma_0
  \: \}
\end{eqnarray}
Analysis according to Section~\ref{sec:optimization} gives
minimum and maximum solutions
\begin{math}
  \alpha^* =
0.000
\end{math}
and
\begin{math}
  \beta^* =
1.000
\end{math},
indicating $\Phi_2 \subseteq [0,1]$ with $0\in\Phi_2$ and
$1\in\Phi_2$.  Therefore, subject to the given evidence, it is
possible but not necessary that the constraint
\begin{math}
  \prob{A=\tval,B=\tval}=\prob{A=\tval}
\end{math}
from Equation~\ref{eq:material} defining the affirmative
\rname{material} conditional is satisfied.  Hence the \rname{material}
interpretation of the indicative conditional I4 is consistent with the
supplied evidence (though not required by the evidence).  Note what
the \rname{material} interpretation of the indicative I4 says:
\begin{quote}
  \raggedright
  Either $X$ was not subject to moderate deforming pressure at time
  $t$ or it was significantly deformed.
\end{quote}

In order to investigate the \rname{existential} interpretation of I4
subject to the provided evidence $\prob{B=\tval}=0$, we can evaluate
the following solution set $\Phi_3$:
\begin{eqnarray}
  \Phi_3 & \Leftarrow &
  \{ \:
  \prob{A=\tval} \: : \:
  \prob{A=\tval,B=\tval}=\prob{A=\tval}, \:
  \prob{B=\tval}=0, \:
  \Gamma_0
  \: \}
\end{eqnarray}
Analysis according to Section~\ref{sec:optimization} gives
minimum and maximum solutions
\begin{math}
  \alpha^* =
0.000
\end{math}
and
\begin{math}
  \beta^* =
0.000
\end{math},
indicating $\Phi_3 \Rightarrow \{0\}$.  In other words, it is not
feasible for $\prob{A=\tval}$ to be greater than zero.  Therefore,
subject to the provided evidence, it is not possible to satisfy
Equation~\ref{eq:existential} with $k=1$ indicating an affirmative
\rname{existential} conditional, because the constraint
$\prob{A=\tval}>0$ is violated.  Hence the \rname{existential}
interpretation of the indicative conditional I4 is not compatible with
the supplied evidence.

By the foregoing calculations, probabilistic analysis has demonstrated
that the \rname{subjunctive} conditional Su4 is compatible with the
provided evidence, as is the \rname{material} interpretation of the
indicative I4.  However the \rname{existential} interpretation of the
indicative conditional I4 is not compatible with the stated evidence.
Neither can the \rname{feasibility} interpretation of the indicative
conditional I4 hold, unless the stated evidence $\prob{B=\tval}=0$ is
ignored.

\subsubsection{A Brown Bird in the Bush}

Finally we consider the last example from \cite{adams-subjunctive}.
We are given ``a natural law that all ravens are black'' and the
observation of ``a brown bird in the distance,'' along with two
propositions:
\begin{enumerate}
  \item[(Su5)] If that were a raven, it would be black
  \item[(I5)] If that is a raven, it is black
\end{enumerate}
We recycle the probability model from Equation~\ref{eq:ab-input} with
the referents of the variables modified to the following:
\begin{quote}
  $A:$ It is a raven,
  \qquad
  $B:$ It is black
\end{quote}
This example is similar to the last one except that the truth of the
subjunctive conditional is now given as evidence.  Using
Equation~\ref{eq:subjunctive} with $k=1$ to indicate an affirmative
statement, we render the subjunctive conditional Su5 as the
probability equation
\begin{math}
  \condprob[0]{B=\tval}{A=\tval} = 1
\end{math},
which yields the polynomial equation
\begin{math}
y
  = 1
\end{math}.
We interpret the natural law about ravens being black as the assertion
that the subjunctive Su5 is true (hence as the constraint
\begin{math}
y
  = 1
\end{math}).  
We also have as evidence the assertion
\begin{math}
  \prob{B=\tval}=0
\end{math}
that the observed bird was certainly not black; this gives the
polynomial constraint
\begin{math}
z + x y = x z
\end{math}.

In order to check whether the affirmative \rname{material}
interpretation of I5 might hold given the provided evidence, we can
evaluate the following solution set $\Phi_5$:
\begin{eqnarray}
  \Phi_5 & \Leftarrow &
  \{ \:
  \prob{A=\tval}-\prob{A=\tval,B=\tval} \: : \:
  \condprob[0]{B=\tval}{A=\tval}=1, \:
  \prob{B=\tval}=0, \:
  \Gamma_0
  \: \}
\end{eqnarray}
Analysis according to Section~\ref{sec:optimization} gives
minimum and maximum solutions
\begin{math}
  \alpha^* =
0.000
\end{math}
and
\begin{math}
  \beta^* =
1.000
\end{math}.
Hence $\Phi_5 \subseteq [0,1]$ with $0\in\Phi_5$ and $1\in\Phi_5$ and
it is possible but not necessary that the \rname{material}
interpretation of the indicative conditional I5 holds, subject to the
provided evidence.

For the affirmative \rname{existential} interpretation of I5, we
evaluate the following solution set $\Phi_6$:
\begin{equation}
  \{ \:
  \prob{A=\tval} \: : \:
  \prob{A=\tval}=\prob{A=\tval,B=\tval}, \:
  \condprob[0]{B=\tval}{A=\tval}=1, \:
  \prob{B=\tval}=0, \:
  \Gamma_0
  \: \}
\end{equation}
Analysis according to Section~\ref{sec:optimization} gives
minimum and maximum solutions
\begin{math}
  \alpha^* =
0.000
\end{math}
and
\begin{math}
  \beta^* =
0.000
\end{math},
hence $\Phi_6 \Rightarrow \{0\}$ and the \rname{existential}
interpretation of I5 is inconsistent with the given evidence (because
the constraint $\prob{A=\tval}>0$ cannot be satisfied).

\subsection{Planes, Politicians, and Passing Exams}

Finally we address Adams's earlier problems from
\cite{adams-conditionals}.  Rather than trying to figure out when the
propositional calculus may be `safely' used, we can dispense with it
altogether and model these problems directly with conditional
probabilities.  These problems point out two kinds of fallacies:
first, confusion about division by zero; second, mishandling
specializations of general rules.  We address these by calling out
explicitly when a quotient has the indefinite value $0/0$; and by
understanding that general conditional statements may have exceptions
when they are specialized to account for new antecedent terms.  Note
that such revision of conditional statements can be done within the
propositional calculus too.  However, proper treatment of impossible
antecedents would be cumbersome in the propositional calculus.

All but one of the problems from \cite{adams-conditionals} can be
solved using the \rname{Boolean-feasibility} conditionals along with
the probability and linear optimization methods described in
Section~\ref{sec:prob-boole}.  The exception is \ref{ex:f4} which
illustrates a different point: sometimes an initial conditional
statement ought to be revised in order to account for new information
(or alternatively it should have been stated in a defeasible manner in
the first place, to allow for future exceptions).  Stubborn adherence
to premature generalizations is a special kind of fallacy.

\renewcommand{\theenumi}{F\arabic{enumi}}
\renewcommand{\labelenumi}{\theenumi.}
\begin{enumerate}
  \raggedright
  \item \label{ex:f1} \emph{John will arrive on the 10 o'clock plane.
    Therefore, if John does not arrive on the 10 o'clock plane, he
    will arrive on the 11 o'clock plane.}

    Variables:
    \begin{itemize}
      \item[$A:$] John arrives on the 10:00 plane
      \item[$B:$] John arrives on the 11:00 plane
    \end{itemize}
    Using the \rname{truth-functional} interpretation, \ref{ex:f1}
    asserts that the following conditional is true:
    \begin{equation}
      \tcond{A}{(\tcond{\neg A}{B})}
    \end{equation}
    Using the \rname{Boolean-feasibility} interpretation we consider
    the corresponding conditional:
    \begin{equation}
      \bcond{A}{(\bcond{\neg A}{B})}
    \end{equation}
    Unlike the familiar \rname{truth-functional} conditional, the
    \rname{Boolean-feasibility} conditional defined by
    Equation~\ref{eq:boolean-feasibility} is false when its antecedent
    is false.  So, using the \rname{Boolean-feasibility}
    interpretation, the inner conditional $\tcond{\neg A}{B}$ must be
    false subject to the constraint $A=\tval$ imposed by the outer
    conditional.  It follows that the outer conditional is false too
    (true antecedent, false consequent).

    We can arrive at the same conclusion regarding the
    \rname{Boolean-feasibility} interpretation of \ref{ex:f1} using
    the probability and linear-programming approach from
    Section~\ref{sec:prob-boole}.  This has the advantage of offering
    some explanation as to why the conditional stated as \ref{ex:f1}
    is incorrect.  Following Equations \ref{eq:animplb} and
    \ref{eq:aandb} we regard the inner conditional $\tcond{\neg A}{B}$
    as equivalent to the following pair of constraints concerning a
    parametric probability model including the primary variables $A$
    and $B$:
    \begin{eqnarray}
      \prob{\embed{\neg A \wedge \neg B}=\tval} & = & 0
      \\
      \prob{\embed{\neg A \wedge B}=\tval} & > & 0 
    \end{eqnarray}
    The antecedent $A$ of the outer conditional signifies the
    following singleton set of constraints $\Gamma$:
    \begin{eqnarray}
      \Gamma & \Leftarrow & \{ \: \prob{A=\tval}=1 \:\}
    \end{eqnarray}
    In order to determine whether the outer conditional
    $\bcond{A}{(\bcond{\neg A}{B})}$ holds we evaluate the following
    solution sets $\Phi$ and $\Psi$:
    \begin{eqnarray}
      \Phi & \Leftarrow & \{ \:
      \prob{\embed{\neg A \wedge \neg B}=\tval} \: : \:
      \Gamma \: \}
      \\
      \Psi & \Leftarrow & \{ \:
      \prob{\embed{\neg A \wedge B}=\tval} \: : \:
      \prob{\embed{\neg A \wedge \neg B}=\tval} = 0, \:
      \Gamma \: \}
    \end{eqnarray}
    Using the probability model from Section~\ref{sec:prob-boole}
    supplemented by additional embedded propositional-calculus
    functions, and including the general probability constraints
    $\Gamma_0=\{0 \le x_i \le 1, \sum_i x_i = 1\}$, symbolic
    probability inference yields:
    \begin{eqnarray}
      \Phi & \Rightarrow & \left\{ \:
x_{4}
      \: : \:
x_{1} + x_{2}
      = 1; \:
      0 \le x_i \le 1, \: \textstyle\sum_i x_i = 1
      \: \right\}
      \\
      \Psi & \Rightarrow & \left\{ \:
x_{3}
      \: : \:
x_{4}
      = 0, \:
x_{1} + x_{2}
      = 1; \:
      0 \le x_i \le 1, \: \textstyle\sum_i x_i = 1
      \: \right\}
    \end{eqnarray}
    Linear optimization computes $\Phi \Rightarrow \{0\}$ and $\Psi
    \Rightarrow \{0\}$.  It follows that the
    \rname{Boolean-feasibility} conditional
    \begin{math}
      \bcond{A}{(\bcond{\neg A}{B})}
    \end{math}
    cannot hold.  These solution sets provide some
    explanation for this impossibility.  Subject to the constraint
    $\Gamma$ representing the antecedent $A$ of the outer conditional,
    the result $\Phi=\{0\}$ shows that the \emph{negative} requirement
    of the inner conditional is always satisfied: it is infeasible for
    the consequent $B$ and the antecedent $\neg A$ to be true
    simultaneously.  However the result $\Psi=\{0\}$ indicates that
    the \emph{positive} requirement of the inner conditional cannot be
    satisfied: there must be zero probability that the consequent $B$
    is true and the antecedent $\neg A$ is true.  Therefore the inner
    conditional $\bcond{\neg A}{B}$ is necessarily false, subject to
    the constraint that the antecedent $A$ of the outer conditional is
    true.
    In other words the set-comprehension expression
    \begin{math}
      \{ \: (\bcond{\neg A}{B}) \: : \: A=\tval \: \}
    \end{math}
    evaluates to $\{\fval\}$.  Because this solution set is not
    $\{\tval\}$, Equation~\ref{eq:boolean-feasibility} tells us that
    the outer conditional
    \begin{math}
      \bcond{A}{(\bcond{\neg A}{B})}
    \end{math}
    is false.

    For the solution set $\Phi$, linear optimization gives solutions
    \begin{math}
      \alpha^* =
0.000
    \end{math}
    and
    \begin{math}
      \beta^* =
0.000
    \end{math}
    to the optimization problems:
    \begin{equation}
\begin{array}{r@{\quad}l}
\mbox{Minimize}: & x_{4} \\
\mbox{subject to}:
& x_{1} + x_{2} + x_{3} + x_{4} = 1 \\
& x_{1} + x_{2} = 1 \\
\mbox{and}:

& 0 \le x_{1} \le 1 \\
& 0 \le x_{2} \le 1 \\
& 0 \le x_{3} \le 1 \\
& 0 \le x_{4} \le 1 \\

\end{array}
      \qquad
\begin{array}{r@{\quad}l}
\mbox{Maximize}: & x_{4} \\
\mbox{subject to}:
& x_{1} + x_{2} + x_{3} + x_{4} = 1 \\
& x_{1} + x_{2} = 1 \\
\mbox{and}:

& 0 \le x_{1} \le 1 \\
& 0 \le x_{2} \le 1 \\
& 0 \le x_{3} \le 1 \\
& 0 \le x_{4} \le 1 \\

\end{array}
    \end{equation}

    For the solution set $\Psi$, linear optimization gives solutions
    \begin{math}
      \alpha^* =
0.000
    \end{math}
    and
    \begin{math}
      \beta^* =
0.000
    \end{math}
    to the optimization problems:
    \begin{equation}
\begin{array}{r@{\quad}l}
\mbox{Minimize}: & x_{3} \\
\mbox{subject to}:
& x_{1} + x_{2} + x_{3} + x_{4} = 1 \\
& x_{4} = 0 \\
& x_{1} + x_{2} = 1 \\
\mbox{and}:

& 0 \le x_{1} \le 1 \\
& 0 \le x_{2} \le 1 \\
& 0 \le x_{3} \le 1 \\
& 0 \le x_{4} \le 1 \\

\end{array}
      \qquad
\begin{array}{r@{\quad}l}
\mbox{Maximize}: & x_{3} \\
\mbox{subject to}:
& x_{1} + x_{2} + x_{3} + x_{4} = 1 \\
& x_{4} = 0 \\
& x_{1} + x_{2} = 1 \\
\mbox{and}:

& 0 \le x_{1} \le 1 \\
& 0 \le x_{2} \le 1 \\
& 0 \le x_{3} \le 1 \\
& 0 \le x_{4} \le 1 \\

\end{array}
    \end{equation}

    Note these inputs and outputs:
    \begin{equation}
\begin{tabular}[c]{ll|l}\hline
\multicolumn{1}{l}{$A$} & 
\multicolumn{1}{l|}{$B$} & 
\multicolumn{1}{l}{$\prob[0]{{A, B}}$} \\ \hline\hline
$\state{T}$ & 
$\state{T}$ & 
$x_{1}$ \\ \hline
$\state{T}$ & 
$\state{F}$ & 
$x_{2}$ \\ \hline
$\state{F}$ & 
$\state{T}$ & 
$x_{3}$ \\ \hline
$\state{F}$ & 
$\state{F}$ & 
$x_{4}$ \\ \hline
\end{tabular}
      \qquad
\begin{tabular}[c]{l|l}\hline
\multicolumn{1}{l|}{$\embed{\neg A \wedge \neg B}$} & 
\multicolumn{1}{l}{$\prob{{\embed{\neg A \wedge \neg B}}}$} \\ \hline\hline
$\state{T}$ & 
$x_{4}$ \\ \hline
$\state{F}$ & 
$x_{1} + x_{2} + x_{3}$ \\ \hline
\end{tabular}
      \qquad
\begin{tabular}[c]{l|l}\hline
\multicolumn{1}{l|}{$\embed{\neg A \wedge B}$} & 
\multicolumn{1}{l}{$\prob{{\embed{\neg A \wedge B}}}$} \\ \hline\hline
$\state{T}$ & 
$x_{3}$ \\ \hline
$\state{F}$ & 
$x_{1} + x_{2} + x_{4}$ \\ \hline
\end{tabular}
    \end{equation}

    Alternatively, in order to evaluate the inner conditional
    $\bcond{\neg A}{B}$ subject to the constraint $A=\tval$ from the
    outer conditional, we could evaluate the following solution set:
    \begin{eqnarray}
      \Upsilon & \Leftarrow &
      \{\:
      \condprob{B=\tval}{A=\fval} \: : \: \prob{A=\tval}=1
      \:\}
    \end{eqnarray}
    These are the resulting fractional linear optimization problems:
    \begin{equation}
\begin{array}{r@{\quad}l}
\mbox{Minimize}: & \left(x_{3}\right) / \left(x_{3} + x_{4}\right) \\
\mbox{subject to}:
& x_{1} + x_{2} + x_{3} + x_{4} = 1 \\
& x_{1} + x_{2} = 1 \\
\mbox{and}:

& 0 \le x_{1} \le 1 \\
& 0 \le x_{2} \le 1 \\
& 0 \le x_{3} \le 1 \\
& 0 \le x_{4} \le 1 \\

\end{array}
      \qquad
\begin{array}{r@{\quad}l}
\mbox{Maximize}: & \left(x_{3}\right) / \left(x_{3} + x_{4}\right) \\
\mbox{subject to}:
& x_{1} + x_{2} + x_{3} + x_{4} = 1 \\
& x_{1} + x_{2} = 1 \\
\mbox{and}:

& 0 \le x_{1} \le 1 \\
& 0 \le x_{2} \le 1 \\
& 0 \le x_{3} \le 1 \\
& 0 \le x_{4} \le 1 \\

\end{array}
    \end{equation}
    The author's solver \cite{norman-nlp} returns the status
\texttt{infeasible}
    for both problems, because the denominator of their common
    objective function has been constrained to equal zero.  Hence the
    solution set $\Upsilon\Rightarrow\emptyset$ meaning that, subject
    to the constraint $A=\tval$ from the outer conditional, the inner
    conditional $\bcond{\neg A}{B}$ cannot be true.  It follows that
    the outer conditional is incorrect also.

  \item \label{ex:f2} \emph{John will arrive on the 10 o'clock plane.
    Therefore, if John misses his plane in New York, he will arrive on
    the 10 o'clock plane.}

    Updated variables:
    \begin{itemize}
      \item[$A:$] John arrives on the 10:00 plane
      \item[$B:$] John misses his plane in New York
    \end{itemize}

    The \rname{truth-functional} interpretation of \ref{ex:f2} is the
    correct assertion that the following propositional-calculus
    formula is true:
    \begin{equation}
      \tcond{A}{(\tcond{B}{A})}
    \end{equation}
    Let us evaluate the corresponding \rname{Boolean-feasibility}
    conditional:
    \begin{equation}
      \bcond{A}{(\bcond{B}{A})}
    \end{equation}
    Following Equations \ref{eq:animplb} and \ref{eq:aandb} we
    interpret the inner conditional $\bcond{B}{A}$ as this pair of
    constraints:
    \begin{eqnarray}
      \prob{\embed{B \wedge \neg A}=\tval} & = & 0
      \\
      \prob{\embed{B \wedge A}=\tval} & > & 0 
    \end{eqnarray}
    In order to evaluate the satisfiability of these constraints given
    the antecedent $A$ of the outer conditional, we consider these two
    solution sets:
    \begin{eqnarray}
      \Phi & \Leftarrow & \{ \:
      \prob{\embed{B \wedge \neg A}=\tval} \: : \:
      \prob{A=\tval}=1 \: \}
      \\
      \Psi & \Leftarrow & \{ \:
      \prob{\embed{B \wedge A}=\tval} \: : \:
      \prob{\embed{B \wedge \neg A}=\tval} = 0, \:
      \prob{A=\tval}=1 \: \}
    \end{eqnarray}
    Symbolic probability inference with the probability model from
    Section~\ref{sec:prob-boole} yields the algebraic expressions:
    \begin{eqnarray}
      \Phi & \Rightarrow & \left\{ \:
x_{3}
      \: : \:
x_{1} + x_{2}
      = 1; \:
      0 \le x_i \le 1, \: \textstyle\sum_i x_i = 1
      \: \right\}
      \\
      \Psi & \Rightarrow & \left\{ \:
x_{1}
      \: : \:
x_{3}
      = 0, \:
x_{1} + x_{2}
      = 1; \:
      0 \le x_i \le 1, \: \textstyle\sum_i x_i = 1
      \: \right\}
    \end{eqnarray}
    Note these inputs and outputs:
    \begin{equation}
\begin{tabular}[c]{ll|l}\hline
\multicolumn{1}{l}{$A$} & 
\multicolumn{1}{l|}{$B$} & 
\multicolumn{1}{l}{$\prob[0]{{A, B}}$} \\ \hline\hline
$\state{T}$ & 
$\state{T}$ & 
$x_{1}$ \\ \hline
$\state{T}$ & 
$\state{F}$ & 
$x_{2}$ \\ \hline
$\state{F}$ & 
$\state{T}$ & 
$x_{3}$ \\ \hline
$\state{F}$ & 
$\state{F}$ & 
$x_{4}$ \\ \hline
\end{tabular}
      \qquad
\begin{tabular}[c]{l|l}\hline
\multicolumn{1}{l|}{$\embed{B \wedge \neg A}$} & 
\multicolumn{1}{l}{$\prob{{\embed{B \wedge \neg A}}}$} \\ \hline\hline
$\state{T}$ & 
$x_{3}$ \\ \hline
$\state{F}$ & 
$x_{1} + x_{2} + x_{4}$ \\ \hline
\end{tabular}
      \qquad
\begin{tabular}[c]{l|l}\hline
\multicolumn{1}{l|}{$\embed{B \wedge A}$} & 
\multicolumn{1}{l}{$\prob{{\embed{B \wedge A}}}$} \\ \hline\hline
$\state{T}$ & 
$x_{1}$ \\ \hline
$\state{F}$ & 
$x_{2} + x_{3} + x_{4}$ \\ \hline
\end{tabular}
    \end{equation}

    For the solution set $\Phi$, linear optimization gives solutions
    \begin{math}
      \alpha^* =
0.000
    \end{math}
    and
    \begin{math}
      \beta^* =
0.000
    \end{math}
    to the optimization problems:
    \begin{equation}
\begin{array}{r@{\quad}l}
\mbox{Minimize}: & x_{3} \\
\mbox{subject to}:
& x_{1} + x_{2} + x_{3} + x_{4} = 1 \\
& x_{1} + x_{2} = 1 \\
\mbox{and}:

& 0 \le x_{1} \le 1 \\
& 0 \le x_{2} \le 1 \\
& 0 \le x_{3} \le 1 \\
& 0 \le x_{4} \le 1 \\

\end{array}
      \qquad
\begin{array}{r@{\quad}l}
\mbox{Maximize}: & x_{3} \\
\mbox{subject to}:
& x_{1} + x_{2} + x_{3} + x_{4} = 1 \\
& x_{1} + x_{2} = 1 \\
\mbox{and}:

& 0 \le x_{1} \le 1 \\
& 0 \le x_{2} \le 1 \\
& 0 \le x_{3} \le 1 \\
& 0 \le x_{4} \le 1 \\

\end{array}
    \end{equation}

    For the solution set $\Psi$, linear optimization gives solutions
    \begin{math}
      \alpha^* =
0.000
    \end{math}
    and
    \begin{math}
      \beta^* =
1.000
    \end{math}
    to the optimization problems:
    \begin{equation}
\begin{array}{r@{\quad}l}
\mbox{Minimize}: & x_{1} \\
\mbox{subject to}:
& x_{1} + x_{2} + x_{3} + x_{4} = 1 \\
& x_{3} = 0 \\
& x_{1} + x_{2} = 1 \\
\mbox{and}:

& 0 \le x_{1} \le 1 \\
& 0 \le x_{2} \le 1 \\
& 0 \le x_{3} \le 1 \\
& 0 \le x_{4} \le 1 \\

\end{array}
      \qquad
\begin{array}{r@{\quad}l}
\mbox{Maximize}: & x_{1} \\
\mbox{subject to}:
& x_{1} + x_{2} + x_{3} + x_{4} = 1 \\
& x_{3} = 0 \\
& x_{1} + x_{2} = 1 \\
\mbox{and}:

& 0 \le x_{1} \le 1 \\
& 0 \le x_{2} \le 1 \\
& 0 \le x_{3} \le 1 \\
& 0 \le x_{4} \le 1 \\

\end{array}
    \end{equation}

    These optimization results reveal $\Phi \Rightarrow \{0\}$ and
    $\Psi \subseteq [0,1]$ with $0 \in \Psi$ and $1 \in \Psi$.
    According to the criteria of Section~\ref{sec:prob-boole} these
    results indicate that, assuming the outer antecedent $A$ is true,
    the inner conditional $\bcond{B}{A}$ is possibly true, but not
    necessarily true.  Because $\Psi$ does not exclude $0$ there is
    ambiguity about whether the positive requirement $\prob{\embed{B
        \wedge A}}=\tval$ is met.

    You can see by inspection of the objective function
    \begin{math}
x_{1}
    \end{math}
    of the linear programs used to compute $\Psi$ that the inner
    conditional $\bcond{B}{A}$ must be false (subject to the outer
    antecedent $A=\tval$) precisely when
    \begin{math}
x_{1}
      = 0
    \end{math},
    that is when the input
    \begin{math}
      \prob[0]{A=\tval,B=\tval} = 0
    \end{math}.
    Conversely if
    \begin{math}
      \prob[0]{A=\tval,B=\tval} > 0
    \end{math}
    (note the strict inequality) then the solution set $\Psi$ must
    exclude zero and hence the inner conditional is necessarily true.

    In terms of the example \ref{ex:f2}, what would the constraint
    $\prob[0]{A=\tval,B=\tval}=0$ mean?  Well, this assertion states
    that it is impossible \emph{a priori} for John to arrive on the
    10:00 plane and also to have missed his plane in New York.  This
    assertion would be quite correct if the plane from New York is the
    same one that is scheduled to arrive at 10 o'clock.

    We can analyze \ref{ex:f2} in a slightly different way using a
    modified truth table to compute which sets of valuations satisfy
    the \rname{Boolean-feasibility} conditionals involved.  Consider
    the following view of the inner conditional $\bcond{B}{A}$:
    \begin{equation}
      \begin{array}{cc|c} \hline
        A & B & \{A:B=\tval\} \\ \hline\hline
        \tval & \tval & \{\tval\} \\ \hline
        \tval & \fval & \emptyset \\ \hline
        \fval & \tval & \{\fval\} \\ \hline
        \fval & \fval & \emptyset \\ \hline
      \end{array}
    \end{equation}
    This indicates that $4$ of the $16$ possible sets of $(A,B)$
    valuations satisfy the requirements for $\bcond{B}{A}$:
    \begin{equation}
      \{ \: (\tval,\tval) \: \}, \qquad
      \{ \: (\tval,\tval), \: (\tval,\fval) \: \}, \qquad
      \{ \: (\tval,\tval), \: (\fval,\fval) \: \}, \qquad
      \{ \: (\tval,\tval), \: (\tval,\fval), \: (\fval,\fval) \: \}
    \end{equation}
    All $4$ sets include the valuation $(A,B)=(\tval,\tval)$.  On the
    other hand it is possible to satisfy the premise $A=\tval$ without
    including this valuation $(\tval,\tval)$.  These are the $3$
    valuations satisfying the premise $A=\tval$, equivalently
    $\bcond{\tval}{A}$ (meaning $\{A:\tval=\tval\}=\{\tval\}$):
    \begin{equation}
      \{ \: (\tval,\tval) \: \}, \qquad
      \{ \: (\tval,\fval) \: \}, \qquad
      \{ \: (\tval,\tval), \: (\tval,\fval) \: \}
    \end{equation}
    Note that the probability constraint $\prob{A=\tval}>0$ gives the
    algebraic constraint
    \begin{math}
x_{1} + x_{2}
      > 0
    \end{math}
    which is satisfied by either
    \begin{math}
x_{1}
      > 0
    \end{math}
    or
    \begin{math}
x_{2}
      > 0
    \end{math}
    or both: these strict inequalities indicate which valuations must
    be included in the sets of satisfactory valuations.

    The valuation set
    \begin{math}
      \{ \: (\tval,\fval) \: \}
    \end{math}
    satisfies the outer antecedent $A$ but not the outer consequent
    $\bcond{B}{A}$.  Therefore $\bcond{B}{A}$ is not a consequence of
    $A$.  In order to achieve consequentiality we would have to
    declare that the valuation $(A,B)=(\tval,\tval)$ must be possible:
    that it is feasible for John to miss his plane in New York and
    still arrive on the 10:00 plane (presumably not the same
    aircraft!).

    For yet another approach we can evaluate the inner conditional
    $\bcond{B}{A}$ subject to the constraint $A=\tval$ from the outer
    conditional using the following solution set:
    \begin{eqnarray}
      \Upsilon & \Leftarrow &
      \{\:
      \condprob{A=\tval}{B=\tval} \: : \: \prob{A=\tval}=1
      \:\}
    \end{eqnarray}
    These are the resulting fractional linear optimization problems:
    \begin{equation}
\begin{array}{r@{\quad}l}
\mbox{Minimize}: & \left(x_{1}\right) / \left(x_{1} + x_{3}\right) \\
\mbox{subject to}:
& x_{1} + x_{2} + x_{3} + x_{4} = 1 \\
& x_{1} + x_{2} = 1 \\
\mbox{and}:

& 0 \le x_{1} \le 1 \\
& 0 \le x_{2} \le 1 \\
& 0 \le x_{3} \le 1 \\
& 0 \le x_{4} \le 1 \\

\end{array}
      \qquad
\begin{array}{r@{\quad}l}
\mbox{Maximize}: & \left(x_{1}\right) / \left(x_{1} + x_{3}\right) \\
\mbox{subject to}:
& x_{1} + x_{2} + x_{3} + x_{4} = 1 \\
& x_{1} + x_{2} = 1 \\
\mbox{and}:

& 0 \le x_{1} \le 1 \\
& 0 \le x_{2} \le 1 \\
& 0 \le x_{3} \le 1 \\
& 0 \le x_{4} \le 1 \\

\end{array}
    \end{equation}
    The computed solutions
    \begin{math}
      \alpha^* =
1.000
    \end{math}
    and
    \begin{math}
      \beta^* =
1.000
    \end{math}
    indicate $\Upsilon\Rightarrow\{1\}$.  However in this case it is
    feasible for the denominator
    \begin{math}
x_{1} + x_{3}
    \end{math}
    of the objective function to equal zero.  A separate optimization
    problem would confirm this.  Thus we amend the results to say that
    either $\Upsilon=\{1\}$ or $\Upsilon=\emptyset$.

  \item \label{ex:f3} \emph{If Brown wins the election, Smith will
    retire to private life.  Therefore, if Smith dies before the
    election and Brown wins it, Smith will retire to private life.}

    Variables:
    \begin{itemize}
      \item[$A:$] Brown wins the election
      \item[$B:$] Smith retires to private life
      \item[$C:$] Smith dies before the election
    \end{itemize}
    The \rname{Truth-functional} interpretation of \ref{ex:f3} gives
    antecedent $\tcond{A}{B}$ and consequent $\tcond{(A \wedge
      C)}{B}$, for an integrated conditional statement:
    \begin{equation}
      \tcond{(\tcond{A}{B})}{(\tcond{(A \wedge C)}{B})}
    \end{equation}
    As a statement of material implication, this is tautologically
    true.  However the content of the example suggests that
    propositions $B$ and $C$ cannot be true simultaneously: Smith
    could not possibly retire to private life if he had already died
    before the election.  In this circumstance it is unintuitive to
    regard $B$ as a consequence of $A$ and $C$.

    Using the \rname{Boolean-feasibility} interpretation provides a
    different result.  Consider the conditional:
    \begin{equation}
      \bcond{(\bcond{A}{B})}{(\bcond{(A \wedge C)}{B})}
    \end{equation}
    Following Equations \ref{eq:animplb} and \ref{eq:aandb} we
    interpret the inner antecedent conditional $\bcond{A}{B}$ as the
    pair of probability constraints, designated as the set $\Gamma$:
    \begin{eqnarray}
      \prob{\embed{A \wedge \neg B}=\tval} & = & 0 \\
      \prob{\embed{A \wedge B}=\tval} & > & 0 
    \end{eqnarray}
    Likewise we interpret the inner consequent conditional $\bcond{(A
      \wedge C)}{B}$ as the pair of constraints:
    \begin{eqnarray}
      \prob{\embed{(A \wedge C) \wedge \neg B}=\tval} & = & 0 \\
      \prob{\embed{(A \wedge C) \wedge B}=\tval} & > & 0 
    \end{eqnarray}
    In order to compute whether the second pair of constraints follows
    from the first, we evaluate these two solution sets:
    \begin{eqnarray}
      \Phi & \Leftarrow & \{ \:
      \prob{\embed{(A \wedge C) \wedge \neg B}=\tval} \: : \:
      \Gamma
      \: \}
      \\
      \Psi & \Leftarrow & \{ \:
      \prob{\embed{(A \wedge C) \wedge B}=\tval} \: : \:
      \prob{\embed{(A \wedge C) \wedge \neg B}=\tval} = 0, \:
      \Gamma
      \: \}
    \end{eqnarray}
    Using the constraint set
    \begin{math}
      \Gamma \Leftarrow \{ \:
      \prob{\embed{A \wedge \neg B}=\tval} = 0, \:
      \prob{\embed{A \wedge B}=\tval} > 0 
      \: \}
    \end{math}.
    The probability model from Section~\ref{sec:transitivity} adds the
    general constraints
    \begin{math}
      \Gamma_0 \Leftarrow 
      \left\{ \: 0 \le y_i \le 1, \; \sum_i y_i = 1 \: \right\}
    \end{math}.

    For the solution set $\Phi$ linear optimization yields minimum and
    maximum feasible values
    \begin{math}
      \alpha^* =
0.000
    \end{math}
    and
    \begin{math}
      \beta^* =
0.000
    \end{math}
    indicating $\Phi\Rightarrow\{0\}$.
    For the solution set $\Psi$ linear optimization yields minimum and
    maximum feasible values
    \begin{math}
      \alpha^* =
0.000
    \end{math}
    and
    \begin{math}
      \beta^* =
1.000
    \end{math}
    indicating $\Psi \subseteq [0,1]$.  These results $\Phi=\{0\}$ and
    $\Psi\subseteq[0,1]$ indicate that
    \begin{math}
      \bcond{(A \wedge C)}{B}
    \end{math}
    is a possible but not necessary consequence of the premise
    \begin{math}
      \bcond{A}{B}
    \end{math}.

    Inspection of the linear programming problems used to calculate
    $\Psi$ helps to characterize when this deduction is correct and
    when it is not.  These are the problems:
    \begin{equation}
\begin{array}{r@{\quad}l}
\mbox{Minimize}: & y_{1} \\
\mbox{subject to}:
& y_{1} + y_{2} + y_{3} + y_{4} + y_{5} + y_{6} + y_{7} + y_{8} = 1 \\
& y_{3} = 0 \\
& y_{3} + y_{4} = 0 \\
& y_{1} + y_{2} \ge \epsilon \\
\mbox{and}:

& 0 \le y_{1} \le 1 \\
& 0 \le y_{2} \le 1 \\
& 0 \le y_{3} \le 1 \\
& 0 \le y_{4} \le 1 \\
& 0 \le y_{5} \le 1 \\
& 0 \le y_{6} \le 1 \\
& 0 \le y_{7} \le 1 \\
& 0 \le y_{8} \le 1 \\
& \epsilon = 0.001
\end{array}
      \qquad
\begin{array}{r@{\quad}l}
\mbox{Maximize}: & y_{1} \\
\mbox{subject to}:
& y_{1} + y_{2} + y_{3} + y_{4} + y_{5} + y_{6} + y_{7} + y_{8} = 1 \\
& y_{3} = 0 \\
& y_{3} + y_{4} = 0 \\
& y_{1} + y_{2} \ge \epsilon \\
\mbox{and}:

& 0 \le y_{1} \le 1 \\
& 0 \le y_{2} \le 1 \\
& 0 \le y_{3} \le 1 \\
& 0 \le y_{4} \le 1 \\
& 0 \le y_{5} \le 1 \\
& 0 \le y_{6} \le 1 \\
& 0 \le y_{7} \le 1 \\
& 0 \le y_{8} \le 1 \\
& \epsilon = 0.001
\end{array}
    \end{equation}
    You can see that the constraint
    \begin{math}
y_{1}
      = 0
    \end{math}
    would force the solutions to both optimization problems to zero,
    hence producing $\Psi=\{0\}$ and rendering the outer conditional
    false.  On the other hand the constraint
    \begin{math}
y_{1}
      > 0
    \end{math}
    (implemented as 
    \begin{math}
y_{1}
      \ge \epsilon
    \end{math}
    with a small numerical constant $\epsilon$ for the optimization
    solver) would force the minimum solution $\alpha^*$ to be strictly
    greater than zero, hence producing $0 \notin \Psi$ (yet still
    $\Psi\neq\emptyset$) and thereby rendering the outer conditional
    true.
    
    In terms of the example \ref{ex:f3} these results say the
    following.  It is correct to deduce the conclusion `If Smith dies
    before the election and Brown wins it, Smith will retire to
    private life' from the premise `If Brown wins the election, Smith
    will retire to private life'---\emph{unless} it is known \emph{a
      priori} that it would be impossible for all three events to
    occur together (Brown's victory, Smith's retirement, and Smith's
    death), in which case the deduction is incorrect.  If the prior
    probability of these three events is not known (meaning that it
    could be anywhere between zero and one, as limited by the laws of
    probability), then it is possible but not necessary that the
    stated conclusion is a consequence of the stated premise.  For the
    content of this example it would be appropriate to impose the
    constraint $\prob{B=\tval,C=\tval}=0$ to specify that it would be
    impossible for Smith to retire after having died.  This would
    yield the algebraic constraint
    \begin{math}
y_{1} + y_{5}
      = 0
    \end{math}
    which would force
    \begin{math}
y_{1}
      = 0
    \end{math}
    (since the laws of probability constrain
    \begin{math}
      0 \le
y_{1}
      \le 1
    \end{math}
    and
    \begin{math}
      0 \le
y_{5}
      \le 1
    \end{math}).
    Equivalently the valuations $(A,B,C)=(\tval,\tval,\tval)$ and
    $(\fval,\tval,\tval)$ could be declared to be impossible.

  \item \label{ex:f4} \emph{If Brown wins the election, Smith will
    retire to private life.  If Smith dies before the election, Brown
    will win it.  Therefore, if Smith dies before the election, then he
    will retire to private life.}

    We continue to use the variables from \ref{ex:f3}:
    \begin{itemize}
      \item[$A:$] Brown wins the election
      \item[$B:$] Smith retires to private life
      \item[$C:$] Smith dies before the election
    \end{itemize}
    The \rname{truth-functional} interpretation asserts the following
    formula of the propositional calculus:
    \begin{equation}
      \tcond{((\tcond{A}{B})\wedge(\tcond{C}{A}))}{(\tcond{C}{B})}
      \label{eq:f4-truth-functional}
    \end{equation}
    And indeed this formula is tautologically true.  For this problem
    \ref{ex:f3}, the corresponding \rname{Boolean-feasibility}
    conditional is also correct:
    \begin{equation}
      \bcond{((\bcond{A}{B})\wedge(\bcond{C}{A}))}{(\bcond{C}{B})}
      \label{eq:f4-boolean-feasibility}
    \end{equation}
    This compound conditional statement is simply a permutation of the
    transitivity example from Section~\ref{sec:transitivity}.  For
    this problem \ref{ex:f4} the issue is not the way the conditionals
    involved are interpreted and analyzed; the issue is the content of
    those conditionals.  In this case it is essential to revise the
    formal model as it is being developed, in order to account for
    information that is revealed during the course of the problem
    description.

    Let us examine the process of formulating conditional statements
    to model \ref{ex:f3}.  We invoke the idea of an agent called `the
    Analyst' which develops the formal model.  First, the Analyst
    splits the problem description \ref{ex:f3} into three
    English-language sentences, each containing a conditional
    statement:
    \begin{itemize}
      \item[$\mathcal{E}_1:$] If Brown wins the election, Smith will
        retire to private life.
      \item[$\mathcal{E}_2:$] If Smith dies before the election, Brown
        will win it.
      \item[$\mathcal{E}_3:$] If Smith dies before the election, then
        he will retire to private life.
    \end{itemize}
    The Analyst then sets out to model each English sentence as a
    formal conditional statement.  Since the type of conditionals is
    not important here, the familiar \rname{truth-functional}
    interpretation is used.  Upon reading the first sentence
    $\mathcal{E}_1$ the Analyst defines the variable $A$ to mean that
    Brown wins, and $B$ to mean that Smith retires.  The Agent then
    introduces the following conditional $\mathcal{S}_1$ (with
    antecedent $A$ and consequent $B$) to represent $\mathcal{E}_1$:
    \begin{itemize}
      \item[$\mathcal{S}_1:$] \tcond{A}{B}
    \end{itemize}
    Upon reading $\mathcal{E}_3$ the Agent adds the variable $C$ to
    mean that Smith dies.  The Agent then introduces the following
    conditionals $\mathcal{S}_2$ and $\mathcal{S}_3$ to represent the
    English sentences $\mathcal{E}_2$ and $\mathcal{E}_3$:
    \begin{itemize}
      \item[$\mathcal{S}_2:$] \tcond{C}{A}
      \item[$\mathcal{S}_3:$] \tcond{C}{B}
    \end{itemize}
    Assembling these three conditionals $\mathcal{S}_1$,
    $\mathcal{S}_2$, and $\mathcal{S}_3$ according to the connective
    `therefore' included in the original problem statement yields the
    compound conditional shown in
    Equation~\ref{eq:f4-truth-functional} which is tautologically true
    (though perhaps unintuitive).

    Upon reflection, the Agent now realizes that there is a problem
    with the initial conditional $\mathcal{S}_1$.  It was revealed in
    the second sentence $\mathcal{E}_2$ that Smith could die before
    the election.  The Agent has prior knowledge that a dead person
    could not possibly retire to private life.  Therefore the Agent
    revises the initial statement $\mathcal{S}_1$, which was
    $\tcond{A}{B}$, to account for the exceptional circumstance of
    Smith's death:
    \begin{itemize}
      \item[$\mathcal{S}'_1:$] \tcond{(A \wedge \neg C)}{B}
    \end{itemize}
    This revised conditional $\mathcal{S}'_1$ reflects an updated
    English sentence $\mathcal{E}'_1$, which might be stated as:
    \begin{itemize}
      \item[$\mathcal{E}'_1:$] If Brown wins the election, and Smith
        does not die before the election, then Smith will retire to
        private life.
      \item[$\mathcal{E}'_1:$] If Brown wins the election, then Smith
        will retire to private life---unless Smith dies before the
        election.
    \end{itemize}
    Assembling this revised conditional $\mathcal{S}'_1$ with
    $\mathcal{S}_2$ and $\mathcal{S}_3$ results in the following
    compound conditional:
    \begin{equation}
      \tcond{((\tcond{(A \wedge \neg C)}{B})\wedge(\tcond{C}{A}))}{(\tcond{C}{B})}
    \end{equation}
    Unlike the formula in Equation~\ref{eq:f4-truth-functional}, this
    revised statement of material implication is not tautologically
    true.  In fact, if $A$ and $C$ are true and $B$ is false (Brown
    won, Smith died and did not retire) then the antecedent
    \begin{math}
      (\tcond{(A \wedge \neg C)}{B})\wedge(\tcond{C}{A})
    \end{math}
    is true but the consequent
    \begin{math}
      \tcond{C}{B}
    \end{math}
    is false.
    
    The Agent, now worried about other exceptions, may wish to revise
    $\mathcal{S}_2$ and the understanding of $\mathcal{E}_2$ also to
    state that Brown could not win the election even after Smith had
    died, if some exceptional event $Z$ happens that
    renders Brown unable to win:
    \begin{itemize}
      \item[$\mathcal{S}'_2:$] $\tcond{(C \wedge \neg Z)}{A}$
      \item[$\mathcal{E}'_2:$] If Smith dies before the election, then
        Brown will win it---unless something exceptional happens.
    \end{itemize}
    For the content of this example it would make sense to add some
    new premises that say that Smith cannot retire after he has died
    (regardless of whether or not Brown won), and that the
    hypothetical exception $Z$ is indeed incompatible with Brown's
    victory:
    \begin{itemize}
      \item[$\mathcal{S}_4:$] $\neg(B \wedge C)$
      \item[$\mathcal{S}_5:$] $\neg(Z \wedge A)$
    \end{itemize}
    
    So there are at least two ways to handle potential exceptions to
    conditional statements: first, to revise preexisting conditional
    statements in order to account for specific exceptions, once those
    exception have been discovered (as in revising $\tcond{A}{B}$ into
    $\tcond{(A \wedge \neg C)}{B}$ to account for the exception $C$);
    second, to express conditional statements in a defeasible manner
    in the first place, in order to account for unknown exceptions (as
    in starting with $\tcond{(C \wedge \neg Z)}{A}$ instead of
    $\tcond{C}{A}$, anticipating that some exception $Z$ might later
    be discovered).  After including unknown exceptions such as $Z$ it
    may be appropriate to analyze each problem twice, once with the
    constraint that $Z$ is true and once with the constraint that $Z$
    is false, in order to see how the potential exception would affect
    the results of inference.
    
  \item \label{ex:f5} \emph{If Brown wins, Smith will retire.  If
    Brown wins, Smith will not retire.  Therefore, Brown will not
    win.}
    
    Interpreted using \rname{truth-functional} conditionals we have
    the tautologically-true formula:
    \begin{equation}
      \tcond{((\tcond{A}{B}) \wedge (\tcond{A}{\neg B}))}{\neg A}
    \end{equation}
    \rname{Boolean-feasibility} conditionals behave differently.
    We consider the corresponding formula:
    \begin{equation}
      \bcond{((\bcond{A}{B}) \wedge (\bcond{A}{\neg B}))}{\neg A}
      \label{eq:f5-boolean-feasibility}
    \end{equation}
    Following Equation~\ref{eq:boolean-feasibility} each inner
    conditional is defined by an equation about sets of truth values:
    \begin{eqnarray}
      \bcond{A}{B} & \equiv & \{\: B \: : \: A=\tval \:\} = \{\tval\} \\
      \bcond{A}{\neg B} & \equiv & \{\: \neg B \: : \: A=\tval \:\} = \{\tval\}
    \end{eqnarray}
    Applying the negation operator both sides of the second equation
    gives the revised form
    \begin{math}
      \{\: B \: : \: A=\tval \:\} = \{\fval\}
    \end{math}
    for the conditional $\bcond{A}{\neg B}$.  The same solution set
    cannot equal $\{\tval\}$ and $\{\fval\}$ simultaneously, hence the
    antecedent of the outer conditional in
    Equation~\ref{eq:f5-boolean-feasibility} is necessarily false.  A
    \rname{Boolean-feasibility} conditional with a false antecedent is
    false itself; hence the deduction described by \ref{ex:f5} is
    incorrect using the \rname{Boolean-feasibility} interpretation of
    its conditional statements.

  \item \label{ex:f6} \emph{Either Dr.~A or Dr.~B will attend the
    patient.  Dr.~B will not attend the patient.  Therefore, if Dr.~A
    does not attend the patient, Dr.~B will.}
    
    Here are the variables:
    \begin{itemize}
    \item[$A:$] Dr.~A attends the patient
    \item[$B:$] Dr.~B attends the patient
    \end{itemize}
    Using \rname{truth-functional} conditionals gives the
    tautologically-true formula:
    \begin{equation}
      \tcond{((A \vee B) \wedge \neg B)}{(\tcond{\neg A}{B})}
    \end{equation}
    Using the \rname{Boolean-feasibility} interpretation we consider
    the isomorphic formula:
    \begin{equation}
      \bcond{((A \vee B) \wedge \neg B)}{(\bcond{\neg A}{B})}
      \label{eq:f6-boolean-feasibility}
    \end{equation}
    This new conditional statement is necessarily false because a
    \rname{Boolean-feasibility} conditional with a false premise is
    not considered true; in this case the truth of the outer
    antecedent
    \begin{math}
      (A \vee B) \wedge \neg B
    \end{math}
    requires the falsity of the inner antecedent $\neg A$.  One way to
    confirm this result is with probability and linear programming.
    Following Equations \ref{eq:animplb} and \ref{eq:aandb} the inner
    conditional $\bcond{\neg A}{B}$ requires that these two
    constraints are both satisfied:
    \begin{eqnarray}
      \prob{\embed{\neg A \wedge \neg B}=\tval} & = & 0 \\
      \prob{\embed{\neg A \wedge B}=\tval} & > & 0
    \end{eqnarray}
    Thus in order to test the conditional in
    Equation~\ref{eq:f6-boolean-feasibility} we evaluate these two
    solution sets:
    \begin{eqnarray}
      \Phi & \Leftarrow & \{\:
      \prob{\embed{\neg A \wedge \neg B}=\tval} \: : \:
      \Gamma
      \:\} \\
      \Psi & \Leftarrow & \{\:
      \prob{\embed{\neg A \wedge B}=\tval} \: : \:
      \prob{\embed{\neg A \wedge \neg B}=\tval} = 0, \:
      \Gamma
      \:\}
    \end{eqnarray}
    with
    \begin{math}
      \Gamma \Leftarrow \{\:
      \prob{\embed{(A \vee B) \wedge \neg B}=\tval} = 1
      \:\}
    \end{math}
    representing the antecedent of the outer conditional.  

    Using the probability model from Section~\ref{sec:prob-boole},
    linear programming gives solutions
    \begin{math}
      \alpha^* =
0.000
    \end{math}
    and
    \begin{math}
      \beta^* =
0.000
    \end{math}
    for $\Phi$ and
    \begin{math}
      \alpha^* =
0.000
    \end{math}
    and
    \begin{math}
      \beta^* =
0.000
    \end{math}
    for $\Psi$, indicating both $\Phi\Rightarrow\{0\}$ and
    $\Psi\Rightarrow\{0\}$.  This pattern means that, when the outer
    antecedent
    \begin{math}
      (A \vee B) \wedge \neg B
    \end{math}
    from Equation~\ref{eq:f6-boolean-feasibility} is true, the inner
    conditional $\bcond{\neg A}{B}$ must be false (because it is not
    feasible that the inner antecedent and consequent are
    simultaneously true, as the inner conditional requires).
    Therefore the outer conditional in
    Equation~\ref{eq:f6-boolean-feasibility} must always be false.

  \item \label{ex:f7} \emph{It is not the case that if John passes
    history, he will graduate.  Therefore John will pass history.}

    The variables are designated as follows:
    \begin{itemize}
    \item[$A:$] John passes history
    \item[$B:$] John graduates
    \end{itemize}
    Using \rname{truth-functional} conditionals gives this
    tautologically true formula of the propositional calculus:
    \begin{eqnarray}
      \tcond{\neg(\tcond{A}{B})}{A}
    \end{eqnarray}
    We consider the isomorphic statement using
    \rname{Boolean-feasibility} conditionals:
    \begin{eqnarray}
      \bcond{\neg(\bcond{A}{B})}{A}
      \label{eq:f7-boolean-feasibility}
    \end{eqnarray}
    The insight here is that it is possible to invalidate the inner
    conditional $\bcond{A}{B}$ without requiring the truth of
    proposition $A$.  Recall from
    Equation~\ref{eq:boolean-feasibility} that the conditional
    $\bcond{A}{B}$ is defined by the constraint that truth is the only
    feasible value for $B$, subject to the constraint that $A$ is
    true:
    \begin{eqnarray}
      \bcond{A}{B} & \equiv & \{\: B \::\: A=\tval \:\} = \{\tval\}
    \end{eqnarray}
    There are three different ways to achieve the negation
    \begin{math}
      \{\: B \::\: A=\tval \:\} \neq \{\tval\}
    \end{math}
    of this equation about sets of truth values: the solution set
    could equal $\{\tval,\fval\}$, or $\{\fval\}$, or $\emptyset$.  In
    neither of these three cases is it required that truth is the one
    and only feasible value for $A$.

    We can evaluate the possible truth values of the antecedent and
    consequent of the outer conditional in
    Equation~\ref{eq:f7-boolean-feasibility}, considering various sets
    of valuations of the variables $A$ and $B$ to be possible.  Let us
    begin with the assumption that both valuations
    $(A,B)=(\fval,\tval)$ and $(\fval,\fval)$ are possible.  If no
    other valuations are possible this gives the set
    \begin{math}
      \mathbf{V}_1 \Leftarrow
      \{\: (\fval,\tval), \: (\fval,\fval) \:\}
    \end{math}
    of possible valuation vectors.  Considering these two possible
    valuations: the set
    \begin{math}
      \{\: B \::\: A=\tval \:\}
    \end{math}
    evaluates to the empty set, and $A$ has possible values
    $\{\fval\}$.  Because 
    \begin{math}
      \{\: B \::\: A=\tval \:\} = \emptyset
    \end{math}
    the \rname{Boolean-feasibility} conditional $\bcond{A}{B}$ defined
    by Equation~\ref{eq:boolean-feasibility} is false in the
    circumstance that the set $\mathbf{V}_1$ gives the possible
    valuations of $(A,B)$.  In terms of the example \ref{ex:f7}, if it
    is known \emph{a priori} that John will not pass history, then the
    \rname{Boolean-feasibility} interpretation of `If John passes
    history, he will graduate' must be false whereas the proposition
    `John will pass history' must be false; hence the overall
    conditional statement \ref{ex:f7} is incorrect.

    Similarly, assuming the set
    \begin{math}
      \mathbf{V}_2 \Leftarrow
      \{\: (\tval,\fval), \: (\fval,\tval), \: (\fval,\fval) \:\}
    \end{math}
    of possible valuations for $(A,B)$ yields the results
    \begin{math}
      \{\: B \::\: A=\tval \:\} \Rightarrow \{\fval\}
    \end{math}
    and
    \begin{math}
      \{\: A \:\} \Rightarrow \{\tval,\fval\}
    \end{math}.
    In this case $\bcond{A}{B}$ must be false, hence the outer
    antecedent $\neg(\bcond{A}{B})$ in
    Equation~\ref{eq:f7-boolean-feasibility} must be true.  Yet the
    outer consequent $A$ may be either true or false.  In other words
    the solution set
    \begin{math}
      \{\: A \: : \: \neg(\bcond{A}{B})=\{\tval\} \:\}
      \Rightarrow \{\tval,\fval\}
    \end{math},
    which by Definition~\ref{eq:boolean-feasibility} means that the
    \rname{Boolean-feasibility} conditional in
    Equation~\ref{eq:f7-boolean-feasibility} is not correct.

    Finally, assuming the set
    \begin{math}
      \mathbf{V}_3 \Leftarrow
      \{\: (\tval,\tval), \: (\tval,\fval), \: 
      (\fval,\tval), \: (\fval,\fval) \:\}
    \end{math}
    of possible valuations for $(A,B)$ yields the results
    \begin{math}
      \{\: B \::\: A=\tval \:\} \Rightarrow \{\tval,\fval\}
    \end{math}
    and
    \begin{math}
      \{\: A \:\} \Rightarrow \{\tval,\fval\}
    \end{math}.
    Again $\bcond{A}{B}$ is necessarily false while $A$ could be
    either true or false; the conditional in
    Equation~\ref{eq:f7-boolean-feasibility} fails.

    In terms of the content of this example \ref{ex:f7}, these results
    from the valuation sets $\mathbf{V}_2$ and $\mathbf{V}_3$ say
    that, if it is possible that John could not pass history, and
    furthermore possible that John could pass history and not
    graduate, then it is not correct to deduce from the premise `It is
    not the case that if John passes history, he will graduate' the
    conclusion that `John will pass history'.  Whereas the premise is
    true given either set of possibilities $\mathbf{V}_2$ or
    $\mathbf{V}_3$, in both of these cases the conclusion could be
    either true or false.
    
  \item \label{ex:f8} \emph{If you throw both switch S and switch T,
    the motor will start.  Therefore, either if you throw switch S the
    motor will start, or if you throw switch T the motor will start.}

    These are the variables:
    \begin{itemize}
    \item[$A:$] You throw switch S
    \item[$B:$] You throw switch T
    \item[$C:$] The motor starts
    \end{itemize}
    Using the \rname{truth-functional} interpretation of conditionals
    gives this statement of material implication, which is
    tautologically true:
    \begin{eqnarray}
      \tcond{(\tcond{(A \wedge B)}{C})}{(\tcond{A}{C})\vee(\tcond{B}{C})}
    \end{eqnarray}
    We consider the isomorphic \rname{Boolean-feasibility} formula:
    \begin{eqnarray}
      \bcond{(\bcond{(A \wedge B)}{C})}{(\bcond{A}{C})\vee(\bcond{B}{C})}
      \label{eq:f8-boolean-feasibility}
    \end{eqnarray}
    Each conditional statement included in
    Equation~\ref{eq:f8-boolean-feasibility} constrains the set of
    possible valuations of the variables $(A,B,C)$.  Here we will use
    a table-based approach to computing and describing the relevant
    sets of valuations.  

    \begin{table}
      \begin{math}
        \begin{array}{ccc|cccc} \\ \hline
          A & B & C & A \wedge B &
          \{C:(A \wedge B)=\tval\} &
          \{C:A=\tval\} &
          \{C:B=\tval\}
          \\ \hline\hline
          \tval & \tval & \tval & \tval & \{\tval\} & \{\tval\} & \{\tval\}
          \\ \hline
          \tval & \tval & \fval & \tval & \{\fval\} & \{\fval\} & \{\fval\}
          \\ \hline
          \tval & \fval & \tval & \fval & \emptyset & \{\tval\} & \emptyset
          \\ \hline
          \tval & \fval & \fval & \fval & \emptyset & \{\fval\} & \emptyset
          \\ \hline
          \fval & \tval & \tval & \fval & \emptyset & \emptyset & \{\tval\}
          \\ \hline
          \fval & \tval & \fval & \fval & \emptyset & \emptyset & \{\fval\}
          \\ \hline
          \fval & \fval & \tval & \fval & \emptyset & \emptyset & \emptyset
          \\ \hline
          \fval & \fval & \fval & \fval & \emptyset & \emptyset & \emptyset
          \\ \hline
        \end{array}
      \end{math}
      \caption{Valuation worksheet for \ref{ex:f8}, indicating which
        sets of valuations satisfy each inner conditional statement
        in Equation~\ref{eq:f8-boolean-feasibility} as described in
        the text.}
      \label{tbl:f8-worksheet}
    \end{table}

    Table~\ref{tbl:f8-worksheet} shows a worksheet that describes the
    sets of valuations that describe each inner conditional that is
    included in Equation~\ref{eq:f8-boolean-feasibility}.  We begin
    with the conditional
    \begin{math}
      \bcond{(A \wedge B)}{C}
    \end{math}
    which is defined by Equation~\ref{eq:boolean-feasibility} as the
    following equation about sets of truth values:
    \begin{eqnarray}
      \{\: C \: : \: (A \wedge B)=\tval \:\} & = & \{\tval\}
      \label{eq:f8-cond1}
    \end{eqnarray}
    Considering binary variables $A$, $B$, and $C$ there are $2^{2^3}$
    or $256$ possible sets of $(A,B,C)$ valuation vectors; let us use
    $\mathbf{V}^*$ to designate this set of valuations.  The relevant
    column in Table~\ref{tbl:f8-worksheet} includes the values of the
    solution set 
    \begin{math}
      \{\: C \: : \: (A \wedge B)=\tval \:\}
    \end{math}
    assuming that only the valuation indicated by each row is
    possible.  The tabulated results describe the features of the set
    $\mathbf{V}_1 \subseteq \mathbf{V}^*$ of valuations that satisfy
    Equation~\ref{eq:f8-cond1}: each valuation set must contain
    $(A,B,C)=(\tval,\tval,\tval)$; it must omit
    $(A,B,C)=(\tval,\tval,\fval)$; and it may or may not contain any
    of the other valuations.  There are $2^6$ or $64$ such sets of
    valuations in the set $\mathbf{V}_1$, for example the sets
    \begin{math}
      \{\: (\tval,\tval,\tval) \:\}
    \end{math}
    and
    \begin{math}
      \{\: (\tval,\tval,\tval), \: (\fval,\fval,\fval) \:\}
    \end{math}.

    Regarding the conditional $\bcond{A}{C}$,
    Equation~\ref{eq:boolean-feasibility} gives the defining equation
    \begin{eqnarray}
      \{\: C \: : \: A=\tval \:\} & = & \{\tval\}
      \label{eq:f8-cond2}
    \end{eqnarray}
    Table~\ref{tbl:f8-worksheet} shows that each member of the set
    $\mathbf{V}_2$ of valuations satisfying Equation~\ref{eq:f8-cond2}
    must contain either $(A,B,C)=(\tval,\tval,\tval)$ or
    $(\tval,\fval,\tval)$ or both; it must omit $(\tval,\tval,\fval)$
    and $(\tval,\fval,\fval)$; and it may or may not contain any of
    the other valuations.  There are $48$ such sets of valuations.

    Finally for the conditional $\bcond{B}{C}$,
    Equation~\ref{eq:boolean-feasibility} gives the defining equation
    \begin{eqnarray}
      \{\: C \: : \: B=\tval \:\} & = & \{\tval\}
      \label{eq:f8-cond3}
    \end{eqnarray}
    Table~\ref{tbl:f8-worksheet} shows that each member of the set
    $\mathbf{V}_3$ of valuation-sets satisfying
    Equation~\ref{eq:f8-cond3} must contain either
    $(A,B,C)=(\tval,\tval,\tval)$ or $(\fval,\tval,\tval)$ or both; it
    must omit $(\tval,\tval,\fval)$ and $(\fval,\tval,\fval)$; and it
    may or may not contain any of the other valuations.  There are
    $48$ such sets of valuations.

    Integrating these results for the valuation sets $\mathbf{V}_1$,
    $\mathbf{V}_2$, and $\mathbf{V}_3$, it is clear how to construct a
    set of valuations that satisfies the conditional
    \begin{math}
      \bcond{(A \wedge B)}{C}
    \end{math}
    from the antecedent of Equation~\ref{eq:f8-boolean-feasibility}
    but neither of the conditionals $\bcond{A}{C}$ nor $\bcond{B}{C}$
    from the consequent.  Such a valuation-set must include
    $(A,B,C)=(\tval,\tval,\tval)$, $(\tval,\fval,\fval)$, and
    $(\fval,\tval,\fval)$; it must omit $(\tval,\tval,\fval)$; and it
    may or may not contain the other valuations.  There are $16$ such
    sets of valuations, including for example:
    \begin{eqnarray}
      (A,B,C) & \in &
      \{\: 
      (\tval,\tval,\tval), \: (\tval,\fval,\fval), \: (\fval,\tval,\fval)
      \:\}
    \end{eqnarray}
    In terms of the content of this example \ref{ex:f8}, there are
    many conceivable circumstances in which the antecedent of the
    outer conditional in Equation~\ref{eq:f8-boolean-feasibility}
    holds but the consequent does not; in other words using the
    \rname{Boolean-conditional} interpretation, the stated deduction
    is incorrect.  In each of the counterexamples it is possible to
    throw switch S alone without starting the motor (the valuation
    $(A,B,C)=(\tval,\fval,\fval)$ is possible) and it is also possible
    to throw switch T alone without starting the motor
    ($(A,B,C)=(\fval,\tval,\fval)$ is possible); yet if both switches
    are thrown the motor must start ($(A,B,C)=(\tval,\tval,\tval)$ is
    possible and $(\tval,\tval,\fval)$ is impossible).

    \item \label{ex:f9} \emph{If John will graduate only if he passes
      history, then he won't graduate.  Therefore, if John passes
      history he won't graduate.}

      We return to the variables from \ref{ex:f7}:
      \begin{itemize}
      \item[$A:$] John passes history
      \item[$B:$] John graduates
      \end{itemize}
      Let us interpret `$B$ only if $A$' to mean `It is not the case
      that $B$ is true and $A$ is false'.  Using the
      \rname{truth-functional} interpretation of conditionals, this
      problem \ref{ex:f9} asserts the following formula of the
      propositional calculus:
      \begin{equation}
        \tcond{(\tcond{\neg(B \wedge \neg A)}{\neg B})}{(\tcond{A}{\neg B})}
      \end{equation}
      Indeed the above formula is tautologically true.  We consider
      the isomorphic formula using \rname{Boolean-feasibility}
      conditionals:
      \begin{equation}
        \bcond{(\bcond{\neg(B \wedge \neg A)}{\neg B})}{(\bcond{A}{\neg B})}
        \label{eq:f9-boolean-feasibility}
      \end{equation}
      
      \begin{table}
        \begin{math}
          \begin{array}{cc|cccc} \\ \hline
            A & B & \neg B & \neg(B \wedge \neg A) &
            \{\neg B : \neg(B \wedge \neg A)=\tval\} &
            \{\neg B : A=\tval\}
            \\ \hline\hline
            \tval & \tval & \fval & \tval & \{\fval\} & \{\fval\}
            \\ \hline
            \tval & \fval & \tval & \tval & \{\tval\} & \{\tval\}
            \\ \hline
            \fval & \tval & \fval & \fval & \emptyset & \emptyset
            \\ \hline
            \fval & \fval & \tval & \tval & \{\tval\} & \emptyset
            \\ \hline
          \end{array}
        \end{math}
        \caption{Valuation worksheet for \ref{ex:f9}, indicating which
          sets of valuations satisfy the inner conditionals in
          Equation~\ref{eq:f9-boolean-feasibility}.}
        \label{tbl:f9-worksheet}
      \end{table}

      The worksheet in Table~\ref{tbl:f9-worksheet} shows that the
      first inner conditional
      \begin{math}
        \bcond{\neg(B \wedge \neg A)}{\neg B}        
      \end{math}
      requires that either valuation $(A,B)=(\tval,\fval)$ or
      $(\fval,\fval)$ or both must be possible, and that
      $(A,B)=(\tval,\fval)$ must be impossible, in order to satisfy
      the equation
      \begin{math}
        \{ \: \neg B \: : \: \neg(B \wedge \neg A)=\tval \: \}
      \end{math}
      that defines this particular \rname{Boolean-feasibility}
      conditional according to Equation~\ref{eq:boolean-feasibility}.
      There are $6$ such sets of valuations.
      Also, Table~\ref{tbl:f9-worksheet} shows that the second inner
      conditional $\bcond{A}{\neg B}$ requires that the valuation
      $(A,B)=(\tval,\tval)$ must be possible and that
      $(A,B)=(\tval,\fval)$ must be impossible.  There are $4$ such
      sets of valuations.  There are altogether $16$ possible sets of
      valuations for the binary variables $(A,B)$.

      There two sets of possible valuations for $A$ and $B$ that
      satisfy the antecedent inner conditional in
      Equation~\ref{eq:f9-boolean-feasibility} but not the consequent.
      The requisite condition is that $(A,B)=(\fval,\fval)$ is
      possible, while both $(A,B)=(\tval,\tval)$ and $(\tval,\fval)$
      are impossible.  The remaining valuation $(A,B)=(\fval,\tval)$
      may be possible or impossible.  Here are the two matching sets
      of valuations:
      \begin{eqnarray}
        \mathbf{V}_1 & \Leftarrow & \{\: (\fval,\fval) \:\} \\
        \mathbf{V}_2 & \Leftarrow & \{\: (\fval,\tval), \: (\fval,\fval) \:\}
      \end{eqnarray}
      In terms of the content of \ref{ex:f9}, these results say that
      it is incorrect to deduce the conclusion `If John passes history
      he won't graduate' from the premise `If John will graduate only
      if he passes history, then he won't graduate'.  The reason is
      that in the case that it is impossible for John to pass history
      (both valuations $(A,B)=(\tval,\tval)$ and $(\tval,\fval)$ are
      impossible), and also possible for him not to graduate (the
      valuation $(A,B)=(\fval,\fval)$ is possible), then the stated
      premise would be satisfied yet the stated conclusion would be
      necessarily false (because the inner antecedent $A$ itself would
      be necessarily false).

      Analysis would have produced the same ultimate result, had we
      interpreted `$B$ only if $A$' to mean `$B$ if and only if $A$',
      as in the propositional-calculus formula $A \leftrightarrow B$
      (with the biconditional operation).

\end{enumerate}

\clearpage
\section{Discussion}

\subsection{Respecting Diversity Among Conditional Statements}

There are many different types of conditional statements.  The
distinction between subjunctive and indicative conditionals is
important.  But then among the indicatives there still several
distinct types.  Besides the principal types there are many potential
specializations of conditional statements, for example to account for
factual evidence that has been provided.  Anyway with appropriate use
of probability we can express each semantically distinct type and
subtype of conditional with a syntactically distinct algebraic
expression.

There are yet other kinds of conditional statements.  Recurrence
relations are an important omission from this document as they are out
of scope.  Indeed there are other mathematical means of deduction too,
which generate different (data) types of solutions (arithmetical,
algebraic, probabilistic, dynamical).

\subsection{Curating the Contents of Conditionals}

Independently of the type of each conditional statement, it is
important to get its content right.  That is, to ensure that the
formal model expresses the information that the modeler actually
intends.  It is a subtle but terrible fallacy to use certain pieces of
information during the informal phase of analysis, but then discard
that very information during the formal phase of analysis.  If you
think in English that it would be impossible for someone to arrive on
an airplane that he had failed to board in the first place, then you
should say that in formal language.  If you discover some exceptional
event that you had not previously considered, then you should revise
your earlier formal model to account for that event.  Or, if you
insist on never retracting any formal statement, then you should be
make those statements in a cautious and limited way that allows for
future exceptions.

\subsection{Polylogicism}

Appreciating the diversity of conditional statements and using the
methods of algebra, probability, and optimization to compute
deductions provides a new perspective on mathematical logic.  Let us
call this `polylogicism': an updated view of the close relationship
between logic and algebra, and the essentially mathematical nature of
logic.

\begin{raggedright}
  \small
  \bibliographystyle{plain}
  \addcontentsline{toc}{section}{\refname}
  \bibliography{../decisions}
\end{raggedright}

\clearpage
\appendix

\section{Optimization with Fractional Objectives}
\label{sec:fractional}

We consider also polynomial optimization problems with fractional
objectives.  Let us update the problem template from
Equation~\ref{eq:pp} by modifying the objective function to be a
quotient of polynomials:
\begin{equation}
  \begin{array}{rl}
    \mbox{Maximize}: & f(\mathbf{x}) / h(\mathbf{x})
    \\
    \mbox{subject to}: &
    g_1(\mathbf{x}) \ge 0, \ldots,
    g_m(\mathbf{x}) \ge 0
    \\
    \mbox{and}: &
    \alpha_1 \le x_1 \le \beta_1, \ldots,
    \alpha_n \le x_n \le \beta_n
  \end{array}
  \label{eq:fpp}
\end{equation}
Here the functions $f$, $h$, and each $g_j$ are polynomials in
$\Re[\mathbf{x}]$.  As before each variable $x_i$ is bounded by finite
limits $\alpha_i$ and $\beta_i$.  The author's solver accommodates
polynomial optimization problems with fractional objectives by
applying the Charnes-Cooper transformation \cite{charnes} to the
reformulated and linearized problems.  

The possibility of division by zero requires special attention.  Using
the Charnes-Cooper transformation, it is considered infeasible for the
denominator of the objective function to have the exact value zero.
However it can be detected when the value of a fractional objective is
unbounded as its denominator approaches zero.  These behaviors are
perhaps best explained by a few examples, in which the numerical
results returned by the author's solver have been annotated with the
status computed for each optimization problem:
\begin{eqnarray}
\begin{array}{r@{\quad}l}
\mbox{Minimize}: & y / x \\
\mbox{subject to}:
& 0 \le x \le 1 \\
& 0 \le y \le 1 \\

\end{array}
  & \Rightarrow & [
0
  ,
0
  ] \quad
\mathtt{optimal}
  \\
\begin{array}{r@{\quad}l}
\mbox{Maximize}: & y / x \\
\mbox{subject to}:
& 0 \le x \le 1 \\
& 0 \le y \le 1 \\

\end{array}
  & \Rightarrow &
\mathtt{unbounded}
  \\
\begin{array}{r@{\quad}l}
\mbox{Minimize}: & y / x \\
\mbox{subject to}:
& x = 0 \\
& 0 \le y \le 1 \\

\end{array}
  & \Rightarrow &
\mathtt{infeasible}
  \\
\begin{array}{r@{\quad}l}
\mbox{Minimize}: & x y / x \\
\mbox{subject to}:
& 0 \le x \le 1 \\
& 0 \le y \le 1 \\

\end{array}
  & \Rightarrow & [
0
  ,
0
  ] \quad
\mathtt{optimal}
  \\
\begin{array}{r@{\quad}l}
\mbox{Maximize}: & x y / x \\
\mbox{subject to}:
& 0 \le x \le 1 \\
& 0 \le y \le 1 \\

\end{array}
  & \Rightarrow & [
1
  ,
1
  ] \quad
\mathtt{optimal}
  \\
\begin{array}{r@{\quad}l}
\mbox{Minimize}: & x y / x \\
\mbox{subject to}:
& x = 0 \\
& 0 \le y \le 1 \\

\end{array}
  & \Rightarrow &
\mathtt{infeasible}
\end{eqnarray}
This behavior regarding division by zero requires two practical
considerations.  First, if a problem with a fractional objective has
been found to be infeasible, then it may be necessary to solve a
second problem with the same constraints and an arbitrary
non-fractional objective in order to distinguish whether the original
infeasibility occurred because the denominator of the fractional
objective must be zero, or because the other constraints in the
problem are infeasible.  Second, if it is important to establish
whether some particular indeterminate form such as $0/0$ or $1/0$ is a
feasible value for a fractional objective
$f(\mathbf{x})/h(\mathbf{x})$, then a separate problem must be solved
in order to determine this (for example, using the numerator
$f(\mathbf{x})$ alone as a non-fractional objective and adding the
constraint $h(\mathbf{x})=0$).

\clearpage

\begin{footnotesize}
  \tableofcontents
\end{footnotesize}

\end{document}